\documentclass[12pt]{amsart}
\usepackage{hyperref}
\usepackage{tikz-cd}
\numberwithin{equation}{section}
\begin{document}

\title{Tower multitype and global regularity of 
the $\bar\partial$-Neumann operator}
\author{Dmitri Zaitsev
}
\dedicatory{Dedicated to the memory of Professor Joe Kohn}
\subjclass[2020]{32W05, 32T25, 32V15, 32V35, 32F10, 35N15}
\maketitle
\tableofcontents

\def\Label#1{\label{#1}}

\def\1#1{\overline{#1}}
\def\2#1{\widetilde{#1}}
\def\3#1{\widehat{#1}}
\def\4#1{\mathbb{#1}}
\def\5#1{\frak{#1}}
\def\6#1{{\mathcal{#1}}}
\def\p#1{{\em #1}}

\def\C{{\4C}}
\def\R{{\4R}}
\def\N{{\4N}}
\def\Z{{\4Z}}
\def\K{{\4K}}

% Standard sets

\def\cn{{\C^n}}
\def\cnn{{\C^{n'}}}
\def\ocn{\2{\C^n}}
\def\ocnn{\2{\C^{n'}}}

% Abbreviations

\def\dist{{\rm dist}}
\def\const{{\rm const}}
\def\rk{{\rm rank\,}}
\def\ord{{\rm ord\,}}
\def\id{{\sf id}}
\def\aut{{\sf aut}}
\def\Aut{{\sf Aut}}
\def\CR{{\rm CR}}
\def\GL{{\sf GL}}
\def\Re{{\sf Re}}
\def\Im{{\sf Im}}
\def\span{\text{\rm span}}
\def\Arg{{\sf Arg}}
\def\Res{{\sf Res}}
\def\Ann{{\sf Ann}}

\def\codim{{\rm codim}}
\def\crd{\dim_{{\rm CR}}}
\def\crc{{\rm codim_{CR}}}

\def\phi{\varphi}
\def\eps{\varepsilon}
\def\d{\partial}
\def\a{\alpha}
\def\b{\beta}
\def\g{\gamma}
\def\G{\Gamma}
\def\D{\Delta}
\def\de{\delta}
\def\Om{\Omega}
\def\om{\omega}
\def\k{\kappa}
\def\l{\lambda}
\def\L{\Lambda}
\def\z{{\bar z}}
\def\w{{\bar w}}
\def\t{\tau}
\def\th{\theta}
\def\p{\phi}
\def\de{\delta}
\def\la{\langle}
\def\ra{\rangle}
\def\r{\rho}
\def\s{\sigma}

\emergencystretch15pt
\frenchspacing

\newtheorem{Thm}{Theorem}[section]
\newtheorem{Cor}[Thm]{Corollary}
\newtheorem{Con}[Thm]{Conjecture}
\newtheorem{Pro}[Thm]{Proposition}
\newtheorem{Lem}[Thm]{Lemma}

\theoremstyle{definition}\newtheorem{Def}[Thm]{Definition}

\theoremstyle{remark}
\newtheorem{Rem}[Thm]{Remark}
\newtheorem{Rems}[Thm]{Remarks}
\newtheorem{Exa}[Thm]{Example}
\newtheorem{Exs}[Thm]{Examples}

\def\bl{\begin{Lem}}
\def\el{\end{Lem}}
\def\bp{\begin{Pro}}
\def\ep{\end{Pro}}
\def\bt{\begin{Thm}}
\def\et{\end{Thm}}
\def\bc{\begin{Cor}}
\def\ec{\end{Cor}}
\def\bcn{\begin{Con}}
\def\ecn{\end{Con}}
\def\bd{\begin{Def}}
\def\ed{\end{Def}}
\def\br{\begin{Rem}}
\def\er{\end{Rem}}
\def\brs{\begin{Rems}}
\def\ers{\end{Rems}}
\def\be{\begin{Exa}}
\def\ee{\end{Exa}}
\def\bpf{\begin{proof}}
\def\epf{\end{proof}}
\def\bca{\begin{cases}}
\def\eca{\end{cases}}
\def\ben{\begin{enumerate}}
\def\een{\end{enumerate}}
\def\beq{\begin{equation}}
\def\eeq{\end{equation}}

\def\bex{\begin{Exa}}
\def\eex{\end{Exa}}
\def\bexs{\begin{Exs}}
\def\eexs{\end{Exs}}

\def\then{\Longrightarrow}

\def\bal{\begin{align*}}
\def\eal{{\end{align*}}}

\begin{abstract}
A new approach is given to property $(P_q)$ defined by Catlin for $q=1$ in a global and by Sibony in a local context, subsequently extended by Fu-Straube for $q>1$. This property is known to imply compactness and global regularity in the $\bar\partial$-Neumann problem as well as condition $R$ by Bell-Ligocka. In particular, we provide a self-contained proof of property $(P_q)$ for pseudoconvex hypersurfaces of finite D'Angelo $q$-type, the case originally studied by Catlin. Moreover, our proof covers more general classes of hypersurfaces inspired by a recent work of Huang-Yin. Proofs are broken down into isolated steps, some of which do not require pseudoconvexity.

  Our tools include: a new multitype invariant based on distinguished nested sequences of $(1,0)$ subbundles, defined in terms of derivatives of the Levi form; real and complex formal orbits; $k$-jets of functions relative to pairs of formal submanifolds; relative contact orders generalizing the usual contact orders; a new notion of supertangent vector fields having higher than expected relative contact orders; and a formal variant of a result by Diederich-Forn\ae ss arising as a key step in their proof of Kohn's ideal termination in the real-analytic case.
%A new approach
%is given
% to property $(P_q)$ 
% defined
%by Catlin for $q=1$ in a global
%and
%by
%Sibony in a local context,
%subsequently extended by
%Fu-Straube for $q>1$.
%This property
% is known to imply  compactness
% and
%  global regularity in the $\bar\partial$-Neumann problem
%  by a result of Kohn-Nirenberg,
%as well as condition $R$ by Bell-Ligocka.
%In particular, we provide
%a self-contained proof of property $(P_q)$ 
%for pseudoconvex hypersurfaces of finite D'Angelo $q$-type,
%the case originally studied by Catlin.
%Moreover, our proof covers
% more general classes of hypersurfaces
%inspired by a recent work of Huang-Yin.
%Proofs are broken down into isolated steps,
%some of which do not require pseudoconvexity.
%
%Our tools include: a new multitype invariant
%based on distinguished nested sequences of $(1,0)$ subbundles,
%defined in terms of derivatives
%of the Levi form;
% real and complex formal orbits;
% $k$-jets of functions relative to pairs of formal submanifolds;
% relative contact orders 
% generalizing the usual contact orders;
%a new notion of supertangent vector fields
% having higher than expected relative contact orders;
% and a formal variant
% of a
% result by
%  Diederich-Forn\ae ss
%  arising as 
%a key step in their proof 
%of Kohn's ideal termination
%in the real-analytic 
%case.
\end{abstract}

\section{Introduction}
\subsection{Origin of the paper}
\Label{intro}
%A fundamental property of
%a boundary value problem for 
%a system of Partial Differential Equations
%is {\em global regularity} stating that
%weak solutions are smooth whenever 
%the data of the problem are smooth.
%Global regularity of the
% {\em $\bar\d$-Neumann
%problem},
%due to its implication 
%via Kohn's formula
%\cite{K84}
%of
%{\em condition $R$}~\cite{BeL80, Be81},
%is
% of further importance
%for the study of smooth boundary extension
%of proper holomorphic maps,
%generalizing the famous Fefferman's theorem~\cite{Fe74},
%and other geometric applications, 
%see e.g.~ \cite{Hu93}.
A fundamental question 
for
a boundary value problem for 
a system of Partial Differential Equations
is finding conditions for
{\em global regularity},
i.e.\ for
(weak) solutions being smooth whenever 
the data of the problem are smooth.
Global regularity of the
 {\em $\bar\d$-Neumann
problem},
due to its implication 
via Kohn's formula
\cite{K84}
of
{\em condition $R$}~\cite{BeL80, Be81},
is
 of further importance
for the study of smooth boundary extension
of proper holomorphic maps,
generalizing the famous Fefferman's theorem~\cite{Fe74},
and other geometric applications, 
see e.g.~ \cite{Hu93}.
For background material on the $\bar\d$-Neumann problem,
that was formulated in the fifties by Spencer as a means
to generalize Hodge theory to non-compact
complex manifolds,
the reader is referred 
to \cite{FK72,D93,BS99,CS01,FS01,Za08,St10,Has14},
in particular,
to
the excellent monograph by Straube \cite{St10}
for specific aspects of the theory used here.

In his influential work,
Catlin 
\cite{C84proc,C84ann}
developed highly elaborate
geometric and potential-theoretic
techniques 
to tackle the problem of
global regularity
of the $\bar\d$-Neumann problem
for
 {\em pseudoconvex bounded
 domains in $\cn$ of finite D'Angelo type}
 \cite{D82}.
More specifically,
Catlin's property $(P)$
inspired vast
applications 
(see \S\ref{pq})
and
continuing research
on
{\em compactness} in the $\bar\d$-Neumann problem,
that is known to imply global regularity
by the celebrated work of Kohn-Nirenberg~\cite{KN65}.

In a sharp contrast,
there seems to be a lack 
in understanding
the orginal work~\cite{C84ann},
which
still,  after 40 years, remains
the only available source establishing 
%\cite[Theorem~3]{C84proc}
%on
%{\em Catlin's multitype},
%needed to obtain
property $(P)$
(and hence
compactness, global regularity of the $\bar\d$-Neumann operator
 and condition $R$)
for general
smoothly
bounded
{\em  pseudoconvex domains of finite type in any dimension}.
%remains Catlin's original work \cite{C84ann}.
Despite its significance
and numerous requests by experts in the field
(apart from some few more special situations
 \cite{Z17, BNZ18,A24}),
 there had not been
any alternative exposition or simplification
of the techniques 
%of \cite{C84ann}
covering the general case.

\subsection{Some new techniques introduced here}
\Label{new-tech}
One of the present paper's goals
is to fill this gap and provide
a new approach 
using a {\em different multitype invariant},
called the {\em tower multitype},
and a different set of new techniques.
In comparison with approach in \cite{C84ann},
our tower multitype is
based on
minimization of numbers of vector fields
defining certain
{\em distinguished nested sequences of $(1,0)$ subbundles},
called {\em towers},
rather than maximization
of rational weights
over choices of local holomorphic coordinates.
We don't use weighted
truncations
of geometric objects and their 
coordinate representations
(cf.~\cite[\S4-\S10]{C84ann}).
Instead,
further new techniques
include:
\ben
%\item
%{\em special subbundles}
%defined in terms of derivatives of the Levi form
%along vector fields;
\item
 {\em formal (Nagano type) orbits} $O$
and their {\em complexifications} $V$
for certain 
{\em special subbundles}
defined in terms of derivatives
of the Levi form,
that arise as obstruction
to the finiteness of the tower multitype;
\item
{\em formal CR property}
for orbits $O$ as above, that
we call the {\em formal Huang-Yin property}
in analogy with its convergent variant
crucial in their work \cite{HY21};
a key consequence is the
 regularity of the complexifications $V$ of $O$;
 \item
 {\em $k$-jets relative to a pair $(O,V)$,
 or briefly $(k,O,V)$-jets},
 where $O,V$ are as above
 or, more generally, any pair of formal submanifolds;
 \item
 {\em relative contact orders} 
 of real hypersurfaces with pairs $(O,V)$
 of formal submanifolds,
 generalizing usual 
{\em contact orders of with complex submanifolds};
% (the definition extends to arbitrary pairs
% of singular formal varieties or
% pairs of ideals in rings of formal power series);
% \item
% proving that relative contact orders 
% with tangent pairs
% are always {\em even} for
% pseudoconvex hypersurfaces,
% generalizing the well-known property
% of the usual contact orders with
% complex submanifolds;
 \item
 {\em supertangent vector fields}
 having higher than expected relative contact orders,
 and
 {\em complex-supertangent vector fields}
 $L$, such that both $L$ and $JL$
 are supertangent,
for which we establish an important
 {\em Lie algebra property}
 (here $J$ is the complex structure);
% for
% the {\em complex-supertangent} 
% vector fields $L$
% such that both $L$ and $\1L$ are supertangent.
\item
infiniteness of relative contact orders
with $(O,V)$ 
when $O$ is {\em complex-tangential},
(i.e.\ tangent to the complex distribution of the hypersurface);
this may be regarded as
a formal variant of 
the important result by 
Diederich-Forn\ae ss
\cite[Proposition~3]{DiF78}
---
a difficult key step in their proof 
of Kohn's ideal termination
in the real-analytic 
case, also
Siu 
%for a geometric proof
~
\cite[p~1170 and Part IV]{S10}.
\een

%
%
%It can be instead compared to 
%Catlin's elaborate
%{\em commutator multitype}
%that is associated with boundary systems
%and choices of lists of vector fields.
%still based on maximization of weights.
%i
%\cite[\S2]{C84ann}
%which still 

\subsection{New finiteness conditions}
\Label{new-fin}
Another goal is this paper
is to provide a range of 
{\em new finiteness
conditions}
 implying global regularity.
Our inspiration here
comes from the recent
breakthrough work by Huang-Yin~\cite{HY21}
on the Bloom conjecture
comparing 3 different approaches
to (regular) types
going back to Kohn's fundamental paper \cite{K72}
--- respectively based on:
\ben
\item
contact orders with complex submanifolds;
\item
iterated commutators of sections in subbundles;
\item orders of nonzero derivatives of the Levi form along subbundles.
\een

See the survey \cite{HY23} for a thorough 
discussion and comparison of 
respective regular types 
defined via these approaches
and for the recent progress on the
 Bloom and D'Angelo conjectures.

Using each of the $3$ approaches,
we formulate $3$ finiteness conditions
on a real hypersurface and show
that each condition leads to global regularity.
Simpler special cases of each of these conditions
 are formulated in Theorem~\ref{reg},
while their refined more general versions in 
Theorem~\ref{refined}.
The difference is that
in Theorem~\ref{reg},
 {\em orders of contact, commutators and derivatives}
 are calculated for arbitrary
 submanifolds or subbundles,
 while in Theorem~\ref{refined},
 {\em only special submanifolds and subbundles}
 are used, that are
 defined in terms of the geometry of the Levi form.
 
Both Theorems~\ref{reg} and \ref{refined} will be reduced
to Theorem~\ref{main}
asserting the same conclusion
under the assumption
of the {\em $q$-finite tower multitype}
(Definition~\ref{q-finite-def}).
This reduction is obtained via Theorem~\ref{reduce}
asserting that assumptions in
either of Theorems~\ref{reg} or \ref{refined}
imply $q$-finite tower multitype.
We also provide
Example~\ref{t-type},
where the tower multitype is finite,
while other finiteness conditions fail to hold.

On the other hand, by the celebrated irregularity result
of Christ~\cite{Ch96} using the work of Barrett~\cite{Ba92},
global regularity does not hold
for general pseudoconvex 
smoothly bounded domains,
hence some conditions are needed.

More broadly,
providing conditions implying
global regularity in the 
$\bar\d$-Neumann problem
is the subject
of an active 
and vast research area,
for which we refer
to the sources
mentioned in \S\ref{intro}
and references therein.

\subsection{Closer inspection of the role of pseudoconvexity}
Because of the wider interest in global regularity
for {\em non-pseudoconvex domains},
especially for forms of higher degree $q>1$,
we
break down the entire proof into
isolated pieces
and localise those needing pseudoconvexity,
which happens only in certain crucial steps.
Such steps are:
\ben
\item
application of Sibony's $B$-regularity theory
\cite{Si87}
and its generalization for $q>1$
($B_q$-regularity)
 by Fu-Straube~\cite{FS01}
in Proposition~\ref{m-p};
\item
the
complex-tangential property 
for real orbits of subbundles of infinite Levi type
in Corollary~\ref{orb-tan};
\item
the infiniteness of the contact order
with complexifications $V$ of complex-tangential
formal submanifolds in Theorem~\ref{df},
whose proof uses pseudoconvexity 
twice:
\ben
\item
 in Corollary~\ref{must-even}
to conclude that the 
{\em relative contact order $k$ is even};
\item
in Corollary~\ref{every-super}
to conclude
that all vector fields tangent to $V$
are {\em supertangent}.
\een
\een

The exact parts of the proof of Theorem~\ref{df}
using pseudoconvexity
are those where Lemma~\ref{pos} is applied,
which
concludes that respective relative contact orders are even.

In particular, 
Proposition~\ref{reg-conclude}
on the existence of
generalized stratifications
and the part of Theorem~\ref{reduce}
using the Levi type conditions (3) and (3$'$)
in Theorems~\ref{reg} and \ref{refined}
do not assume pseudoconvexity.
Consequently,
extending 
$B_q$-regularity theory
to non-pseudoconvex domains
will provide 
global regularity
under conditions (3) and (3$'$).

\subsection{Path from finiteness conditions to global regularity}
\Label{path}
The scheme of implications
starting with geometric finiteness conditions
and ending with global regularity and condition R
used here
can be illustrated as follows:
$$
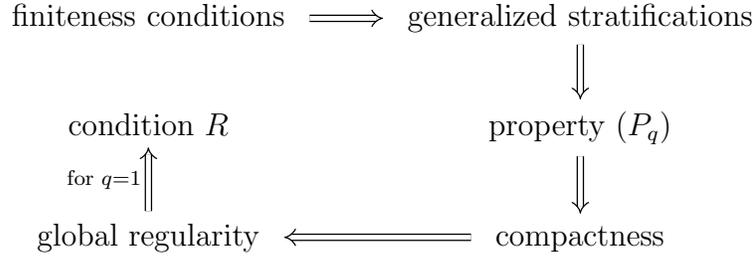
\begin{tikzcd}
	\text{ finiteness conditions } \rar[Rightarrow]
	%A \Rightarrow
	& \text{ generalized stratifications } \dar[Rightarrow]
	\\
	\text{ condition $R$ } & \text{ property $(P_q)$ } 
	\dar[Rightarrow]
	\\
	\text{ global regularity }\uar[Rightarrow]{\text{for }q=1}
	\rar[Leftarrow]
	& \text{ compactness }
\end{tikzcd}
$$

All our finiteness conditions will be reduced 
to ($q$-)finiteness of the tower multitype
using Theorem~\ref{reduce}.
In turn, $q$-finiteness of the tower multitype
implies a countable generalized stratification property
by means of Proposition~\ref{reg-conclude}.
Using $B_q$-regularity theory
\cite{Si87,FS01}, such generalized stratifications
lead to property $(P_q)$
(see Definition~\ref{p}), 
as formulated in Proposition~\ref{m-p}.

Other implications in the above diagram
starting from property $(P_q)$
are well-understood.
See
\cite{C84proc,St97} and
\cite[Theorem~4.8]{St10}
for property $(P_q)$ implying compactness,
\cite{KN65} and \cite[Theorem~4.6]{St10}
for compactness implying global regularity,
and \cite{K84}, \cite[(2.71)]{St10}, 
for global regularity implying {\em condition R}.
See
the monograph
by Straube \cite{St10}
and other sources mentioned above
for more detailed and thorough discussions.

\subsection{Property $(P_q)$}
\Label{pq}
A key step towards global regularity in \cite{C84proc}
is establishing 
%Catlin's 
%celebrated 
the
potential-theoretic 
{\em property~$(P)$}.
It has been extended to arbitrary compact sets 
by Sibony~\cite{Si87} 
and to $q>1$ by Fu-Straube~\cite{FS01}
as respectively 
{\em $B$-} and {\em $B_q$-regularity}
%see also McNeal~\cite{M02} for a more general variant.
(see also 
\cite{M02} for a variant of this condition and
 \cite{AFL21} for a related notion of {\em $p$-convexity} in Riemannian geometry):

\bd
\Label{p}
(Following \cite[\S4.4]{St10}.)
A compact subset $K\subset \cn$ is said to 
satisfy 
{\em property $(P_q)$}, where $1\le q\le n$ is an integer,
%or to {\em satisfy P} for brevity,
 if for every $C>0$, there exists an open neighborhood $U$ of $K$ in $\cn$
 and a real $C^2$ function $\l\colon U\to\R$ such that $0\le\l\le 1$
 and the sum of any $q$ eigenvalues of the complex hessian matrix
 $ (\l_{z_j\1z_k})_{j,k=1}^n$ is greater than $C$.
% $$
% 0\le\l\le 1,
% \quad 
% \sum_{j,k=1}^n \frac{\d^2\l}{\d z_j \d\bar z_k}(z) a_j \bar a_k
% \ge C \|a\|^2,
% \quad
% z\in U,
% \quad
% a=(a_1,\ldots,a_n)\in\cn.
%% \d\bar\d\l(L,\bar L) \ge M\|L\|^2.
% $$
\ed

Property $(P_q)$ has found vast applications, including:
{\em compactness} of the
{\em $\bar\d$-Neumann operator $N_q$} 
(i.e.\ the inverse of the
{\em $\bar\d$-Neumann Laplacian} $\square=\bar\d\bar\d^*+\bar\d^*\bar\d$ for $(0,q)$ forms)
\cite{C84proc,St97,ChF05,Ha11,Has11,AS15,Zi21},
\cite[Theorem~4.8]{St10},
implying 
{\em global and exact regularity}~\cite{KN65}, \cite[Theorem~4.6]{St10}, and subsequently,
regularity of the 
{\em Bergman projection $P_{q-1}$ of $(0,q-1)$-forms}  via Kohn's formula 
$P_{q-1}=\id - \bar\d^*N_q\bar\d$ \cite[(2.71)]{St10}, 
in particular {\em condition $R$} mentioned above;
regularity of the 
{\em Szeg\"o projection} \cite{B87},
estimates for the 
{\em Bergman kernels} \cite{BB81}, 
the {\em Szeg\"o kernels}
 and for the Bergman, Kobayashi metrics \cite{CF11,CF12};
compactness of the 
{\em complex Green operator} 
(i.e.\ the inverse of the Kohn Laplacian $\square_b=\bar\d_b\bar\d_b^*+\bar\d_b^*\bar\d_b$) \cite{N06,RS08,R10,KhPZ12,St12,BiS16};
the 
{\em Diederich-Fornaess} exponent and non-existence of Stein domains with Levi-flat boundaries \cite{Si87, FSh16}; 
%{\em division problems} for holomorphic functions with $C^\infty$ boundary values \cite{BM87};
{\em approximation by plurisubharmonic}
and by {\em holomorphic functions} in neighborhoods of domain's closure,
constructions of {\em bounded exhaustion functions}~\cite{Si87};
existence of {\em Stein neighborhood bases}~\cite{Si87};
construction of {\em plurisubharmonic peak and bumping functions}~\cite{Si87,BhS09,Bh12};
Holomorphic Morse inequalities and eigenvalue asymptotics for
the $\bar\d$-Neumann Laplacian~\cite{FJ10}.

Property $(P_q)$ is even proved to be equivalent to 
compactness of the $\bar\d$-Neuman operator by 
Fu-Straube~\cite{FS98,FS01}
for {\em localy convexifiable} domains
and by Christ-Fu~\cite{ChF05}
 for 
{\em complete Hartogs domains} in $\C^2$,
while in general, it is still unknown whether property $(P_q)$ is 
actually equivalent to the compactness.

\subsection{Three finiteness conditions}
We now give precise definitions
of the finiteness conditions
mentioned briefly before in \S\ref{new-fin}
in their simplest form,
while their more general and refined
form appear in Theorem~\ref{refined} below.

As a matter of notation,
when saying ``smooth'' we always mean the $C^\infty$ category.
Let $S\subset\C^n$ be a smooth real hypersurface
with
a {\em local defining function} $r$ at a point $p\in S$,
%in an open subset $U\subset\C^n$ 
i.e.\ a smooth real function
in a neighborhood $U$ of $p$
with $S\cap U=\{ r=0\}$ 
and $dr\ne 0$ at every point of $U$.
It is an easy consequence of the implicit function theorem
that any two local defining functions 
%$r_1,r_2\colon U\to \R$ of $S$
differ by a nonzero smooth function factor.
%in a neighborhood of $p$.
% For a smooth function $r$ in a neighborhood of $p$ in $\C^n$,
We write $\C TS$ for the complexified tangent bundle of $S$, 
$HS\subset TS$ and $\C HS\subset \C TS$
for the complex tangent bundle and its complexification,
and
 $H^{10}S, H^{01}S=\1{H^{10}S}\subset \C TS$ for the smooth
 complex subbundles
 of the tangent $(1,0)$ and $(0,1)$ vectors respectively,
 so that
  $\C HS = H^{10} S \oplus \1{H^{10} S}$.
% When the choice of $S$ is clear, we write
% $H^{10}=$H^{10}S

%A {\em formal holomorphic map}
%$\g
% = \sum_\a c_\a z^\a \
%\colon (\C^q,0)\to (\cn,p),
%$
%given by formal power series
%$$
%g=(g_1,\ldots,g_n)\in (\C[[t]])^n, \quad t\in \C^q,
%$$
%is any formal power series 
%
%We write $j^\infty_p r(z,\z)$ for the $\infty$-jet of $r$ at $p$ 
% represented by its formal Taylor series in $\C[[z-p,\1{z-p}]]$.
% $$
% j^\infty_p r
% = \sum_{\a\b\in \N^m}  
%\frac1{\a!\b!} \frac{\d^{|\a|+|\b|} f}{\d z^\a\d \1z^\b}(p)
% (z-p)^\a(\1{z-p})^\b
%\in \C[[z-p,\1{z-p}]].
%$$
%A {\em formal complex $q$-dimensional submanifold}
%$V\subset $
%is 

%The classes of domains for which we establish property $(P_q)$ are defined by means of three fundamental tools 
%
% 
We begin
with special cases of our
 finiteness conditions
implying property $(P_q)$
 that are the easiest to formulate
and will
provide their more refined versions
in Theorem~\ref{refined}
 below.

%For a smooth real hypersurface $S\subset \C^n$,
%
%$$
%t(E,p)
%$$
%there exist $m\ge2$ vector fields 
%$L^m,\ldots,L^1\in E\oplus \1B$ with

%\section{Finiteness conditions}

\bt[Special case of Theorem~\ref{refined}]
\Label{reg}
Let $D\subset\cn$, $n\ge2$, be a pseudoconvex bounded domain 
with smooth 
$(C^\infty)$ boundary $S=\d D$,
and $1\le q\le n-1$ be an integer. Assume that
 for every boundary point $p\in S$,
% with a local defining function $r$ at $p$,  
at least one of the following conditions holds:
\ben
\item 
there does not exist a 
%the contact order of $S$ with every 
formal holomorphic immersion 
$$
\g
% = \sum_\a c_\a z^\a \
\colon (\C^q,0)\to (\cn,p)
$$
that is tangent to $S$ of infinite order, 
i.e. the composition of the Taylor series of
a defining function of $S$ at $p$ with $\g$
vanishes identically;
%i.e.\ $j_p^\infty r\circ (\g,\1\g) = 0$
%for the Taylor series $Tr$ at $p$ of 
%for some (and hence any) smooth local defining function $r$ of $S$ 
%in a neighborhood of $p$;
%\item for every $q$-dimensional smooth complex subbundle 
%$E\subset H^{10}S$;
\item 
there does not exist a smooth complex subbundle 
$E\subset H^{10}S$ of rank $q$ in a neighborhood of $p$
such that
$$
[L^m,\ldots,[L^{2},L^1]\ldots](p)\in 
\C HS
%H^{10}_p S + \1{H_p^{10} S}
$$
holds
for all $m\ge2$ and all vector fields 
$L^m,\ldots,L^1\in \G(E)\cup \G(\1E)$;
% of infinite commutator type $t(E,p)$;
%for some (and hence any) smooth local defining function $r$ of $S$ 
%in a neighborhood of $p$.
\item 
there does not exist a smooth
 complex subbundle 
$E\subset H^{10}S$
of rank $q$
 in a neighborhood of $p$
such that
$$
L^m \cdots L^{3} \d r([L^{2},L^1])(p)= 0
$$
holds
for all $m\ge2$ and all vector fields 
$L^m,\ldots,L^1\in \G(E)\cup \G(\1E)$,
where $r$ is a local defining function of $S$ at $p$.
\een
Then 
$S$ satisfies property $(P_q)$.
In particular,
 the $\bar\d$-Neumann operator of $D$ on $(0,q)$ forms is
 {\em compact},
hence is {\em globally regular}
and $D$ 
{\em satisfies condition $R$}.
\et

%See \cite[Theorem~4.8]{St10} and \cite[Theorem~4.6]{St10}
%respectively 
%for proofs that $(P_q)$ implies compactness
%and that compactness implies global (and exact) regularity
%for $(0,q)$ forms.

\br
It is easy to see that conditions in the theorem do not depend
on the choice of the local defining function $r$.
\er

%Here $H_pS = T_p S \cap i T_p S$ and 
%$H^{10}_pS = \{X + iJX : X\in H_pS\} \subset \C\otimes H_p S$
%are respectively the complex tangent and $(1,0)$ subbundles on $S$.
\br
Condition (1) in Theorem~\ref{reg} 
is satisfied when $S$ is of {\em finite regular (contact) $q$-type} at $p$,
in particular,
 also when $S$ is of {\em finite (singular) $q$-type}
(or simply {\em finite D'Angelo type} for $q=1$)
 as defined in  \cite{D82,D93}.
%hence also when $S$ is of {\em finite D'Angelo type} at $p$.
Recall \cite[\S4.3.3]{D93} that 
the {\em regular $q$-type of $S$ at $p\in S$}
is the supremum of {\em contact orders} at $p$ of $S$ with 
%complex $q$-dimensional sub
germs
of holomorphic immersions 
$\g\colon (\C^q,0)\to (\cn,p)$,
i.e.\ of vanishing orders at $0$ of $r\circ \g$.
%where $r$ is a local defining function of $S$.
%A hypersurface of finite regular 
%$q$-type at every point 
%need not be of finite (singular) $q$-type, 
%and the regular type
%need not even be locally bounded, 
%unlike finite $q$-type \cite{D82}.
%
Indeed, if (1) in Theorem~\ref{reg} were violated
for a 
formal holomorphic immersion 
$\g
\colon (\C^q,0)\to (\cn,p)
$
tangent to $S$ of infinite order,
truncating $\g$ to arbitrary high orders
 would yield complex $q$-dimensional submanifolds
of arbitrarily high contact orders with $S$,
 contradicting the finiteness of the regular $q$-type.
\er

To relate contact orders in condition (1)
with subbundles, we introduce
the following notions:

\bd
\Label{E-contact}
Let $S\subset \cn$ be a smooth real submanifold
and $p\in S$ a fixed point.
\ben
\item
The {\em complex formal orbit}
$\6O^\C_E(p)$
of a smooth (real or complex) subbundle $E\subset \C TS$
at $p$
is the complex formal variety
given by the complex formal power series ideal
$$
I = \{ f\in \C[[z-p]] 
: 
L^t \cdots L^1 f(p)=0,
\;
L^t,\ldots,L^1\in \G(E)\cup\G(\1E),
\;
t\ge 0
 \}.
$$
\item
The {\em contact order}
at $p$
between $S$ 
and a formal variety $V$ given by
an ideal 
$I(V)\subset \C[[z-p]]$
in the ring of formal power series  
is
$$
\max 
\{ 
k
:
I_p(S) \subset I(V) + \5m^k
\}
\in \N_{\ge1}\cup\{\infty\}
,
$$
where $I_p(S)\subset \C[[z-p]]$
is the ideal of formal power series vanishing on $S$
(in the formal sense)
and 
$\5m\subset \C[[z-p]]$ is the maximal ideal.
\een
\ed

Note the case $t=0$ in the above definition if $I$ corresponds
to the vanishing of $f(p)$.
We use the above terminology 
to refine condition (1) in Theorem~\ref{reg}.
Conditions (2) and (3) in that theorem
involve respectively
iterated commutators of certain vector fields
and repeated differentiation of the Levi form
along certain $(1,0)$ and $(0,1)$ vector fields.
They
can be equivalently formulated
as the existence of complex rank $q$ subbundles
$E\subset H^{10}S$
with infinite invariants 
$t(E,p)$ and $c(E,p)$
 respectively,
defined below,
where by convention,
the minimum of an empty subset of $\N$ is
$\infty$.

\bd
\Label{types}
Let $S\subset \cn$, $n\ge2$, be a real smooth hypersurface
with local defining function $r$ in a neighbourhood of $p\in S$,
and $E\subset H^{10}S$ a smooth complex subbundle.
\ben
\item
The {\em contact type}
$a(E,p)\in \N_{\ge2}\cup \{\infty\}$
of $E$ at $p$  is the contact order
at $p$ between $S$
and the complex formal orbit
$\6O^\C_E(p)$
(see Definition~\ref{E-contact});
\item
The  {\em (commutator) type 
$t(E,p)\in \N_{\ge2}\cup \{\infty\}$
of $E$ at $p$}  is 
$$
%\begin{multline*}
\min\{
t\ge 2:
\exists
L^t, \ldots, L^1\in  \G(E)\cup \G(\1E),
%\\
[L^m,\ldots,[L^{2},L^1]\ldots](p)
\not\in 
%H^{10}_p S + \1{H_p^{10} S}
\C HS
\}
.
%\end{multline*}
$$

\item
 The {\em Levi (form) type
 $c(E,p)\in \N_{\ge2}\cup \{\infty\}$
  of $E$ at $p$} is 
$$
%\begin{multline*}
%c(E,p):=
\min\{
t\ge 2:
\exists
L^t, \ldots, L^1\in \G(E)\cup \G(\1E),
%\\
L^m \cdots L^{3} \d r([L^{2},L^1])(p)\ne0
\}.
%\in \N_{\ge2}\cup \{\infty\}.
%\end{multline*}
$$
%and is indepenendent on the choice of $\th$.
\een
\ed

With this definition, conditions  (2) and (3)
in Theorem~\ref{reg}
can be restated as

\ben
\item[(2)]
there does not exist a smooth complex subbundle 
$E\subset H^{10}S$ of rank $q$ in a neighborhood of $p$
with $t(E,p)=\infty$, i.e.\ of infinite commutator type;
\item[(3)]
there does not exist a smooth complex subbundle 
$E\subset H^{10}S$ of rank $q$ in a neighborhood of $p$
with $c(E,p)=\infty$, i.e.\ of infinite Levi type.
\een

Finiteness conditions such as in Theorem~\ref{reg} (1)-(3)
and notions of type as in Definition~\ref{types} (2)-(3)
in the study of 
 the $\bar\d$-Neumann problem 
 go back to Kohn \cite{K72,K74}, 
 Bloom-Graham~\cite{BlG77}, 
 Bloom~\cite{Bl81}
 and D'Angelo~\cite{D79,D82,D86a,D86b},
 see also \cite[\S4.3.1]{D93} and \cite{DK99}.
The invariant $t(E,p)$
in Definition~\ref{types}
 coincides
  with the one defined by
  Kohn \cite{K72,K74}
  for rank $1$ subbundles and
   Bloom~\cite[Definition~2.5]{Bl81} in general case,
  while $c(E,p)$ 
  has been shown
by
 Huang-Yin~\cite{HY23} to coincide 
  with the respective invariant defined Bloom~\cite[Definition~2.9]{Bl81}
 for {\em pseudoconvex hypersurfaces}.
  On the other hand, the 
  {\em contact type $a(E,p)$}
 of a subbundle
 in Definition~\ref{types} (1) is new.

In \cite{K74}, Kohn formulated a conjecture
that a condition similar to (2) in Theorem~\ref{reg} is 
{\em necessary and sufficient for subellipticity}
(later disproved by D'Angelo and Catlin).
Remarkably, Theorem~\ref{reg} shows the importance
of Kohn's condition for the weaker properties of
compactness and global regularity rather than subellipticity.

Other commutator conditions,
e.g.~\cite{DeT76, De78,De80,FeK90,De17},
and conditions on derivatives of the Levi form,
e.g.~\cite{HR22},
have been used
 in the study of the
$\bar\d$-Neumann problem.

More broadly, finiteness conditions of a similar nature
have appeared in problems of Geometry and PDEs.
For instance,
a condition similar to (1) in Theorem~\ref{reg}
involving singular formal varieties rather than regular ones
has been used for the problem of
finite jet determination of smooth CR maps in
 \cite{BER00}.
 Condition (3) in Theorem~\ref{reg}
 based on differentiation of the Levi form
 is somewhat parallel
 to {\em finite nondegeneracy} condition
  crucial for rigidity properties of CR maps,
 where the Levi form is differentiated only in 
 $(0,1)$ directions, 
 see \cite{BHR96,E-jdg,BERbook}
and
 the recent survey \cite{LM22} and references therein.
Iterated commutator type conditions
have been prominent in the Analysis of PDEs,
see \cite{NT63,Ho67,RoS76,K05},
and more recently in e.g~\cite{AB12}.

%Note that Kohn's condition in Conjecture~\ref{kohn} and Theorem~\ref{fin-order} is {\em intrinsic} in terms of brackets of vector fields. 
%The minimum number of brakets needed to obtain a nontangential direction
%for a fixed $(1,0)$ vector field $L$ is called the {\em type of $L$} in Bloom-Graham~\cite[\S3]{BlG77} and its 
%{\em supremum over all germs of vector fields} 
%at a point is considered by Bloom~\cite{Bl81}. 
%See also D'Angelo~\cite{D86a, D86b}, \cite[\S4.3.1]{D93} and 
%D'Angelo-Kohn~\cite{DK99} for more recent developments. 
%We note that
%Kohn's original condition in Conjecture~\ref{kohn} and Theorem~\ref{fin-order}
%allows for the types of individual vector fields being finite 
%even when their supremum is infinite.

The exact relationship between 
 contact orders and types 
 in Definition~\ref{types}
 for {\em pseudoconvex hypersurfaces}
remains a difficult problem and
part of the 
 Bloom's and D'Angelo's conjectures \cite{Bl81,D86b},
 both still open in general, 
 with only recent substantial progress made by 
Huang-Yin~\cite{HY21}
that initiated further research in
 \cite{CCY21,HY23,HYY23}. 

\subsection{Special subbundles and tower multitype}
We now turn to 
{\em more refined and stronger versions} of
conditions (1)-(3) in
Theorem~\ref{reg}, where the same conditions
only need to be verified for a more limited class of
%formal immersions $\g$
and subbundles $E$.
While condition (1) was formulated
in terms of formal immersions
rather than subbundles,
we use the terminology
of Definition~\ref{E-contact}
to formulate a refinement
of this condition also in terms of subbundles.

It will be convenient to replace $\d r$ appearing
in condition (3) in Theorem~\ref{reg}
with a more general {\em complex contact form} $\th$
in a neighborhood of $p$ in $S$,
by which
we mean any nonzero $\C$-valued $1$-form 
 vanishing on  $\C HS = H^{10} S \oplus \1{H^{10} S}$. 
For simplicity, we replace $S$ with an open subset
if needed,
where such contact form is defined.
%Working with contact forms on $S$
%makes our definition of a tower
%intrinsically defined in terms 
%of the CR structure 

We shall introduce 
{\em dual forms}
defined canonically on
the $(1,0)$ subbundle $H^{10}S$
of a real hypersurface $S$
(or more generally, on any abstract CR manifold
with a chosen equivalence class of contact forms
up to a multiplication by a nonzero function),
and consider special complex subbundles of $H^{10}S$
defined by linearly independent collections of such forms.
%We then consider formal holomorphic immersions
%$
%\g
%\colon (\C^q,0)\to (\cn,p)
%$
%tangent $S$:

\bd[Dual forms and special subbundles]
\Label{dual-f}
Let $S\subset \C^n$, $n\ge2$, be a smooth real hypersurface
%(or any smooth abstract CR manifold of hypersurface type)
with a complex contact form $\th$.
\Label{special-sub0}
\ben
\item
The {\em $\th$-dual form} of an (ordered)
 list of complex vector fields
$$
%\beq
L^{t},\ldots,L^1 \in \G(H^{10}S)\cup\G(\1{H^{10}S})),
\quad t\ge 1,
%\eeq
$$
is
the complex $1$-form 
$\om_{L^{t},\ldots,L^1;\th}$
on $H^{10}S$
defined 
for $L\in \G(H^{10}S)$, $p\in S$,
by 
%for $k=2$ is uniquely determined by
$$
\bca
%\Label{om-dual}
%\om_r(L) = L r& k=0\\
\om_{L^1;\th}(L_p) := \th ([L, L^1])(p) & t=1\\
\om_{L^{t},\ldots,L^1;\th}(L_p) :=
(Lf_{L^{t},\ldots,L^1;\th})(p), \;\;  
f_{L^{t},\ldots,L^1;\th} := \Re (L^{t} \cdots L^3 \th([L^2, L^1])) & t\ge 2
\eca.
$$
%where for $t\ge2$, we call $f_{L^{t},\ldots,L^1;\th}$ 
%the {\em $\th$-dual function} of the collection $L^{t},\ldots,L^1$.
\item
A complex subbundle $E\subset H^{10}S$
is called {\em special} if it can be defined by
$$
E=
\{
\xi\in H^{10}S :
\om_1(\xi)= \ldots= \om_l(\xi)= 0
\},
\quad
\om_1\wedge\cdots\wedge\om_l\ne 0
\text{ on } (H^{10}S)^l,
$$
where 
 each  $\om_j$, $j=1,\ldots,l$,
is the $\th$-dual $1$-form
$\om_j = \om_{L^{t_j}_j,\ldots,L^1_j}$
for some integers $t_j\ge1$
and vector fields
$
L^{t_j}_j,\ldots,L^1_j
\in \G(H^{10}S)\cup\G(\1{H^{10}S}))
$.
%\item
%A formal holomorphic immersion 
%$$
%\g
%% = \sum_\a c_\a z^\a \
%\colon (\C^q,0)\to (\cn,p)
%$$
%is called {\em special}
%if it is tangent to $S$ of infinite order
%and tangent to a special subbundle $E\subset H^{10}S$
%in the formal sense,
%i.e.\ any vector field $L$ from $E$
%preserves the real formal ideal 
%of the formal
\een
\ed

It is easy to see from the definition that
$
\om_{L_k^{t_k-1},\ldots,L_k^1; h\th} = \om_{hL_k^{t_k-1},\ldots,L_k^1;\th}
$
holds for any complex smooth function $h$.
Since any other complex contact form
is $\th'=h\th$ for some nonzero function $h$,
%$\th'$-dual forms are the same as $\th$-dual forms
%for different choices of vector fields,
the set of all $\th$-dual forms
and the notion of a special subbundle
only depend on the CR structure of $S$
rather than on the choice of $\th$.
%see Remark~\ref{th-dep} (3) below.

With Definitions~\ref{E-contact},~\ref{types},~\ref{special-sub0}, we can now
strengthen Theorem~\ref{reg} by refining
conditions (1)-(3) as follows:

\bt
\Label{refined}
The conclusion of Theorem~\ref{reg} still holds
with conditions \rm{(1)-(3)}
stated only for {\em special subbundles} $E\subset H^{10}S$
of rank $\ge q$,
i.e. when for every $p\in S$, at least one of the following conditions holds:
\ben
\item[(1$'$)]
there does not exist a {\em special} smooth complex subbundle 
$E\subset H^{10}S$ of rank $\ge q$ in a neighborhood of $p$
with $a(E,p)=\infty$, i.e.\ of infinite contact type,
and whose complex formal orbit $\6O^\C_E(p)$ is regular,
i.e.\ defined by a formal manifold ideal
(see Definition~\ref{m-ideal});
\item[(2$'$)]
there does not exist a {\em special} smooth complex subbundle 
$E\subset H^{10}S$ of rank $\ge q$ in a neighborhood of $p$
with $t(E,p)=\infty$, i.e.\ of infinite commutator type;
\item[(3$'$)]
there does not exist a {\em special} smooth complex subbundle 
$E\subset H^{10}S$ of rank $\ge q$ in a neighborhood of $p$
with $c(E,p)=\infty$, i.e.\ of infinite Levi type.
\een
\et

Note that if a subbundle $E$ of rank $\ge q$
is of infinite contact, commutator or Levi type at $p$,
any rank $q$ subbundle of $E$ is automatically
of infinite commutator or Levi type at $p$,
and any $q$-dimensional formal submanifold of the orbit
$\6O^\C_E(p)$ is automatically of infinite contact order.
Hence, conditions (1$'$)-(3$'$) in Theorem~\ref{refined}
 are implied
by the respective conditions (1)-(3) in Theorem~\ref{reg}.
(To see that (1) implies (1$'$), observe that 
any formal submanifold admits a parametrization
by a formal immersion, see \S\ref{mfd-ideals}).
On the other hand, we cannot expect those rank $q$ subbundles
to be also {\em special} (in the sense of Definition~\ref{special-sub0},
which is the reason conditions (2$'$) and (3$'$)
are formulated for rank $\ge q$.

%The inclusion of the assumption in (1$'$) of 
%the complex formal orbit being {\em regular}
%allows to limit the choice of $E$
%to only those 

A yet more general finiteness condition,
of which conditions of Theorem~\ref{reg} and~\ref{refined}
will be proved to be special cases,
can be given in terms of the {\em tower multitype},
which is one of the main new technical tools of this paper.
It will be convenient to replace $n$ with $n+1$, $n\ge 1$.

\bd
\Label{tm}
Let $S\subset\C^{n+1}$, $n\ge1$,
 be a smooth real hypersurface
with a complex contact form $\th$.
\ben
\item
A complex $1$-form $\om$ defined on $H^{10}S$
is called {\em $E$-dual of order $t\in\N_{\ge2}$,
where $E\subset H^{10}S$} is a complex subbundle,
if it is $\th$-dual of a list of $(t-1)$ complex vector fields
$$
%\beq
L^{t-1},\ldots,L^1 \in \G(E)\cup\G(\1{E})).
%\eeq
$$
\item
A {\em tower on $S$ of multi-order}
 $(t_1,\ldots,t_n)\in (\N_{\ge2}\cup\{\infty\})^n$
% where $2\le t_1\le\ldots\le t_n\le\infty$,
is a nested sequence of
complex subbundles
%\beq\Label{tower}
$$
H^{10}S = E_0\supset \ldots \supset E_m,
\quad
 0\le m\le n,
%\eeq
$$
such that $t_{m+1}=\ldots=t_n=\infty$, and
 for each $k=1,\ldots,m$, one has $t_k\in \N_{\ge2}$ and
 there exists an $E_{k-1}$-dual form 
 $
\om_k
$ 
 of order $t_k$ with
%that
%does not vanish on
%(fibers of) $E_{k-1}$,
%and such that $E_k$ satisfies
$$
E_{k}=E_{k-1} \cap \{\om_k=0\},
\quad
\omega_k|_{E_{k-1}}\ne 0.
$$
\item
The {\em tower multitype} of $S$ at $p\in S$
is the CR invariant
$$
\6T(p) 
%= \min\{(t_1,\ldots,t_n)\} 
\in (\N_{\ge2}\cup \{\infty\})^n
$$ 
defined as
the {\em lexicographically minimum multi-order}
$(t_1,\ldots,t_n)$
 of a tower 
% $H^{10}S = E_0\supset \ldots \supset E_m$
on a neighborhood of $p$ in $S$. 
%\item
%The tower multitype $\6T(p)$
%\een
\een
\ed

\brs
\ben
\item
In part (2) of Definition~\ref{tm}, in particular,
the second subbundle $E_1\subset E_0= H^{10}S$
can be defined by any $1$-form 
$\om=\om_{L^{t-1},\ldots,L^1}$
as in Definition~\ref{dual-f}, part (1), that does not vanish on
(each fiber of) $E_0$.
\item
Here the
{\em lexicographic order} is defined in the standard way by
%\begin{multline*}
$$
(t_1,\ldots,t_n) < (t'_1,\ldots,t'_n)
\iff
%\\
\exists 
\;
k\in \{1,\ldots,n-1\},
(t_1,\ldots,t_{k-1})=(t'_1,\ldots,t'_{k-1}),
\;
 t_k < t'_k.
 $$
%\end{multline*}
Taking the lexicographic order in (3)
guarantees that $\6T(p)$ is an invariant
only depending on the CR structure of $S$.
%similarly to Catlin's multitype~\cite{C84ann}.
It follows directly from the definitions
that all subbundles $E_k$ in a tower
are {\em special}
in the sense of Definition~\ref{special-sub0}.
\item
The choice of $t-1$ for the number of vector fields
in the list is made so that the tower multitype
in $\C^2$, which is a single number,
coincides with all other types.
\een
\ers

The following theorem provides
new finiteness conditions 
in terms of the tower multitype
guaranteeing property $(P_q)$:
\bt
\Label{main}
Let $D\subset\cn$ be a pseudoconvex bounded domain 
with smooth 
$(C^\infty)$ boundary $S=\d D$,
and $1\le q\le n-1$ be an integer. Assume that
 for every boundary point $p\in S$,
the tower multitype $\6T(p)=(t_1,\ldots,t_n)$ satisfies
\beq
\Label{t<q}
\# \{k: t_k = \infty\} < q,
\eeq
where $\#$ is the number of elements in the set.
Then $S$ satisfied property $(P_q)$.
\et

Theorem~\ref{main}
will be proved
at the end of \S\ref{towers-section},
where we use the terminology of {\em $q$-finiteness}
for the assumption \eqref{t<q}.
It is one of the reasons
behind our definition of the tower multitype
to make this part of the proof
particularly transparent and only taking one section.

Once Theorem~\ref{main} is established,
both
Theorems~\ref{reg} and~\ref{refined}
will be reduced to Theorem~\ref{main}
by showing that their conditions
are in fact special cases of the one in Theorem~\ref{main}:

\bt
\Label{reduce}
Let $D\subset\C^{n+1}$, $n\ge1$,
 be a pseudoconvex bounded domain 
with smooth 
$(C^\infty)$ boundary $S=\d D$,
and $1\le q\le n$ be an integer.
Then for any fixed $p\in S$, each of the conditions {\rm (1)-(3)}
in Theorem~\ref{reg}
and {\rm (1$'$)-(3$'$)}
in Theorem~\ref{refined}
imply \eqref{t<q}
for the tower multitype
$\6T(p)=(t_1,\ldots,t_n)$.
\et

The scheme of implications between various
conditions in Theorems~\ref{reg}, \ref{refined} and \ref{main}
can be illustrated by the diagram:
$$
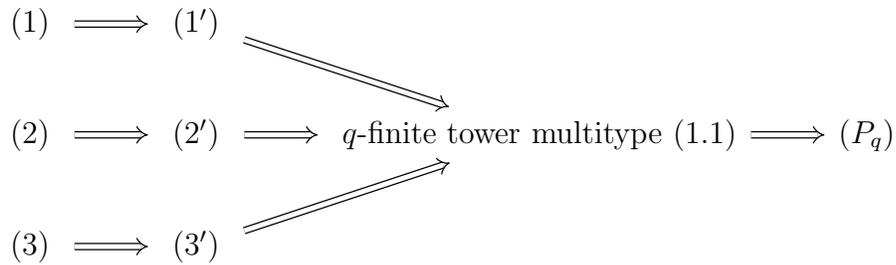
\begin{tikzcd}
%[column sep=huge]
	\text{ (1) } \rar[Rightarrow]
	& \text{ (1$'$) } \drar[Rightarrow]
	&
	&
	\\
	\text{ (2) }\rar[Rightarrow] 
	& \text{  (2$'$) } \rar[Rightarrow] 
	& \text{ $q$-finite tower multitype } \eqref{t<q} 
	\rar[Rightarrow]
	%{\text{ Theorem~\ref{main}}}
	& \text{$(P_q)$}
	\\
	\text{ (3) } \rar[Rightarrow] 
	& \text{  (3$'$) } \urar[Rightarrow]
	&
	&
\end{tikzcd}
$$
Implications in the first column are simple
consequences of the definitions
as mentioned above after Theorem~\ref{refined}.
It is clear that both Theorems~\ref{reg} and~\ref{refined}
are consequences of  Theorems~\ref{reduce} and~\ref{main},
corresponding respectively to the 2nd and 3rd column implications
in the above diagram.

Proving Theorem~\ref{reduce}
is
taking up the rest of the paper starting from \S\ref{inf}.
The proof requires several new tools 
outlined before in \S\ref{new-tech}.
%taking formal Nagano orbits of subbundles
%violating condition (3) in Theorem~\ref{reg},
%and on a formal version of the celebrated result by Diederich-Fornaess~\cite{DiF78} relating finite type with
%Kohn's {\em holomorphic dimension} \cite{K79}. 
However, the number of steps involved
in proving Theorem~\ref{reduce}
under different conditions (1$'$), (2$'$) or (3$'$)
varies substantially.

The Levi type condition (3$'$) (and hence (3))
is the closest one to Definition~\ref{tm},
for which the proof 
is completed already at the end of
\S\ref{inf}.
The proof for
the commutator type condition (2$'$) (and hence (2))
takes some more effort and 
tools such as formal orbits $O$
and careful study of the Levi form along $O$.
The proof of this part is completed
at the end of \S\ref{form-orb}.
Finally, the proof
for the contact type condition (1$'$) (and hence (1))
requires the most efforts and 
the entire list of new tools outlined 
in \S\ref{new-tech},
and is only completed in \S\ref{pf-reg}.

\bigskip

We conclude this section by an example
(a modification of D'Angelo's~\cite[Example~6]{D86a})
demonstrating that 
the condition in Theorem~\ref{main}
is, in fact, more general than all conditions of
 Theorems~\ref{reg} and~\ref{refined}.
 
\bex
\Label{t-type}
Let $S\subset\C^3_{z_1,z_2,w}$ be
the pseudoconvex real hypersurface defined by
$r=0$, where
$$
r = -2\Re w + |z_1z_2|^2 + f(z_1,z_2),
$$
where $f$ is any smooth 
real function in $\C^2$
that is
strictly plurisubharmonic 
away from $0$ and
 vanishes to infinite order at $0$.
Then, for $p\ne 0$,
 $S$ is strictly pseudoconvex,
hence it is easy to see that the tower multitype $\6T(p)=(2,2)$.

Further, we claim that $\6T(0)=(4,4)$.
Indeed, fixing the contact form $\th=\d r$
and the vector field $L\in H^{10}S$
of the form
$$
L = \d_{z_1} + \d_{z_2} + a(z_1,z_2)\d_w
$$
with suitable function $a$ chosen to make it tangent to $S$,
we have
$ \om_1|_{E_0}\ne 0$ for $E_0=H^{10}S$ and
the $1$-form
$
\om_1=\om_{\1L,L,\1L;\th}
$
in the notation of Definition~\ref{special-sub0},
i.e.\
$$
\om_1=d\Re \1L \d r([L,\1L)
= d \Re \1L L\1L r
= c(dz_1+dz_2),
$$
where $c\ne0$ is a constant.
Hence $\om_1 = c(dz_1+dz_2)$
and
$t_1=4$ is the minimum choice 
for the first entry of the multiorder of
a tower, hence of
the tower multitype 
$\6T=(t_1,t_2)=(4,t_2)$.
Setting 
$$
E_1=E_0\cap\{\om_1=0\},
\quad E_2=0,
$$ 
we obtain a tower $E_0\supset E_1\supset E_2$
of the minimum multi-order $(t_1,t_2)=(4,4)$.
Thus \eqref{t<q} is satisfied for all points $p$
for $q=1$.

On the other hand, 
since $S$ is tangent of infinite order
to the complex lines $w=z_1=0$
and $w=z_2=0$,
condition (1) in Theorem~\ref{reg}
is violated for $q=1$.
Taking the subbundle $E\subset H^{10}S$
defined by $dz_2=0$ and $L'=\d_{z_1}+a'(z_1,z_2)\d w$,
it can be seen that conditions (2) and (3) are also violated.
Moreover, $E$ can be seen to be special,
hence also conditions 
{\rm (2$'$)-(3$'$)}
in Theorem~\ref{refined}
are violated.
\eex

%The proof of Theorem~\ref{reg} is completed in \S\ref{pf-reg}.

%as defined by Kohn~\cite{K72} in $\C^2$ and D'Angelo~\cite{D82} in general $\C^n$.
%Conditions (2) and (3) are intrinsic and 
%Note that Theorem~\ref{reg} immediately implies
% the same conclusion for $(p,q)$ forms with $0\le p\le n$ instead of $(0,q)$ forms,
%since $\bar\d$ acts on the coefficients of the $dz^\a$ differentials.
%See \cite[Proposition~4.2]{St10} for conditions equivalent to the compactness
%of the $\bar\d$-Neumann operator.

%whose notable applications include
% celebrated {\em Condition R} and resulting smooth extension
%of proper holomorphic maps \cite{BeL80,Be81,BeC82}.
%Unlike subellipticity, for which the finite type is necessary
%by another work of Catlin~\cite{C83},
%global regularity is known to hold for more general classes 
%of domains, see e.g.\ 
%\cite{Si87,B88, BS91,HI97,FS98,K99,ChF05,HeM06,Ha11,PZ14,BS99, FS01,St10}.

\subsection{Countable generalized stratifications}
\Label{reg-sec}
In analogy with Catlin's argument in \cite{C84proc}
for establishing property $(P_q)$,
we employ a ``stratification'' of the given hypersurface
into subsets (not necessarily submanifolds, 
hence our cautious term ``generalized stratification'') 
having certain convexity properties
relative to the Levi form.
See also the diagram in \S\ref{path}
for implications using generalized stratifications.
%Our main technical tool is the 
%{\em tower multitype} (see Definition~\ref{tmt} below)
%whose level sets define a of the 
%real hypersurface $S\subset\C^n$
%similarly 
%that plays a role similar to Catlin's multitype
%in defining a 
However, to handle the more general conditions
in Theorems~\ref{reg} and \ref{refined},
we need to allow for 
{\em countable}
rather than finite stratifications.
In doing so, we are led to the following
generalization of 
the {\em weakly regular} boundaries
by Catlin
\cite{C84proc} and 
{\em regular domains} by Diederich-Fornaess 
\cite{DiF77stein} 
(corresponding to finite ``stratifications'' and $q=1$).
As before, we shall fix an integer $1\le q \le n-1$.

\bd
\Label{reg-def}
%[countably $q$-regular hypersurfaces]
We call a smooth real hypersurface $S\subset\C^n$
 {\em countably $q$-regular} if
it can be represented as a 
countable disjoint union 
$
S = \cup_{k=1}^\infty S_k
$
of locally closed subsets $S_k\subset S$
(``strata'')
such that for each $k$ and $p\in S_k$,
there exists a CR submanifold $M\subset S$ 
satisfying the following properties:
\ben
\item
$M$ contains an open neighborhood of $p$ in $S_k$
(in the relative topology);
\item
$\dim_\C(H^{10}_x M\cap K^{10}_x)<q$ for all $x\in M$,
where $K^{10}_x\subset H^{10}_x S$ 
is the kernel of the Levi form of $S$.
\een
\ed

\br
The dimension $d(x):=\dim_\C(H^{10}_x M\cap K^{10}_x)$
was used by Kohn
in his important notion of ``holomorphic dimension''
of $M$
as the {\em minimum} of $d(x)$ for all $x$
\cite[Definition~6.16]{K79},
arising in the study of Kohn's multiplier ideals.
Later Catlin~\cite{C84ann}
considered the {\em maximum}
of $d(x)$, also called ``holomorphic dimension'',
a closely related but different notion
arising in the context of generalized stratifications.
Condition (2) in Definition~\ref{reg-def}
is equivalent to Catlin's holomorphic dimension
being less than $q$.

Note that we also consider $d(x)$ 
for non-pseudoconvex domains 
when the Levi kernel $K^{10}_x$
may not coincide with the Levi
null cone $\{X:\d\bar\d r(X,X)=0\}$
(where $r$ is a defining function of $S$ as before).
\er

\br
In a private communication,
Dall'Ara and Mongodi provided an argument
showing that 
a real hypersurface $S$
is {\em countably $1$-regular}
if and only if it has 
{\em empty Levi core}
in their terminology~\cite{DaM23}.
%for $q=1$,
%Definition~\ref{reg-def}
%is equivalent to their 
This can be combined
with a recent work of Treuer~\cite{T23}
showing that empty Levi core
implies property $(P_1)$.
\er

{\em Countably
$q$-regular}
 hypersurfaces still satisfy property $(P_q)$:
\bp
\Label{reg-pq}
Let $S\subset \C^n$ be
 a pseudoconvex smooth real hypersurface.
 Assume that $S$ is 
 {\em contably $q$-regular}.
 Then $S$ satisfies property $(P_q)$.
\ep

\br
This is the only step using pseudoconvexity
in proving Theorem~\ref{main}.
A proof of this proposition under more general assumptions
than pseudoconvexity
will lead to respective extensions of 
Theorem~\ref{main} to non-pseudoconvex situations.
\er

Once Proposition~\ref{reg-pq} is established,
the proof of
%Theorem~\ref{reg} and Theorem~\ref{refined}
Theorem~\ref{main}
will be completed by the following statement:

\bp
\Label{reg-conclude}
Let $D$ satisfy the assumptions 
of 
%either 
%Theorem~\ref{reg} or 
%Theorem~\ref{refined}
%or 
Theorem~\ref{main},
{\em except pseudoconvexity}.
Then the boundary $\d D$ is 
countably $q$-regular.
In fact, the sets $S_k$
in Definition~\ref{reg-def}
can be chosen to be 
the level sets of the tower mult-type function
$\6T$ (see Definition~\ref{tm}).
\ep

We mention that
% the
%(generalized)
strata of constant multitype
also occur naturally in other geometric
problems such as constructions
of convergent normal forms
and invariant connections, see
\cite{KoZ15,KoZ19}.

%This seems to provide the first known approach
%to {\em global regularity for general 
%smoothly bounded pseudoconvex domains of finite type}
%that does not make use of results in \cite{C84ann}.

The proof Proposition~\ref{reg-conclude}
will be completed at the end of \S\ref{towers-section},
while
the proof of
Proposition~\ref{reg-pq}
%deviates from Catlin's original argument in
%\cite{C84proc}[Proof of Theorem 1]
%and is instead 
is based on the following results
by Sibony~\cite{Si87} for $q=1$
and their generalizations by Fu-Straube~\cite{FS01} for $q>1$:

\bp[Proposition~4.15 in \cite{St10}; for $q=1$ in \cite{Si87}, also implicit in \cite{C84proc}]
\Label{m-p}
Let $D \subset\cn$ be a pseudoconvex domain with smooth boundary $S=\d D$
and $M\subset S$ be a smooth submanifold such that 
$$
\dim_\C(H_p^{10}M\cap K_p^{10})<q, 
\quad p\in M,
$$
where $K_p^{10}\subset H^{10}_pS$
is the kernel of
the Levi form of $S$.
Then any compact subset of $M$ satisfies property $(P)_q$.
\ep

%We shall use the following basic properties
%\cite[Proposition~1.9]{Si87},
%\cite[Lemma~4.12 and Corollary~4.13]{St10}:

\bp[Corollary~4.14 in \cite{St10},
Proposition~3.4 in \cite{FS01}; for $q=1$ in \cite{Si87}]
\Label{countable-union0}
Let $K=\cup_{m=1}^\infty K_m$ be a countable union of compact subsets of $\C^n$, 
each satisfying $(P_q)$. Assume that $K$ is compact.
Then $K$ also satisfies $(P_q)$.
%\een
\ep

\bpf[Proof of Proposition~\ref{reg-pq}]
Assume 
the given hypersurface $S\subset\C^n$ 
to be countably $q$-regular as in the theorem,
and let 
$
S = \cup_{k=1}^\infty S_k
$
be a countable decomposition as in Definition~\ref{reg-def}.
%It is obvious from the definition of property $(P_q)$
%Hence,
In view of Proposition~\ref{countable-union0},
it suffices to show that each stratum $S_k$ 
is covered by a countable union of compact subsets of $S$
satisfying $(P_q)$.

To construct such a cover of $S_k$,
consider the countable set $\6B$ of all closed balls 
$B=\{z:|z-a|\le R\}\subset\C^n$
with both center $a$ and radius $R>0$  rational
such that the closure $\1{S_k\cap B}$ is contained
in a submanifold $M=M_B\subset S$ satisfying property (2)
in Definition~\ref{reg-def}.
Each such compact set $\1{S_k\cap B}$ satisfies $(P_q)$ by Proposition~\ref{m-p}.
Finally, 
property (1) in in Definition~\ref{reg-def} guarantees that
$S_k$ is covered by the countabily many sets $\1{S_k\cap B}$
with $B\in\6B$, completing the proof.
%
%By Corollary~\ref{strat}, the boundary $\d\Om$
%is a countable union $\cup M_k$ of submanifolds $M_k\subset\d\Om$
%as in Proposition~\ref{m-p}. Intersecting with closed balls
%with rational center and radius, each such $M_k$
%can be expressed as a countable union of its compact subsets,
%each of which satisfies P by Proposition~\ref{m-p}.
%Then $\d\Om$ is a union of countably many compact
%satisfying P, hence $\d\Om$ itself satisfies P by 
%Proposition~\ref{countable-union}.
\epf

\section{Tower multitype}
\Label{towers-section}
To prove Proposition~\ref{reg-conclude},
we introduce a new tool called
{\em tower multitype}
that assigns an $n$-tuple
$$
\6T(p)\in (\N_{\ge2}\cup\{\infty\})^n
$$
to any point $p$
of any smooth (abstract)
CR manifold $S$ of hypersurface type
of real dimension $2n+1$,
whose level sets yield
a stratification 
as in Definition~\ref{reg-def}.
For simplicity, we shall focus 
on the embedded case $S\subset\C^{n+1}$.

We shall also need the notion of 
a smooth {\em CR submanifold} $M\subset S$,
which is any real submanifold of $S$
such that  the dimension of 
$$
H_x^{10}M:= \C T_xM \cap H_x^{10}S, 
\quad
x\in S,
$$ 
is constant on $M$.
%for $x\in M$, where $H_x^{10}\cn$ is the space of all
%$(1,0)$ tangent vectors $\sum_j a_j \d_{z_j}$ at $x\in M$.
%We work with the complexified tangent bundle 
%$\C TM=\C\otimes TM$ and its complex subbundles.
By $H^{10}M \subset \C HM$ denote the $(1,0)$ bundle
by $HM=\Re H^{10}M\subset TM$ its real part,
hence $\C HM = H^{10}M\oplus \1{H^{10}M}$.
%We shall regard all forms on $M$
%extended to
%the complexified bundle $\C TM$ by $\C$-(multi-)linearity.
We shall omit the manifold $M$ whenever it is unambiguous,
i.e. write $H^{10}=H^{10}M$, $\C T=\C TM$ etc.
As a slight abuse of notation,
we shall write $L\in E$ (instead of $L\in \G(E)$)
for a complex vector field $L$
whenever $L$ is a section of a subbundle $E\subset \C TM$,
and write $L\in E_1\cup E_2$
when $L$ is a section of either $E_1$ or $E_2$.
%(defined over an open subset of $M$).

As before, by a {\em complex contact form} $\th$ on $S$
we mean any nonzero $1$-form on $S$
(i.e.\ nonzero at each point of $S$), 
which vanishes 
on $\C HS = H^{10}S + \1{H^{10}S}$. 
%Here a $k$-form  is always extended to the complexified
%tangent space by complex (multi-)linearity.
When $M\subset\C^{n+1}$ is a real hypersurface
with defining function $r$, one can choose $\th$
to be any nonzero complex multiple of $\d r|_M$.
We shall frequently abbreviate $H^{10}S$ to $H^{10}$ for brevity.

\subsection{Towers}
Our main new tool is a certain 
{\em nested sequence of complex 
subbundles} $E_s\subset H^{10}$, $s=0,1,\ldots$, with 
$$
\rk E_s=\rk H^{10}-s,
$$
 that we call a {\em tower}.
%It can be defined for any (abstract) CR manifold
%and a choice of a complex contact form,
%while we only need it here for smooth real hypersurfaces
%in $\C^{n}$.
For convenience, we include 
the definition of $\th$-dual forms (Definition~\ref{special-sub0} (1))
in part (1) below.

\bd[Tower]
\Label{tower-def}
Let $S\subset\C^{n+1}$ be a smooth real hypersurface with
a complex contact form $\th$.
\ben
\item
The {\em $\th$-dual form} of 
an (ordered) collection of complex vector fields
$$L^{t},\ldots,L^1 \in H^{10}\cup \1{H^{10}},\quad t\ge 1,$$
is
the complex $1$-form
$\om_{L^{t},\ldots,L^1;\th}$
on $H^{10}$
defined 
for $L\in H^{10}$, $p\in S$,
by 
$$
\bca
%\Label{om-dual}
%\om_r(L) = L r& k=0\\
\om_{L^1;\th}(L_p) := \th ([L, L^1])(p) & t=1\\
\om_{L^{t},\ldots,L^1;\th}(L_p) :=
Lf_{L^{t},\ldots,L^1;\th}(p), \;\;  
f_{L^{t},\ldots,L^1;\th} := \Re (L^{t} \cdots L^3 \th([L^2, L^1])) & t\ge 2
\eca,
$$
where for $t\ge2$, we call $f_{L^{t},\ldots,L^1;\th}$ 
the {\em $\th$-dual function} of the collection $L^{t},\ldots,L^1$.

\item
%Let $(t_1,\ldots,t_n)\in (\N_{\ge2}\cup\{\infty\})^n$
%satisfy $2\le t_1\le\ldots\le t_n\le\infty$.
A {\em tower on $S$ of the multi-order}
 $(t_1,\ldots,t_n)\in (\N_{\ge2}\cup\{\infty\})^n$,
% where $2\le t_1\le\ldots\le t_n\le\infty$,
is a nested sequence of
complex subbundles
%\beq\Label{tower}
$$
H^{10} = E_0\supset \ldots \supset E_m,
\quad
 0\le m\le n,
%\eeq
$$
such that $t_{m+1}=\ldots=t_n=\infty$, and
 for each $k=1,\ldots,m$, one has $t_k\in \N_{\ge2}$ and
there exists a collection of $t_k-1$ vector fields 
$$
L_k^{t_k-1}, \ldots, L_k^1 \in E_{k-1}\cup \1E_{k-1},
$$
for which the form
$
\om_k:=\om_{L_k^{t_k-1}, \ldots, L_k^1;\th}
$ 
does not vanish on
(fibers of) $E_{k-1}$,
and such that $E_k$ satisfies
$$E_{k}=E_{k-1} \cap \{\om_k=0\}.$$
\item
%Chose $l$ such that
%$t_k=2$ for $k\le l$ and $t_k>3$ for $k>l$.
For any choice of the vector fields $(L_k^s)$ as above,
collect all $\th$-dual functions as defined in (1)
for all $k$ with $t_k\ge 2$ into the set
$$
\{  f_{L_k^{t_k-1}, \ldots, L_k^1;\th} : t_k\ge 2\}
%\quad
%f_k := f_{L_k^{t_k-1}, \ldots, L_k^1;\th},
%\; k> l,
$$
that we call an {\em associated set of functions} of the given tower.
%we call all $\th$-dual functions
%collecting the $\th$-dual functions $f_{L_k^{t_k-1}, \ldots, L_k^1;\th}$ 
%for all 
%we call the collection 
%$$
%\{f_{l+1},\ldots, f_m\}
%$$
%$$
%(L_1,\ldots,L_l; f_{l+1},\ldots, f_m),
%\quad
%L_k := L_k^1,
%\; k\le l,
%\quad
%f_k := f_{L_k^{t_k-1}, \ldots, L_k^1;\th},
%\; k> l,
%$$
%a {\em defining system} of the given tower.
\een
\ed

\begin{samepage}
\begin{Rems}
\Label{th-dep}
\ben
\item
It is clear from part (2) of the definition that
the multi-order of a tower
$H^{10} = E_0\supset \ldots \supset E_m$
always satisfies
 $$
 (t_1,\ldots,t_n)
 \in
 (\N_{\ge2})^m\times (\{\infty\})^{n-m}.
 $$

\item
Note that for
$k=1,\ldots,l$,
the $1$-form
$\om_k=\th([\cdot,L_k])$
is a contraction of the
(restriction of)
 {\em Levi tensor}
$
\l_\th\colon \C HS \times \C HS \to \C,
$
defined by
$$
%\beq
%\Label{levi-id}
	\l_\th(L^2_p, L^1_p) = \th ([L^2, L^1]_p),
	\quad 
	L^2, L^1 \in \C HS,
	\quad
	p\in S.
%%\eeq
$$
The standard hermitian {\em Levi form}
is obtained by choosing a purely imaginary contact form $\th$
%is the multiple $i\l_\th(L^2,\1{L^1})$, 
for $L^2\in H^{10}S$, $L^1\in H^{01}S = \1{H^{10}S}$,
while $\l_\th$ vanishes on $H^{10}S\times H^{10}S$ and 
$H^{01}S\times H^{01}S$
as consequence of the integrability of the induced CR structure.
When $r$ is
a defining function of $S$ and 
$\th=c\d r$, where $c\ne 0$ is a real constant, we have
$\l_\th= -d\th = c\d\bar\d r$ by the Cartan formula for the exterior derivative.

\item
Even though
the $\th$-dual forms $\om_{L_k^{t_k-1},\ldots,L_k^1;\th}$ depend on
the choice of the contact form $\th$, 
the notion of a tower only depends on the (principal) 
ideal generated by $\th$,
which is uniquely determined by the CR structure for hypersurfaces.
Indeed, let
$
\2\th =  h\th
$,
where $h$ is a nonzero smooth complex function.
Then
$
\om_{L_k^{t_k-1},\ldots,L_k^1; h\th} = 
\om_{L_k^{t_k-1},\ldots,hL_k^1;\th}
$.
Hence, with $\2\th$ instead of $\th$,
we can modify the vector fields $(L_k^s)$ to obtain the same
forms $\om_k$ defining the same tower as claimed.
Since
any two complex contact forms differ by a multiplication with a
complex nonzero smooth function, the tower is defined invariantly independently of the choice of the complex contact form $\th$.
\een
\end{Rems}
\end{samepage}

%\subsection{Special subbundles}

\subsection{First structure property of towers}
The first main property of towers
is the existence of certain CR submanifolds of $S$
with controlled restiction of the Levi kernels.

\bp[First structure property]
\Label{first-structure}
Let $S\subset \C^{n+1}$ be a smooth real hypersurface
and
%(not necessarily pseudoconvex)
$$
H^{10}S = E_0\supset \ldots \supset E_m,
$$
a tower of the multi-order $(t_1,\ldots,t_n)$ on $S$ with 
%a defining system
%$(L_1,\ldots, L_l; f_{l+1},\ldots,f_m)$.
an associated set of functions 
$\{f_1,\ldots,f_l\}$.
Then the following hold:
\ben
\item
the
restrictions to $H^{10}S$ of the differentials 
$d f_1,\ldots,d f_l$ are linearly independent,
in particular, the zero set
$$
M:= \{f_1=\ldots=f_l=0\} \subset S
%\quad c\in \R^{m-s}.
$$
is a smooth CR submanifold;
\item
the kernel distribution (of varying rank)
$K^{10}
% = \{\xi : \l_\th(\xi,\cdot)=0 \}  
\subset H^{10} S$ of the Levi form of $S$
%tensor $\l_\th$
%of $S$ is 
% transversal to $M$
%in the sense that
satisfies
$$
H^{10} M\cap K^{10}   \subset E_m.
$$
%In particular, when the tower multitype $\6T(p)$ is finite, 
%one has $E_{m+1}=\{0\}$ and hence 
%$K^{10}_x\cap H^{10}_xM =\{0\}$ for all $x\in M$.
% the restriction of the Levi form of $S$ to $H^{10}M$ is nondegenerate
%(hence positive definite if $S$ is pseudoconvex).
\een
\ep

\bpf
%Choose vector fields $(L^s_k)$ as 
It follows from Definition~\ref{tower-def} (2)
that the forms 
$\om_1,\ldots,\om_m$
defined there
are linearly independent
when restricted to $H^{10}S$.
If the set $\{f_1,\ldots,f_l\}$ is empty, (1) is void.
Otherwise,
the equality of the sets of the forms
$$
\{
\om_k : t_k\ge 2
\}
=
\{
df_j|_{H^{10}S} : 1\le j\le l
\}
$$
%holds for some $l_k\in\{1,\ldots,l\}$ whenever
%$t_k\ge 2$,
 proves (1).

To show (2), let $\xi \in K^{10}$.
Since $\om_k = \om_{L_k;\th} = \l_\th(\cdot, L_k)$ 
when $t_k=1$
in the notation of Definition~\ref{tower-def} (1),
it follows that $\om_k(\xi)=0$.
On the other hand, when $t_k\ge 2$,
$\xi \in H^{10}M$
implies
 $\om_k(\xi)=df_{l_k}(\xi)=0$ for some $l_k\in \{1,\ldots,l\}$.
Hence 
$$
\xi \in K^{10} \cap H^{10} M
\implies \xi \in H^{10} S \cap \{ \om_1 = \ldots = \om_m = 0\} = E_m
$$
as desired.
\epf

\subsection{Tower multitype}
Using the notion of the tower (Definition~\ref{tower-def}), 
we define
a CR invariant at each point as follows:
\bd[tower multitype]
\Label{tmt}
Let $S\subset\C^{n+1}$ be a smooth real hypersurface.
%\ben
%\item
The {\em tower multitype} of a point $p\in S$
$$
\6T(p) 
%= \min\{(t_1,\ldots,t_n)\} 
\in (\N_{\ge2}\cup \{\infty\})^n,
$$ 
is the {\em lexicographically minimum multi-order}
$(t_1,\ldots,t_n)$
 of a tower 
% $H^{10}S = E_0\supset \ldots \supset E_m$
on a neighborhood of $p$ in $S$. 
%\item
%The tower multitype $\6T(p)$
%\een
\ed

In view of Remark~\ref{th-dep} (3),
the tower multitype $\6T(p)$ is independent
of the contact form $\th$ and hence is a CR invariant of
a germ of $S$ at $p$.
Since the multi-orders $(t_1,\ldots,t_n)$ in Definition~\ref{tmt}
consist of integers
and $\infty$, their lexicographic minimum
always
exists and is realized for certains towers.
Our second main property is for those towers
whose multi-orders realize
their lexicographic minimum ---
the tower multitype:

%
%The tower multitype can be compared to the Catlin multitype \cite{C84ann},
%also defined as lexicographic supremum of certain tuples of 
%{\em rational numbers} or $\infty$.
%In comparison, the {\em tower multitype} $\6T(p)$ 
%is a tuple of integers or 
%$\infty$, making it easy to show that
%$\6T(p)$ is realized as the multi-order of a fixed tower,
%since the lexicographic minimum of a tuple in 
%$(\N_{\ge2}\cup \{\infty\})^n$ is always achieved.
%Any tower on a neighborhood of $p$
%whose multi-order equals $\6T(p)$ is called 
%a {\em minimal tower at $p$}.
%

%In comparison, similar statements for Catlin's multitype are more elaborate due to the entries being rational rather than integers, and may require choices of coordinates given by formal power series.

\bp[Second structure property]
\Label{second-prop}
Let $S\subset \C^{n+1}$ be a smooth real hypersurface,
%(not necessarily pseudoconvex), 
$p\in S$ a point, $U\subset S$
an open neighborhood of $p$, and
$$
H^{10}U = E_0\supset \ldots \supset E_m
$$
a tower
 on $U$,
  whose multi-order equals the multitype $\6T(p)$.
 Choose
  any associated set of functions
$$\{f_1,\ldots,f_l\}$$
of the given tower
(see Definition~\ref{tower-def}  (3)).
Then the following hold:
\ben
\item 
$\6T(p') \le \6T(p)$ for any $p'\in U$;
\item
the tower multitype level set satisfies
$$
\{p'\in U : \6T(p') = \6T(p) \}
\subset  \{f_1=\ldots=f_l=0\}.
$$
\een
\ep

\bpf
Since $\6T(p)=(t_1,\ldots,t_n)$ 
is the multi-order of the given tower on $U$
and $\6T(p')$ is the minimum multi-order for a tower on 
a neighborhood of $p'$, 
(1) is immediate.

To show (2), let $p'\in U$ be such that
$f_j(p')\ne 0$ for some $j=1,\ldots,l$,
%and choose $k$ minimal with this property.
which 
in view of Definition~\ref{tower-def}, is of the form
\beq
\Label{fkq}
f_j(p')=
\Re (L_k^{t_k-1} \cdots L_k^3 \th([L_k^2, L_k^1]))(p')\ne 0,
\quad
L_k^{t_k-1}, \cdots, L_k^1 \in 
E_{k-1} \cup \1E_{k-1}
\eeq
for some $k$ that we choose to be minimal with this property
for any $j$.
%Splitting $L_k^{t_k-1}$ into sum of its components in  
%$E_{k-1}$ and $\1E_{k-1}$ and choosing a nonzero term,
%we may assume that \eqref{fkq} holds with
%either $L_k^{t_k-1}\in E_{k-1}$ or $L_k^{t_k-1}\in \1E_{k-1}$.
If  $L_k^{t_k-1}\in \1E_{k-1}$, taking conjugates of all vector fields
and of $\th$ and replacing $\th$
with $f \th$, where $f$ is a nonzero function,
we may assume that $L_k^{t_k-1}\in E_{k-1}$. In particular,
we obtain
$$
\om'_k|_{(E_{k-1})_{p'}}\ne 0,
$$
where
 for $x\in S$, $L\in \G(H^{10}S)$,
$$
\om'_k(L_x) := \th([L,L_k^1)(x) \text{ for } t_k=3,
$$
or
$$
\om'_k(L_x) := L L_k^{t_k-2} \cdots L_k^3 \th([L_k^2, L_k^1]))(x)
\text{ for } t_k>3.
$$
%where $x\in S$.
In the case $t_k=3$, we have $\om'_k = \om_{L_k^1;\th}$
in the notation of Definition~\ref{tower-def}.
In the case $t_k>3$,
splitting into real and imaginary parts, we obtain
$$
\om'_k(L_{p'})=
L
(\Re f + i\Im f)(p'),
\quad
f := 
L_k^{t_k-2} \cdots L_k^3 \th([L_k^2, L_k^1])).
$$
Taking the term that does not identically vanish for $L_{p'}\in (E_{k-1})_{p'}$
and multiplying $L_k^{t_k-2}$ by $i$ if necessary,
we may assume
$
\om'_k(L_x)=
L\Re f (x)
$, hence
$\om'_k = \om_{L_k^{t_k-2}, \ldots, L_k^1;\th}$.
In both cases, we obtain a new tower
$$
H^{10}U' = E_0 \supset \ldots \supset E_{k-1}\supset E'_{k}
$$
in a neighborhood $U'\subset U$ of $p'$ of
the lexicographically smaller multi-order 
$$
(t_1,\ldots,t_{k-1},t_k-1,\infty,\ldots,\infty)
<
(t_1,\ldots,t_{k-1},t_k,\ldots,t_n),
$$
by setting 
$$
E'_k:= E_{k-1} \cap \{\om'_k =0\}.
$$
In view of Definition~\ref{tmt}, it follows that
$\6T(p')<\6T(p)$, hence $p'$ is not in the level set of $p$,
completing the proof of (2).
\epf

Part (1) of Proposition~\ref{second-prop} immediately yields:

\bc
\Label{l-closed}
For a smooth real hypersurface  $S\subset \C^{n+1}$, 
the following hold:
\ben
\item
the tower multitype $\6T(p)$ is upper-semicontinuous as function of $p\in S$;
\item
any tower multitype level set $\{q\in S : \6T(q) = \const \}$
is locally closed, i.e. closed in its open neighborhood.
\een
\ec

\subsection{Hypersurfaces of $q$-finite tower multitype}
We now define a finiteness condition for the multitype
depending on the integer $q$
that will guarantee that the
decomposition into level sets of the tower multitype
defines
a stratification as in Definition~\ref{reg-def}.

\bd
\Label{q-finite-def}
An $n$-tuple $(t_1,\ldots,t_n)\in (\N_{\ge2}\cup\{\infty\})^n$
is called 
{\em $q$-finite}, where $1\le q\le n$ is an integer,
if the number of the entries $t_k$ equal to $\infty$ is less than $q$.
(In particular, the $n$-tuple is $1$-finite if all entries 
are finite.)
If the tuple is not $q$-finite, we call it 
{\em $q$-infinite}.
\ed

We can now establish 
{\em countable $q$-regularity}
defined in 
Definition~\ref{reg-def}
%\S\ref{reg-sec}
as a consequence of the $q$-finiteness of the tower multitype:

\bp
\Label{q-finite-reg}
Let $S\subset\C^{n+1}$ be a 
(not necessarily pseudoconvex)
smooth hypersurface
whose 
{\em tower multitype is $q$-finite at every point}.
Then $S$ is countably $q$-regular 
in the sense of Definition~\ref{reg-def}
where the ``strata'' $S_k$ 
can be chosen to be the level sets
of the tower multitype function $\6T$.
%In particular, $S$ satisfies property $(P_q)$
%by Proposition~\ref{reg-pq}.
\ep

\bpf
Since the tower multitype of $S$ is 
$q$-finite at every point,
$S$ splits into the countable disjoint union of the $\6T$-level sets
%$$
\beq
\Label{decomp}
S = \bigcup_{(t_1,\ldots,t_n) 
\in (\N_{\ge2}\cup\{\infty\})^n} 
\{p : \6T(p) = (t_1,\ldots,t_n) \}.
\eeq
%$$
By
Corollary~\ref{l-closed},
each level set of $\6T$
is locally closed,
and by the Second structure property (Proposition~\ref{second-prop}), it
 is locally contained in the zero set
 $$
%\beq\Label{M}
M=\{f_1=\ldots=f_l=0\},
%\eeq
$$
where $\{f_1,\ldots,f_l\}$ is an associated 
set of functions
of a tower 
$H^{10}S = E_0\supset \ldots \supset E_m$
on an open subset of $S$.

In view of the $q$-finiteness assumption
and Remark~\ref{th-dep} (1),
we must have
$$\rk E_m = \rk H^{10}S - m=n-m,$$
which implies $\rk E_m<q$.
On the other hand, by the First structure property (Proposition~\ref{first-structure}),
$M$ is a CR submanifold of $S$
satisfying
$$H^{10}M\cap K^{10}\subset E_m,$$
where $K^{10}$ is the distribution of the kernels of the Levi form of $S$, hence
$$
\dim_\C(H^{10}_x\cap K^{10}_x)
\le \rk E_m
<q$$
 for all $x\in M$.
Thus, the decomposition \eqref{decomp}
satisfies all conditions of Definition~\ref{reg-def}
and the proof is complete.
\epf

\bpf[Proof of Proposition~\ref{reg-conclude}]
In view of Defintion~\ref{q-finite-def},
it is clear that Proposition~\ref{q-finite-reg}
contains a restatement of
Proposition~\ref{reg-conclude}
whose proof is therefore completed.
\epf

\bpf[Proof and Theorem~\ref{main}]
Theorem~\ref{main} makes the same assumptions
as Proposition~\ref{reg-conclude} plus pseudoconvexity.
Then Theorem~\ref{main} follows by combining
Proposition~\ref{reg-conclude} with Proposition~\ref{reg-pq},
whose proof was given at the end of \S\ref{reg-sec}.
\epf

\section{Points of $q$-infinite tower multitype}
\Label{inf}
We shall now proceed
towards proving Theorem~\ref{reduce},
 assuming by contradiction
 that \eqref{t<q} is violated, i.e.\
that the multitype at a point $q\in S$ is $q$-infinite,
where $S\subset\cn$ (or $S\subset\C^{n+1}$ whenever
such choice of $n$ is more conevenient) is a pseudoconvex
smooth real hypersurface as before.

\subsection{Subbundles of inifinite Levi type}
%To obtain sufficient conditions
%for
%{\em finite tower multitype},
We thus turn to study points where the multitype is 
{\em $q$-infinite}
% and prove, in particular,
%that such points have infinite regular type.
and begin by establishing the existence
of 
{\em subbundles of infinite Levi type}
along which  the Levi form of the hypersurface
vanishes 
together with its arbitrarily high derivatives
at the given point.
As before, we write $L\in E$ instead of $L\in \G(E)$ for brevity.
%Recall that
% the {\em Levi type of a subbundle 
%$E\subset H^{10}S$ at a point $p\in S$} is defined by
%$$
%c(E):=\min\{
%t\ge 2:
%\exists
%L^t, \ldots, L^1\in E\cup\1{E},
%L^t \cdots L^3 \th([L^2,L^1])(p)\ne0
%\}
%\in \N_{\ge2}\cup \{\infty\},
%$$
%where $\th$ is any complex contact form
%in a neighborhood of $p$ in $S$,
%and $c(E)$ does not depend on the choice of $\th$.
For the reader's convenience,
we restate the definition of the {\em Levi type of a subbundle}
(Definition~\ref{types} (2)):
%\cite{K72,Bl81,D93,HY21,HY23}:

\bd
\Label{flat-sub}
Let $S\subset\C^{n+1}$ be a smooth real hypersurface
with a complex contact form $\th$
in a neighborhood of a point $p\in S$.
\ben
\item
 The {\em Levi type of a complex subbundle 
$E\subset H^{10}S$ at $p$} is 
defined by
$$
c(E,p):=\min\{
t\ge 2:
\exists
L^t, \ldots, L^1\in E\cup\1{E},
L^t \cdots L^3 \th([L^2,L^1])(p)\ne0
\}
\in \N_{\ge2}\cup \{\infty\}.
$$
%and is indepenendent on the choice of $\th$.
\item
If the set inside the minimum above is empty,
i.e. $c(E,p)=\infty$,
the subbundle $E$
is said to be
{\em of infinite Levi type at $p$}, i.e.\
if for any $t\ge 2$ and any smooth vector fields
$L^t, \ldots, L^1\in E\cup\1{E}$,
one has
\beq
\Label{Lt}
L^t \cdots L^3 \th([L^2,L^1])(p)
=0.
\eeq
\een
\ed
Since $\th$ is unique up to a multiple with a nonzero function,
this definition is independent of $\th$.

\br
%Kohn called in \cite{K72}
%the type of $S\in \C^2$ the number $c(H^{10}S,p)-1$.
The invariant $c(E,p)$ in Definition~\ref{flat-sub} slightly
differs from
 the original definition of Bloom
 in that we 
differentiate the Levi form 
evaluated along vector fields
instead of differentiating its trace as
it is done in
 \cite[Definition~2.9]{Bl81}.
However, in their recent work
 Huang-Yin~\cite{HY23}
proved that both versions --- with or without the trace --- 
yield the same values of $c(E,p)$ for 
{\em pseudoconvex hypersurfaces}.
\er

Subbundles  of infinite Levi type arise naturally 
at points of $q$-infinite tower multitype
(see Definition~\ref{q-finite-def}):

\bp
\Label{inf-flat}
Let $S\subset \C^{n+1}$ 
be a smooth real hypersurface a complex contact form $\th$.
Assume that a point $p\in S$ has
{\em $q$-infinite tower multitype $\6T(p)$},
 $U\subset S$ is
an open neighborhood of $p$, and
\beq
\Label{tower-flat}
H^{10}U = E_0\supset \ldots \supset E_m,
%\quad m<n,
\eeq
a tower
 on $U$
  whose multi-order equals the multitype $\6T(p)$.
Then the subbundle $E_m$
satisfies $\rk E_m\ge q$ and
 is of infinite Levi type at $p$.
%for any $t\ge 2$ and any vector fields 
%$L^t, \ldots, L^1\in E_{m}\oplus \1{E}_{m}$,
%one has
%$$
%L^t \cdots L^3 \th([L^2,L^1])(p)
%=0.
%$$
\ep

\bpf
In view of Definition~\ref{tower-def},
the multi-order of the given tower has the form
$$(t_1,\ldots,t_m,\infty,\ldots,\infty)
\in
(\N_{\ge2})^m\times (\{\infty\})^{n-m}$$
with $\rk E_m = n-m.$
Since $\6T(p)$ is $q$-infinite,
Definition~\ref{q-finite-def}
implies $\rk E_m\ge q$.

Next assume by contradiction that
\beq
\Label{nonzero}
L^t \cdots L^3 \th([L^2,L^1])(p)\ne 0
\eeq
for some $t\ge 2$ and some choice of vector fields
$L^t, \ldots, L^1\in E_{m}\cup \1{E}_{m}$.
Then, by repeating the arguments of the proof of Proposition~\ref{second-prop},
we reach a contradiction
constructing another tower
on a neighborhood of $p$ in $S$
of a lexicographically smaller multi-order
$$
(t_1,\ldots,t_m,t_{m+1},\ldots,t_n)<
(t_1,\ldots,t_m,\infty,\ldots,\infty).
$$
%
%
%
%Splitting $L^t$ into components in $E_m$ and $\1{E}_{m}$
%and taking a nonzero term, we may assume
%either $L^t\in E_m$ or $L^t\in \1E_m$
%in \eqref{nonzero}.
%
%Since the multi-order of the given tower equals $\6T(p)$
The contradiction shows that \eqref{nonzero} cannot hold,
hence completing the proof.
\epf

\bpf[Proof of Theorem~\ref{reduce} under conditions \rm{(3)} or \rm{(3$'$)}]
Conditions (3$'$) and hence (3)
on the finiteness of the Levi types are the easiest
to use, since the Levi type is directly related
to the definitions of special subbundles
and towers (Definitions~\ref{special-sub0} and~\ref{tower-def}).

Indeed, assume by contradiction that the conclusion of Theorem~\ref{reduce} fails,
i.e.\
the tower multitype at $p$ is $q$-infinite.
Then taking any tower 
$H^{10}U = E_0\supset \ldots \supset E_m$
realizing the multitype,
it follows
by Proposition~\ref{inf-flat}
that $E_m$
is a subbundle of rank $\ge q$
and infinite Levi type.
Furthermore, 
by Definitions~\ref{special-sub0} and~\ref{tower-def},
$E_m$ is a special subbundle.
Hence both conditions (3$'$) and (3)
are violated, leading to the desired contradiction.
\epf

\section{Formal orbits}
\Label{form-orb}

\subsection{Real formal orbit of a subbundle}
\Label{E-orb}
When working with
vector fields given by convergent power series,
identities \eqref{Lt} for a fixed point $p$
provide some information along 
the so-called Nagano leaf, see e.g.\ \cite{BERbook}.
In the smooth case considered here, we instead have to consider
the so-called {\em formal orbit}.
We use some basic geometric theory of formal power series
that is reviewed in the Appendix \S\ref{formal-theory}.

For every complex vector field $L\in \C T\R^m$
and $p\in \R^m$, consider its
 $\infty$-jet represented by the formal Taylor series
$$
j_p^\infty L \in \C[[x-p]]\otimes \R^m,
\quad
x=(x_1,\ldots,x_m),
$$
%represented by its formal Taylor series at $p$,
where $\C[[x-p]]\otimes \R^m$
consists of complex vector fields
$$
\sum_{j=1}^m a_j(x)\d_{x_j}
=
\sum_{j,\a} a_{j,\a} (x-p)^\a \d_{x_j},
\quad
a_{j,\a}\in \C,
$$
 in $\R^m$
with formal power series coefficients,
where $\a\in \N_{\ge0}^m$ are multi-indices.

Let  $S\subset\R^m$ be a smooth submanifold
and 
 $E\subset \C TS$
a smooth complex subbundle.
To construct the 
real formal orbit of $E$ in $\R^m$,
we have to work with formal vector fields in $\R^m$.
Consider the space $j^\infty_pE$ of $\infty$-jets
at $p$
of all complex vector fields in
% (a neighborhood of $p$ in) 
$\R^m$ 
whose restrictions to $S$ are in $E$:
$$
j_p^\infty E := \{j_p^\infty L: L\in \C T\R^m, L|_S\in E\}
\subset \C[[x-p]]\otimes \R^m.
$$
To obtain {\em real} formal vector fields, consider the real part
$$
\Re j_p^\infty E := \{ \Re L : L \in j_p^\infty E \}
\subset \R[[x-p]]\otimes \R^m.
$$
Note that since $E$ is a complex subbundle,
we have the equality $\Re j_p^\infty E = \Im j_p^\infty E$.
Finally, consider the generated Lie algebra
%$$
%\5g_E(p) \subset \R[[x-p]]\otimes \C
%$$
%to be the complex Lie algebra generated by $j_p^\infty E$,
%i.e.\
%by all Taylor series of smooth complex vector fields 
%whose restrictions to $S$ are in the subbundle $E$.
$$
\5g^\R_E(p):= {\sf Lie} (\Re j^\infty_p E)
\subset
\R[[x-p]]\otimes \R^m.
$$

%
%In our definition of tower, 
%we consider vector fields in subbundles along with their conjugates.
%To define relevant orbits, we consider subbundles of $\C TS$
%that are invariant under conjugation.
%\bl
%Let $S\subset\R^m$ be a smooth submanifold
%and $E\subset \C TS$ a smooth subbundle
%that is invariant under conjugation.
%Then for $p\in S$, the associated Lie algebra
%$\5g_E(p)$ is the complexification of a real Lie algebra
%$\5g^\R_E(p) \subset \R[[x-p]]\otimes \R$,
%i.e.\ $\5g_E(p) = \5g^\R_E(p) \oplus i \5g^\R_E(p)$.
%\el
%
%\bpf
%Since $E$ is invariant under conjugation,
%so is $j^\infty_pE$ and hence the generated Lie algebra $\5g_E(p)$.
%Then $\5g^\R_E(p):=\5g_E(p) \cap \R[[x]]\otimes \R$
%is the real Lie subalgebra satisfying the desired conclusion.
%\epf

%Consider the Lie algebra
%$\5g_{E}\subset \C[[z,\zeta]]$

Now we can use Proposition~\ref{orb-exist}
and the following definition
of the orbit of a Lie algebra of formal vector fields
from the Appendix in the following definition:
\bd\Label{for-orb}
Let $S\subset\R^n$ be a smooth submanifold
and $E\subset \C TS$ a smooth complex subbundle.
The {\em real formal orbit} $\6O^\R_E(p)$ of $E$ at $p\in S$
is defined to be the orbit of the Lie subalgebra
$\5g^\R_E(p)$ of the Lie algebra
of all real formal vector fields in $\R^m$.
\ed

Recall from Proposition~\ref{orb-exist}
that $O=\6O^\R_E(p)\subset\R^n$
is a {\em real formal submanifold},
i.e.\
defined by an ideal $I(O)$
in the ring of formal power series $\R[[x-p]]$,
with $\dim O =\dim \5g^\R_E(p)$,
such that any formal vector field $L\in \5g^\R_E(p)$
is tangent to $O$.
The tangency of $L$ to $O$ is defined as
preservation of the ideal $I(O)$, i.e.\ $L(I(O))\subset I(O)$.
%Also, by the same proposition,
By Remark~\ref{tangent-eq}, 
the tangent space to $O$ at $p$ satisfies
$$
T_p O = \5g^\R_E(p),
$$
where 
the tangent space at $p$ is defined
by 
$$
T_p O = \{df_p=0 : f\in I(O)\},
$$
see Definition~\ref{t-sp} in the Appendix.
Note that the tangency to
$O$ for all vector fields in $\5g^\R_E(p)$
is equivalent to that for all vector fields in 
$\Re j^\infty_pE = \Im j^\infty_pE$,
hence also to the tangency for all vector fields in $j^\infty_pE$.

The orbit $\6O^\R_E(p)$ can be compared
with the $\infty$-jet $j^\infty_pS$ of $S$, which is
the formal submanifold of $\R^m$ defined by 
the ideal $I(j^\infty_pS)\subset\R[[x-p]]$ consisting of all Taylor series at $p$
of smooth functions in $\R^m$ vanishing on $S$.
We have the following inclusion relation between
these two formal submanifolds:

\bl\Label{cont}
Let $S\subset\R^n$ be a smooth submanifold
and $E\subset \C TS$ a smooth complex subbundle.
Then the formal orbit $\6O^\R_E(p)$ at $p\in S$ is 
contained in the $\infty$-jet $j^\infty_p S$
in the sense that $I(\6O^\R_E(p))\supset I(j^\infty_pS)$.
\el

\bpf
Since formal vector fields $L\in j^\infty_pE$
are tangent to the real formal submanifold $O':=j^\infty_pS$,
i.e.\ $L(I(O'))\subset I(O')$,
so are vector fields from the real part $\Re j^\infty_pE$
and hence from the entire Lie algebra $\5g^\R_E(p)$.
Then the last statement in Proposition~\ref{orb-exist}
implies that $O'$ contains the orbit $\6O^\R_E(p)$
of $\5g^\R_E(p)$ as desired.
\epf

\subsection{Orbits of subbundles of infinite Levi type}

We shall use the following basic property of formal orbits:
\bl\Label{orb-vanish}
Let $S\subset\R^m$ be a smooth submanifold through $0$
and $E\subset \C TS$ a smooth complex subbundle.
Let $f$ be a complex smooth function in 
a neighborhood of $0$ in $S$ such that
\beq\Label{f-vanish}
L^t\ldots L^1 f(0) = 0
\eeq
for all $t\ge 0$ and all vector fields 
$L^t, \ldots, L^1 \in E$.
Then the Taylor series $j^\infty_0 f$ of $f$ at $0$ vanishes
along the formal orbit $O=\6O^\R_E(0)$ of $E$ at $0$,
i.e.\  $j^\infty_0 f\in \C I(O)$.
\el

\bpf
Let $\phi\colon (\R^q,0)\to (\R^m,0)$
be a formal parametrization of $O$,
see Definition~\ref{par} in the Appendix.
It follows from the definition of the formal orbit
that there exists formal vector fields
$X_1,\ldots,X_q$ in $\R^q$ that span $T_0\R^q$
and such that each push-forward $\phi_* X_j$
is the Taylor series at $0$ of a vector field in the Lie algebra 
$\5g^\R_E(p)$.
Therefore \eqref{f-vanish} implies 
\beq\Label{x-van}
X_{j_1} \ldots X_{j_k}  (\phi^* f)(0)=0
\eeq
for all $k\ge 0$ and all indices $j_1,\ldots,j_k\in \{1,\ldots,q\}$.
Replacing $X_j$ by their linear combinations over the ring of formal power series, we may assume $X_j = \d_{x_j}$.
But then \eqref{x-van} implies that the Taylor series 
$j^\infty_0(\phi^* f)$ vanishes.
Since $\phi$ is a parametrization of $O$,
this is equivalent to $j^\infty_0f\in \C I(O)$ as claimed.
\epf

Let $S\subset \cn$ be a smooth real hypersurface
with a complex contact form $\th$.
As customary, we identify $\cn$ with $\R^{2n}$
by means of the correspondence $x+iy \mapsto (x,y)$
for $x,y\in\R^n$.
As immediate corollary of Lemma~\ref{orb-vanish}, we obtain:
\bc\Label{levi-orb}
Assume that 
a smooth complex subbundle $E\subset H^{10}S$ is 
of infinite Levi type at $p$
(in the sense of Definition~\ref{flat-sub}).
Then for all
% $p\in S$
%and 
complex vector fields $L^2,L^1$ in a neighborhood of $p$ in $\cn$
with $L^j|_S \in E\cup\1E$, $j=1,2$,
the Taylor series $j^\infty_p\th([L^2, L^1])$ vanishes 
along the formal orbit $\6O^\R_E(p)$, i.e.\
\[
%\beq\Label{o-v}
j^\infty_p\th([L^2, L^1]) \in \C I(\6O^\R_E(p)).
%\eeq
\]
\ec

\bpf
Indeed, the statement follows by applying 
Lemma~\ref{orb-vanish} to $f=\th([L^2, L^1])$.
\epf

\subsection{Consequences of pseudoconvexity}
Recall that a real hypersurface $S\subset\cn$ is
{\em pseudoconvex}
if for every $p\in S$, a purely imaginary 
contact form $\th$ can be chosen 
in a neighborhood of $p$ in $S$ such that
$\th([L,\1L])\ge 0$ holds
for all vector fields $L\in H^{10}S$
in a neighborhood of $p$ in $S$.
In the sequel we assume $\th$
to be smoothly extended to a neighborhood of $p$ in $\cn$
whenever necessary.

We make use of the following well-known lemma:
\bl[positivity of the lowest weight component]
\Label{pos-comp}
Let $f\ge 0$ be a nonnegative smooth function
in a neighborhood of $0$ in $\R^m$.
Fix a collection of positive weights $\mu_j>0$, $j=1,\ldots,m$,
and number $k\ge0$,
and consider a decomposition 
$$
f=f_k + f_{>k},
$$
where $f_k$ is weighted homogeneous of degree $k$,
while the Taylor expansion of $f_{>k}$ at $0$ consists of terms
of weight $>k$.
Then $f_k\ge 0$.
\el

\bpf
The conclusion is obtained by taking the limit of 
$$
t^{-k}f(t^{\mu_1}x_1,\ldots, t^{\mu_m}x_m) \ge0
$$
as $t \to 0$ for each fixed $(x_1,\ldots,x_m)\in\R^m$.
\epf

A well-known result from Linear Algebra
for a positive semi-definite hermitian form $h(\cdot,\cdot)$
asserts that $h(L,L)=0$ for some $L$ implies
$h(L,L')=0$ for any $L'$.
We need the following formal refinement
of this result for vector fields along a formal submanifold:
\bp
\Label{l-l'}
Let $S\subset\cn$ be a pseudoconvex smooth real hypersurface
through $0$
with a contact form $\th$
(extended to a neighborhood of $0$ in $\cn$)
 and $O\subset j^\infty_0 S$ 
a real formal submanifold.
Let $L,L'$ be smooth complex vector fields 
in a neighborhood of $0$ in $\cn$ with 
$L|_S, L'|_S \in H^{10}S$.
Then 
%$$
\beq
\Label{ll0'}
j^\infty_0 \th([L,\1L])\in \C I(O) 
\implies  
j^\infty_0 \th([L,\1{L'}])\in \C I(O).
\eeq
%$$
\ep

\bpf
Since $S$ is pseudoconvex and the conclusion
is independent of the choice of the contact form $\th$,
we can choose $\th$ such that 
$\th([L'',\1{L''}])\ge 0$ for all $L''\in H^{10}S$
in a neighborhood of $0$.
Then
 for all $c\in \C$,
\beq\Label{ll'}
\th([L+cL',\1{L+cL'}])
=
\th([L,\1{L}])
+ 2\Re(\1c\th([L,\1{L'}]))
+c\1c \th([{L'},\1{L'}])
\ge 0
\eeq
holds on $S$ in a neighborhood of $0$.
Since $O\subset j^\infty_p S$,
there exists a smooth map germ 
$\g\colon (\R^q,0)\to (S,0)$ such that
$j^\infty_0\g$ is a formal parametrization of $O$,
see \S\ref{mfd-ideals}.
Taking pullbacks of \eqref{ll'} under $\g$
we obtain that
\beq
\Label{pullbacks}
\g^*\th([L,\1{L}])
+ 2\Re(\1c \g^*\th([L,\1{L'}]))
+c\1c \g^* \th([{L'},\1{L'}])
\ge 0
\eeq
holds in a neighborhood of $0$ in $\R^q$.

The 
right-hand side 
%conclusion 
of the implication
\eqref{ll0'}
%of the proposition
 is equivalent
to $\g^*\th([L,\1{L'}])$ vanishing of infininte order at $0$.
Assume by contradiction that this order is finite and equal
to some $s\in\N_{\ge0}$.
Assign the weight $1$ to $(t_1,\ldots,t_q)\in\R^q$ and 
 the weight $\mu = s+1$ to $(\Re c, \Im c)$.
Then 
$2\Re(\1c \g^*\th([L,\1{L'}]))$ 
has a nonzero lowest weight term 
\beq
\Label{nonzero1}
2\Re(\1c \g^*\th([L,\1{L'}]))_{s+\mu}
\not\equiv 0
\eeq
 of weight $s+\mu$,
while all terms of the expansion of 
$c\1c \g^* \th([{L'},\1{L'}])$ at $0$ have weight $\ge 2\mu>s+\mu$.
On the other hand, since $\g$ is a parametrization of $O$,
 the first term $\g^*\th([L,\1{L}])$
vanishes of infinite order when 
$j^\infty_p \th([L,\1L])\in \C I(O)$.
Therefore, the lowest weight term in \eqref{pullbacks}
is $2\Re(\1c \g^*\th([L,\1{L'}]))_{s+\mu}$,
which satisfies
$$
2\Re(\1c \g^*\th([L,\1{L'}]))_{s+\mu}\ge 0
$$
in view of Lemma~\ref{pos-comp}.
Since $c$ is arbitrary, the left-hand side must vanish,
which contradicts \eqref{nonzero1}, completing the proof.
\epf

\subsection{Identities modulo formal ideals}
As before we work in the ring of $\C[[x]]$
of formal power series 
with complex coefficients
in $x\in \R^m$. 
For an ideal $I\subset \C[[x]]$ 
adopt the standard notation
$$
f=f' \mod I \iff f-f'\in I,
\quad f,f'\in \C[[x]].
$$
A similar notation will be used 
for formal vector fields
$$
L=L' \mod I \iff L-L'\in I \cdot \C[[x]]\otimes \C^m,
\quad L,L\in \C[[x]]\otimes \C^m.
$$
Here 
$I \cdot \C[[x]]\otimes \C^m$
consists of complex formal vector fields
with coefficients in $I$,
while
the notation modulo $I$
is shorter.
Similar notation will be adopted for ideals
in the ring of smooth functions of their germs.

We next study implications of the conclusion of 
Proposition~\ref{l-l'}.
When $L$ is a smooth vector field tangent to $S$
with $\th(L)_p=0$ at a point $p\in S$, then $L_p \in \C H_p S$.
We need a formal version of this conclusion
where $\th(L)$ vanishes along a formal submanifold:

\bl
\Label{th-h}
Let $S\subset\R^m$ be a smooth submanifold
through $0$ with a nonzero complex $1$-form $\th$
(smoothly extended to a neighborhood of $0$ in $\R^m$),
so that $\th=0$ defines a corank $1$ complex subbundle $H\subset \C TS$.
Let $O\subset j^\infty_0 S$
be a real formal submanifold and $L$ a
smooth real vector field
in a neighborhood of $0$ in $\R^m$ that is tangent to $S$.
Then 
$$j^\infty_0 \th(L)=0 \mod \C I(O) \implies j^\infty_0 L \in j^\infty_0 H 
\mod \C I(O).
$$
\el

\bpf
Choose a smooth frame $e_1,\ldots,e_s$ of $H$
and
complete it to a smooth frame $e_1,\ldots,e_{s+1}$ of $\C TS$
with $\th(e_{s+1})=1$,
that we may assume to be
extended smoothly to a neighborhood of $0$ in $\R^m$.
Since $L$ is tangent to $S$, we have
$$
L=\sum_{j=1}^{s+1} a_j e_j
\mod \C I(S)
$$
in a fixed neighborhood of $0$ in $\R^m$,
where $\C I(S)$ is the ideal of complex smooth functions
in that neighborhood vanishing on $S$.
Since $H\subset \C TS$ is defined by $\th=0$,
one has 
$$\th(L)= a_{s+1} \mod  \C I(S).$$
%where $f\ne0$ is a nonzero smooth function.
By definition, the inclusion $O\subset j^\infty_0 S$ 
is equivalent to
$j^\infty_0 I(S) \subset I(O)$, and then
$$j^\infty_0 \th(L)=0 \mod \C I(O) 
\implies 
j^\infty_0 a_{s+1}=0 \mod \C I(O).
$$
The right-hand side implies 
$$
j^\infty_0 L = j^\infty_0 L'  \mod \C I(O),
\quad
L':= \sum_{j=1}^{s} a_j e_j,
$$
which proves the desired conclusion
since $L'|_S\in H$.
\epf

%
%For a real smooth pseudoconvex hypersurface $S\subset\cn$
%as before and $p\in S$, we write
%$j^\infty_pS\subset\cn$ for the real formal submanifold
%defined by the Taylor series of all functions in the ideal $I(S)$
%of smooth real functions vanishing on $S$.
%We also write $D_E$ for smooth sections in a bundle $E$
%and $j^\infty_p D_E$ for
%the space of their Taylor series at $p$.

As per Definition~\ref{formal-tangent} in the Appendix,
for a formal submanifold $O\subset\R^m$,
we write $D_O\subset \R[[x]]\otimes\R^m$
for the submodule of real formal vector fields $L$
satisfying $L(I(O))\subset I(O)$, so that
$$
\C D_O  = \{f+ig: f,g\in D_O\}  =
\{ L\in \C[[x]]\otimes\R^m : L(\C I(O)) \subset \C I(O)) \}  \subset \C[[x]]\otimes\R^m
$$
holds
for the complexification of $D_O$.
We summarize the main properties
of subbundles of infinite Levi type
(Definition~\ref{flat-sub})
 inside 
 {\em pseudoconvex} hypersurfaces
in the following statement:

\bc
\Label{orb-tan}
Let $S\subset\cn$ be a  pseudoconvex smooth real hypersurface
with a smooth complex subbundle $E\subset H^{10}S$
that is of infinite Levi type at $p\in S$.
Denote by $O=O^\R_{E}(p)$ the real formal orbit.
Then the following holds:
\ben
\item
commutators with Taylor series at $p$ of sections in $E\cup \1E$ 
leave invariant
the space of Taylor series at $p$ of sections in $\C HS$
 modulo the ideal of the orbit, i.e.\
$$
[j^\infty_p (E\oplus\1E) , 
j^\infty_p \C HS]
\subset j^\infty_p \C HS
\mod \C I(O);
$$
\item
commutators with
vector fields in $D_O$ tangent to $O$ leave invariant
the space of Taylor series at $p$ of sections in $\C HS$
 modulo the ideal of the orbit, i.e.\
$$
[\C D_O , 
j^\infty_p \C HS]
\subset j^\infty_p \C HS
\mod \C I(O);
%D_{\6O_{E}(p)}
%\subset
%\6D^{10}_{j^\infty_pS} \oplus \1{\6D^{10}_{j^\infty_pS}}
%\mod I(\6O_{E}(p)),
$$
%i.e.\ every formal vector field tangent to
%$O$ is contained, modulo the ideal of $O$,
%in the sum of formal $(1,0)$ and $(0,1)$ vector fields
%tangent to $j^\infty_p S$.
\item
the orbit $O$ is complex-tangential to 
$j^\infty_p S$ in the sense that
$$
\C D_O \subset  j^\infty_p \C HS \mod \C I(O).
$$
\een
\ec

\bpf
Choose a complex contact form $\th$ on $S$ as before,
smoothly extended to a neighborhood of $p$ in $\cn$.
To show (1), 
consider 
formal vector fields in 
${j^\infty_p E}$
and 
${j^\infty_p H^{10}S}$
which are
 Taylor series at $p$ of smooth vector fields $L, L'$
 whose restrictions to $S$ are in $E$ and $H^{10}S$ respectively.
Then, since $E$ is of infinite Levi type, by Corollary~\ref{levi-orb} 
$$
j^\infty_p\th([L, \1L]) \in \C I(\6O),
$$
and, in view of the pseudoconvexity, Proposition~\ref{l-l'}
implies
\beq
\Label{j-v}
j^\infty_p \th([L,\1{L'}])\in \C I(O).
\eeq
On the other hand, 
\eqref{j-v} holds also for $L'$ with $L'|_S\in H^{01}S$,
since $H^{10}S$ is an involutive distribution in $\C TS$
by the integrability of the CR structure.
Then by Lemma~\ref{th-h},
we obtain
$$
[j^\infty_p E , 
j^\infty_p \C HS]
\subset j^\infty_p \C HS
\mod \C I(O),
$$
from which (1) follows,
since the right-hand side is invariant under conjugation.

To show (2), observe that 
(1) along with Jacobi identity implies that
$$
[L , 
j^\infty_p \C HS]
\subset j^\infty_p \C HS
\mod \C I(O)
$$
holds for any $L$
which is an iterated Lie bracket
of formal vector fields in $j^\infty_p(E\oplus \1E)$.
Since $O=\6O_E^\R(p)$
is the orbit of the Lie algebra $\5g=\5g^\R_E$
spanned by the real parts $\Re L =\frac12( L+\1L)$
of iterated Lie brackets,
it follows that 
$$
[\5g, 
j^\infty_p \C HS]
\subset j^\infty_p \C HS
\mod \C I(O).
$$
Finally, by Lemma~\ref{orb-tangent} in the Appendix,
$D_O=\5g \mod I(O)$,
proving (2).

Now, since $E\oplus\1E\subset \C HS$,
by
repeatedly using (1),
we conclude that
any iterated Lie bracket of vector fields
in $j^\infty_p(E\oplus\1E)$
is contained in $j^\infty_p \C HS$ modulo $\C I(O)$.
As before, using the relation 
$D_O=\5g \mod I(O)$
completes the proof of (3).
\epf

\bpf[Proof of Theorem~\ref{reduce} under conditions \rm{(2)} or \rm{(2$'$)}]
Given the
completed Proof of Theorem~\ref{reduce} under condition (3$'$)
at the end of \S\ref{inf},
it suffices to establish the condition (3$'$).

Assume by contradiction condition (3$'$) fails,
i.e.\
there exists
a neighborhood $U$ of $p$ in $S$
and
a 
special (see Definition~\ref{special-sub0})
subbundle
 $E\subset H^{10}U$ of rank $\ge q$,
 which is 
of infinite Levi type at $p$
and 
{\em special} (see Definition~\ref{special-sub0}).
By Corollary~\ref{orb-tan} (3),
the real formal orbit
$O=\6O^\R_E(p)$
is complex-tangential.
In particular,
its tangent space at $p$
is complex-tangential,
i.e.\ 
\beq
\Label{tp}
T_pO \subset H_pS.
\eeq

On the other hand, 
consider
an iterated commutator
$
[L^m,\ldots,[L^{2},L^1]\ldots]
$
with 
$
L^m,\ldots,L^1\in \G(E)\cup\G(\1E)
$
as in condition (2$'$).
Then by the construction of the orbit 
$O=\6O^\R_E(p)$,
extending $L^j$ smoothly to a neighborhood in $\cn$,
it is easy to see that
their formal power series
$j^\infty_p L^j$
along with their respective iterated commutator
$[j^\infty_p L^m,\ldots,[j^\infty_p L^{2},j^\infty_pL^1]\ldots]$
is tangent to $\6O$.
Finally, since the value of the iterated commutator at $p$
only depends on the Taylor series, 
the inclusion \eqref{tp}
implies that both conditions (2$'$) and (2)
are violated, leading to the desired contradiction.
\epf

\section{Special subbundles and a formal Huang-Yin condition}
In their detailed study of different notions of types
for subbundles $E\subset H^{10}S$,
Huang and Yin~\cite{HY21}
formulated an important condition
on the Nagano leaf 
of a weighted truncation of $E$
to be CR.
We shall need an analogous condition
on the real formal orbit
that we call the {\em formal Huang-Yin condition}.

Recall from \S\ref{form-cr} in the Appendix
%(Definition~\ref{CR-formal}
that a real formal submanifold $O\subset\cn$
given by formal power series $f(z,\1z)$
can be identified with its complexification
given by the complexified series $f(z,\zeta)$.
I.e.\ $O$ is identified with
 a complex formal submanifold of $\cn\times \cn$
given by the complex manifold ideal $\C I(O)$
that we identify with its complexification 
$\{f(z,\zeta) : f(z,\bar z)\in \C I(O)\}$
in the ring $\C[[z,\zeta]]$, $z,\zeta\in \cn$,
that we still denote by $\C I(O)$ in a slight abuse of notation.
Then
the $(1,0)$ tangent space $H^{10}_0 O \subset \cn$
is defined by
$$
H^{10}_0 = 
\Big\{
\xi\ : \frac{\d f}{\d z}(0) \xi = 0 \text{ for all } f(z,\zeta)\in \C I(O)
\Big\},
$$
and
the space $D^{10}_O$ of formal vector fields
by 
$$
D^{10}_O = \{ L = \sum_j a_j(z,\zeta) \d_{z_j}  : L(\C I(O)) \subset \C I(O)\},
$$
and following 
Della Salla-Juhlin-Lamel~\cite{DJL12},
$O$ is {\em said to be CR}
if 
$$
H^{10}_0O \subset D^{10}_O(0)
$$
holds,
where $D^{10}_O(0) = \{L_0 : L\in D^{10}_O\}$
is the space of evaluations at $0$.

\bd
\Label{hy}
Let $S\subset\cn$ be a smooth real hypersurface.
We say that a smooth complex subbundle $E\subset H^{10}S$
satisfies the {\em formal Huang-Yin condition} at $p\in S$
if the real formal orbit $O=\6O^\R_E(p)$ is CR,
i.e.\ $H^{10}_pO \subset D^{10}_O(p)$.
\ed

It follows from 
Lemma~\ref{djl} 
in the Appendix
(a consequence of a result from \cite{DJL12})) 
that if $O=\6O^\R_E(p)$ is CR,
the ideal $\C[[z]]\cap \C I(O) \subset \C[[z]]$ 
is a complex manifold ideal
defining the 
{\em intrinsic complexification} of $O$
as a complex formal submanifold of $\cn$.
In other words, the intrinsic complexification of $O$
is the complex formal submanifold
defined by the ideal $\C[[z]]\cap \C I(O)$.

We shall prove that in our situation
the formal Huang-Yin condition always holds.
For this, recall that our subbundle $E$
is a part of the tower
$$
H^{10}S = E_0\supset \ldots \ldots E_m = E,
$$
where all subbundles $E_k\subset H^{10}$
are {\em special}
in the sense of Defintion~\ref{special-sub0},
i.e.\ defined by linearly independent $\th$-dual forms.
%(in fact, under somewhat more general assumption).
%
%
%
%$E$ is in some special way by means of a
%defining system (see Definition~\ref{tower-def}, part (3))
%$$
%(L_1,\ldots,L_l; f_{l+1},\ldots, f_m),
%$$
%and $E=E_m \subset H^{10}S $ is defined by
% simultaneous vanishing of the $1$-forms
%\beq
%\Label{wf-wan}
%E=
%\{
%\om_{L_1;\th}= \ldots= \om_{L_l;\th}
%= \d f_{l+1}= \ldots = \d f_m = 0
%\}.
%\eeq
%This description motivates the following definition:
%\bd[Special subbundle]
%\Label{special-sub}
%We call
%a complex subbundle $E\subset H^{10}M$ {\em special}
%if there exists a collection
%$(L_1,\ldots,L_l; f_{l+1},\ldots, f_m)$,
%where $L_j\in \C HS$ are smooth vector fields 
%and $f_k$ are smooth real functions,
%such that \eqref{syst-non} holds
%and $E$ is given by \eqref{wf-wan}.
%\ed

%As always, the complex contact form $\th$
%is unique up to a nonzero function multiple,
%implying that our definition is independent of 
%the choice of $\th$.

In fact, we establish
the desired formal Huang-Yin property
under somewhat more general assumptions
that may be of independent interest.
The difference is that we assume
that $E\subset H^{10}S$ is defined by 
any collection of linearly independent 
$1$-forms
that are either contractions of the Levi form
or differential of arbitrary smooth functions
--- an assumption that always 
holds for special subbundles.

\bt
\Label{CR}
Let $S\subset\cn$ be a  pseudoconvex smooth real hypersurface
with a complex contact form $\th$ and
a smooth complex subbundle $E\subset H^{10}S$
such that the following holds:
\ben
\item 
$E$ is of infinite Levi type at $p\in S$ 
in the sense of Definition~\ref{flat-sub};
\item
$E$ is defined by 
\beq
\Label{wf-wan}
E=
\{
\om_{L_1;\th}= \ldots= \om_{L_l;\th}
= \d f_{l+1}= \ldots = \d f_m = 0
\},
\eeq
where
$L_j\in \C HS$ are smooth vector fields 
and $f_k$ are smooth real functions
satisfying
$$
%\beq
\Label{syst-non}
\om_{L_1;\th}\wedge \ldots \wedge \om_{L_l;\th}
\wedge \d f_{l+1}\wedge \ldots\wedge \d f_m \ne 0,
\quad
\om_{L_j;\th}(L_x) = \th([L,L_j])(x),
%\eeq
\quad x\in S.
$$
\een
Then the real formal orbit of $E$ at $p$ is CR,
i.e.\ $E$ satisfies the formal Huang-Yin condition at $p$.
\et

\bpf
As before we write $O=\6O^\R_E(p)$ for the real formal orbit of $E$ at $p$. Since $T_pO = D_O(p)$, Corollary~\ref{orb-tan}  implies that $T_pO\subset H_pS$ and the Levi form of $S$ at $p$
vanishes on $T_pO$. In particular, all $1$-forms $\om_{L_j;\th}$
vanish on $T_pO$. 

On the other hand, 
it follows from \eqref{wf-wan}
that
the subbundle $E\oplus \1E$ is tangent to the CR submanifold
$$
M=\{x: f_{s}(x)=f_s(p), \; s=l+1,\ldots,m\}\subset S.
$$
Then, also any iterated Lie bracket of vector fields in $E\oplus \1E$
is tangent to $M$, hence
the formal orbit $O$ is tangent to $j^\infty_p M$,
implying $T_pO\subset T_pM$, and therefore
$$
H^{10}_p O \subset \{\d f_{l+1}= \ldots = \d f_m = 0\}.
$$
Together with vanishing $\om_{L_j;\th}$ obtained above,
this implies
$$
H^{10}_p O \subset \{
\om_{L_1;\th}= \ldots= \om_{L_l;\th}
=\d f_{l+1}= \ldots = \d f_m = 0\},
$$
hence
 $H^{10}_p O\subset E_p$ in view of \eqref{wf-wan}.

On the other hand, 
since formal vector fields in $j^\infty_p E$ are tangent to $O$
and $E\subset H^{10}S$,
 $\infty$-jets of $(1,0)$ vector fields 
 $L$ with $L|_S\in E$
% whose restrictions to $S$
%are contained in $E$ 
are contained in
the module $D^{10}_O$ of formal $(1,0)$ vector fields
tangent to $O$
(recall \eqref{d10} in the Appendix for the definition).
Since any vector in $E_p$ can be extended to such vector field,
it follows that $E_p\subset D^{10}_O(p)$.
Combining with above inclusion $H^{10}_p O\subset E_p$, we obtain
$H^{10}_p O\subset  D^{10}_O(p)$ as claimed
in view of Definition~\ref{hy}.
%Since the opposite inclusion 
%$H^{10}_p O\supset  D^{10}_O(p)$
%always holds, 
\epf

As observed above, assumption (2)
in Theorem~\ref{CR}
is always satisfied 
for special subbundles.
Hence we obtain as a special case of Theorem~\ref{CR}:

\bc
\Label{hy-cor}
Let $S\subset\cn$ be a  pseudoconvex smooth real hypersurface
with a complex contact form $\th$.
Then any special subbundle $E\subset H^{10}S$
of infinite Levi type at $p$
is CR, i.e.\
satisfies the formal Huang-Yin condition at $p$.
\ec

We conclude by providing
an example,
that is communicated to the author
by Lamel and Lebl,
demonstrating
that the formal Huang-Yin condition
is crucial to guarantee
that the 
{\em formal complex orbit},
see Definition~\ref{E-contact} (1),
is regular, i.e.\ a formal submanifold:

\bex[Lamel-Lebl]
\Label{l-l}
Consider the following 
real non-CR submanifold 
of codimension $2$ in $\C^3$:
$$
M := \{ 
(z,w_1,w_2)\in \C^3 : 
w_1=|z|^4, \; w_2=|z|^6
\}.
$$
It is clear that the complex tangent space
has complex dimension $2$
 at points with $z=0$ 
and $1$ otherwise,
hence $M$ is not CR.
Taking any collection $\6L$ of real-analytic vector fields
tangent to $M$ that everywhere span the tangent space,
we can realize $M$ as the  (formal) real  orbit of $\6L$.
It is also clear that
the intrinsic complexification of $M$,
i.e.\
the minimal complex-analytic subvariety
$
V\subset \C^3
$
containing $M$,
is
given by 
$w_1^3-w_2^2=0$,
which is singular at $0$.
Hence the formal complex orbit $V$ of $\6L$
is also singular.
\eex

\section{A formal variant of a result of Diederich-Fornaess}
 In \cite{DiF78}
Diederich and Fornaess proved a fundamental result about
 the equivalence of 
 finite type and zero holomorphic dimension
for 
{\em real-analytic} hypersurfaces,
which together 
with Kohn's celebrated
work \cite{K79}
yields subelliptic estimates under these conditions.
The most difficult part of their proof 
is Proposition~3.1 in \cite{DiF78}
roughly stating 
that intrinsic complexifications
of smooth CR submanifolds of pseudoconvex hypersurfaces $S$,
along which
the Levi form vanishes, are tangent to $S$ of infinite order
at generic points. 
This result has been generalized by Bedford-Fornaess~\cite{BF81}
and more recently by Pinton-Zampieri~\cite{PZ15}.

To complete the proof of Theorem~\ref{reg},
we shall need a 
{\em formal variant} of this important result,
that may be of independent interest.
Instead of a CR submanifold in $S$,
in our setup
we have a formal CR submanifold ---
the formal orbit of a special subbundle
of infinite Levi type at a given point,
where we use Corollary~\ref{hy-cor}
to ensure the CR condition,
which in turn, 
ensures the existence of
the intrinsic complexification 
as a formal complex submanifold,
see Corollary~\ref{cr-cx}
and
Definition~\ref{i-c}
in the Appendix.
We shall denote by 
$J$ the standard complex structure of $\C^{n}$.

\bt\Label{df}
Let $S\subset\cn$ be a
{\em pseudoconvex} 
 smooth
real hypersuface through $0$,
and let
$O, V\subset \cn$ be respectively
%{\em formal CR submanifold of CR-dimension $q\ge 1$} 
a real and a complex formal submanifold with $O\subset V$.
Assume:
\ben
\item
$O$ is generic in $V$ in the sense that
$T_0 O + J T_0 O = T_0 V$;
\item
$O$ is contained in $S$ in the formal sense, i.e.
$
O\subset j_0^\infty S
$;
\item
$O$ is complex-tangential to $j_0^\infty S$ in the sense
that $D_O \subset j_0^\infty HS \mod I(O)$;
\item
$O$ is of finite commiutator type, i.e.\
%$D^{10}_O\oplus \1{D^{10}_O}$ generates
%as Lie algebra all tangent vectors at $0$, i.e. 
${\sf Lie}(D^{10}_O\oplus \1{D^{10}_O})(0) = \C T_0 O.$
%under Lie brackets from $D^{10}_O\oplus \1{D^{10}_O}$.
%$$
%\quad D_O\subset j^\infty_0 HS
%$$
%i.e.\ $O$ is a ``complex tangential'' submanifold of $S$
%in the formal sense.
\een
Then $V$
 is tangent to $S$ of infinite order
 or, equivalently,
 formally contained in $S$, i.e. $V\subset j^\infty_0 S$.
%\item
%the order of contact of $M$ with $V_k$ at $0$ is at least $k$;
%\een
In particular, $S$ is of infinite regular type at $0$.
%
%Assume that $M$ is of finite Bloom-Graham type,
%i.e.\ iterated brackets of (formal) tangent $(1,0)$ vector fields
%on $M$ and their conjugates
%span the tangent space at $0$.
%Then $S$ contains the intrinsic complexification of $M$.
\et

We shall apply Theorem~\ref{df}
at points of 
{\em $q$-infinite tower multitype}
(Definition~\ref{q-finite-def}),
where $O$ is a real formal orbit of 
the corresponding special
 subbundle $E\subset H^{10}S$
 of infinite Levi type
with $\rk E\ge q$
given by Proposition~\ref{inf-flat},
and $V$ is the intrinsic complexification of $O$.
%See Corollary~\ref{i-c-prop} for the genericity condition.

The original proof in \cite{DiF78} (and in \cite{PZ15})
is based on a careful analysis of the tangent vector fields
at {\em generic points} and cannot be transferred to the
present formal setup, where all power series are defined
at a fixed point $p$
and cannot be evaluated generically.
The approach in \cite{BF81} uses attached analytic disks
to the smooth CR submanifold in $S$,
a tool that is not available for formal CR submanifolds.

Thus we need to use a different approach
based on another set of new tools
such as 
{\em jets along formal submanifolds}
in \S\ref{jets-form}
and 
{\em supertangent vector fields}
in \S\ref{sup}.
The proof
of Theorem~\ref{df}
will be completed 
at the end of \S\ref{sup} below.

\subsection{Assumptions of Theorem~\ref{df} hold
with $\dim_\C V\ge q$
at points of $q$-infinite tower multitype}

Before embarking on our proof of Theorem~\ref{df},
we show how its assumptions arise in our context
for a point of $q$-infinite tower multitype
in the sense of Definition~\ref{q-finite-def}.

\bp
\Label{thm-ass}
Let $S$ be a smooth pseudoconvex hypersurface through $0$.
Assume that the tower multitype of $S$ at $0$ is $q$-infinite.
Then there exists a smooth complex subbundle $E\subset H^{10}S$
 in a neighborhood of $0$ with $\rk E\ge q$,
 which is special (see Definition~\ref{special-sub0}),
whose formal orbit $O$ is CR and 
such that $O$ and its intrinsic complexification $V$
satisfy the assumptions (1)--(4) of Theorem~\ref{df},
and one has $\dim_\C V\ge q$.
\ep

\bpf
We begin by applying Proposition~\ref{inf-flat}
establishing the existence of the tower \eqref{tower-flat}
in a neighborhood of $0$ whose last subbundle
$E=E_m$ is of rank at least $q$ and
of infinite Levi type at $p$ (in the sense of Definition~\ref{flat-sub}).
Since $E_m$ is a subbundle from a tower,
it is {\em special} in the sense of Definition~\ref{special-sub0}
as a direct consequence of the definitions.

Let $O$ denote the (real) formal orbit $O$ of $E$,
see Definition~\ref{for-orb}.
Since the subbundle $E$ is {\em special}
in the sense of Definition~\ref{special-sub0},
and hence $O$ is CR
by Corollary~\ref{hy-cor}.
In particular, the {\em intrinsic complexification} $V$ of $O$
is defined as a complex formal submanifold of $\cn$,
see Definition~\ref{i-c} in the Appendix.

We now establish assumptions (1)--(4) of Theorem~\ref{df}.
Assumption (1) is a general property
of the intrinsic complexification, see 
Corollary~\ref{i-c-prop}.
Assumption (2)
follows from Lemma~\ref{cont}.
Assumption (3) follows from 
Corollary~\ref{orb-tan}, part (3),
by taking real parts.
Finally, to show assumption (4),
recal from the definition of the orbit, Definition~\ref{for-orb},
that, in particular,
\beq
\Label{t0O}
T_0 O = \5g^\R_E(0) = {\sf Lie}(\Re j^\infty_0 E),
\eeq
where 
\beq
\Label{L}
\Re j^\infty_0 E = \{ j^\infty_0 \Re L : L \in \C T\C^{n}, L|_S\in E\}.
\eeq
Let $L$ be as in \eqref{L} and
split it into $(1,0)$ and $(0,1)$ parts:
$$
L=L^{10} + L^{01}, 
\quad L^{10} \in H^{10} \cn,
\quad
L^{01}\in H^{01}\cn.
$$
Since $E\subset H^{10}S$, we must have 
$$
L|_S\in E \implies L^{01}|_S =0.
$$
Then
$$
L' := L^{10} + \1{L^{01}} \in H^{10}\cn,
$$
is a $(1,0)$ vector field with
$L'|_S \in E$,
hence $j^\infty_0 L'$ is tangent to the orbit $O$,
which implies
$j^\infty_0 L' \in D^{10}_O$, see \eqref{d10} in the Appendix
for this notation.
Since
$$
2\Re L=
L+\1L =  L^{10} + L^{01} + \1{L^{10}} + \1{L^{01}} =  L'+\1{L'}, 
$$
it follows that $j^\infty_0 \Re L \in D^{10}_O \oplus D^{01}_O$.
Finally, taking arbitrary $L$ in \eqref{L} and using \eqref{t0O},
it follows that 
$$
T_0 O \subset {\sf Lie}(D^{10}_O\oplus D^{01}_O)(0).
$$ 
Since the right-hand side is a complex vector space,
this proves assumption (4) as desired.

Finally, since $O$ is the orbit of the complex subbundle $E$ of rank at least $q$,
both $O$ and its intrinsic complexification $V$
have real dimensions at least $2q$,
in particular, $\dim_\C V\ge q$ as stated.
\epf

\subsection{Approximation by maps into hypersurfaces}
\Label{approx-hyp}
For the proof of Theorem~\ref{df}
we shall need to compare
the Levi form of the hypersurface $S$
with the complex hessian of its defining function $r$
along the complex formal submanifold $V$.
Our approach is to ``approximate'' a formal parametrization of $V$,
or more generally, a formal map $A$ into $V$
with a smooth map into $S$. 
We obtain a fine control of this approximation 
modulo the ideal of the pull-back $(j^\infty_0 r)\circ A$
as follows:

\bl
\Label{param-approx}
Let $Y\subset \R^m$ be a formal hypersurface
defined by a principal ideal $I(Y)$.
Let 
$R\in \R[[x]]$, $x=(x_1,\ldots,x_m)$,
be a generator of $I(Y)$ with
$dR(0)\ne 0$,
and $A\colon (\R^s,0)\to (\R^m,0)$ 
any formal power series map.
Then there is an ``approximating'' formal power series map
$\2 A\colon (\R^s,0)\to (\R^m,0)$ satisfying
\ben
\item
$R \circ \2A =0$;
\item
$\2A-A = 0 \mod (R\circ A)$.
\een
In case $Y=j^\infty_0 S$ for a smooth real hypersurface $S\subset \R^m$ through $0$, the approximating map can be chosen of the form
$\2A = j^\infty_0 a$, where 
$a\colon (\R^s,0)\to (S,0)$ is a (germ of a) smooth map.
\el

\bpf
By the implicit function theorem for formal power series,
we may assume the given generator $R$ of the principal ideal $I(Y)$
to be of the form
$$
R(x', x_m) = x_m - \Phi(x'),
\quad x'=(x_1,\ldots,x_{m-1}),
\quad
\Phi\in \R[[x']],
\quad
(x',x_m)\in \R^{m-1}\times \R,
$$
cf. \eqref{spec-gen}.
Splitting $A$ into components
$$
A=(A',A_m),
\quad
A' \colon (\R^s,0)\to (\R^{m-1},0),
\quad
A_m \colon (\R^s,0)\to (\R,0),
$$
we obtain (1) by setting
$
\2A := (A', \Phi\circ A')
$.
Then 
$$
A - \2A = (0, A_m - \Phi\circ A' )=
(0, R\circ A)
$$
proving (2).

When $Y=j^\infty_0 S$ for a smooth real hypersurface $S\subset \R^m$ through $0$,
% with local defining function $r$, 
a local defining function for $S$ near $0$ can be chosen of the form $r(x',x_m)=x_m - \phi(x')$ with $j^\infty_0 \phi = \Phi$.
Then $R:=j^\infty_0 r$ yields a generator of $I(j^\infty_0 S)$.
By Borel's theorem, there exists a germ of a smooth map
$a'\colon (\R^s,0)\to (\R^{m-1},0)$ with $j^\infty_0 a' = A'$.
Then the smooth germ $a:= (a', \phi\circ a')$ satisfies the desired properties.
\epf

\subsection{Approximation by tangent vector fields}
Similar to the approximation of maps
in \S\ref{approx-hyp},
we approximate a formal $(1,0)$ vector field $L$
with a smooth $(1,0)$ vector field tangent to 
the given real hypersurface.
We control the approximation
by the ideal obtained from applying $L$
to the Taylor series $j^\infty_0 r$ of the defining function $r$.
%A key step
%is to compare the Levi form of a hypersurface $S\subset \cn$
%along a complex tangent vector field $L'$
%with the complex Hessian of the defining function
%along another vector field $\2L$ tangent to the
%complexified formal orbit in $\cn$.
%The following lemma provides a simple
%construction of approximating complex tangent 
%vector fields in the sense of ideals.

\bl\Label{L-approx}
Let $R\in\C[[z,\z]]$, $z=(z_1,\ldots,z_n)$,
be a real formal power series with $dR(0)\ne 0$
and  
$L=\sum a_j(z,\z) \d_{z_j}$ 
 a $(1,0)$ formal  vector field
in $\cn$.
Then there exists an (approximating) $(1,0)$ formal vector field 
$\2L$ in $\cn$ satisfying:
\ben
\item $\2L R= 0$;
\item $\2L=L \mod (LR)$,
where $(LR)$ 
is the principal ideal generated by $LR$.
\een
In case $R=j^\infty_0 r$
for a smooth defining function $r$ 
of a real hypersurface in $\cn$, 
the approximating vector field $\2L$ 
can be chosen of the form
$\2L = j^\infty_0 \2l$, where 
$\2l$ is a germ at $0$ of a smooth $(1,0)$ vector field in $\cn$
with $\2l r=0$.

%In the special case when $m=2l$ and $L$ is $(1,0)$
%in $\C^l\cong \R^{2l}$, 
%the approximating vector field $L'$ satisfying (1) and (2) can also be assumed $(1,0)$.
\el

\bpf
Since
$R$ is real and 
$dR(0)\ne 0$, after possible reordering coordinates, we may assume
$R_{z_n}(0)\ne 0$. 
Consider a $(1,0)$ formal vector field $\2L = L + f\d_{z_n}$,
where $f\in\C[[z,\z]]$ is a formal power series. Then
$$
\2LR = LR + fr_{z_n}=0 \iff f = -\frac{LR}{R_{z_n}}.
$$
Since $R_{z_n}(0)\ne 0$, $\frac1{R_{z_n}}$ is a formal power series.
Hence taking $f = -\frac{LR}{R_{z_n}}$, we obtain
the formal vector field $\2L= L + f\d_{z_n}$
satisyfing both conditions (1) and (2) as desired.
%For the last statement for $L$ of type $(1,0)$ in $\C^l$,
%repeat the same argument with $r_{z_l}(0)\ne 0$ and
%$\2L = L - \frac1{r_{z_l}}\d_{x_l}$.

When $R=j^\infty_0 r$
for a smooth defining function $r$ 
of a real hypersurface in $\cn$,
by Borel's theorem, there exists
a germ $l$ of a smooth vector field in $\cn$ at $0$ 
with $L=j^\infty_0 l$.
Then setting 
$$
\2l := l - \frac{lr}{r_{z_n}} \d_{z_n}
$$
yields the desired properties.
\epf

\section{Jets along formal submanifolds}
\Label{jets-form}
\subsection{$(k,X,Y)$-equivalence and jets}
Our next tool needed for proving Theorem~\ref{df}
provided by a formal variant
of jets of functions restricted to submanifolds
along smaller submanifolds. 
When $X\subset Y\subset \R^m$ are smooth submanifolds through 
$p\in X$,
we can identify germs at $p$ of smooth complex functions $f$ on $Y$ with
germs of smooth complex functions on $\R^m$ modulo
the complexified ideal   $\C I_p(Y)$
of smooth germs at $p$ 
vanishing on $Y$,
and their $k$-jets $j^k_x f$ 
restricted to $X$ as germs at $p$ of functions of $x\in X$
with their equivalence classes modulo the ideal
$
\C I_p(X)^{k+1} + \C I_p(Y)
%, \quad x=(x_1,\ldots,x_m).
$
of smooth function germs at $p$.
This motivates the following formal variant of this concept:

\bd
Let $X\subset Y\subset \R^m$ be
real formal submanifolds given respectively by 
their manifold ideals $ I(Y) \subset I(X)\subset\R[[x]]$.
We say that formal power series 
$$
F,G\in\C[[x]],
\quad
x=(x_1,\ldots,x_m),
$$ 
%have their restrictions to $Y$ $k$-equivalent relative to $X$ or
are 
{\em $(k,X,Y)$-equivalent}
%at $X$ along $Y$
if 
$$
F-G \in \C I(X)^{k+1} + \C I(Y)).
$$
A {\em formal $(k,X,Y)$-jet} is a $(k,X,Y)$-equivalence class of 
formal power series in $\C[[x]]$.
\ed

As a matter of notation, we shall write 
$$
J^k_X \C[[Y]]= \C[[x]]/(\C I(X)^{k+1} + \C I(Y))
$$
for the ring of all $(k,X,Y)$-jets of power series in  $\C[[x]]$
and 
$$
j^k_{X,Y}F \in J^k_X \C[[Y]], 
\quad F\in \C[[x]],
$$
for $(k,X,Y)$-jets of individual formal power series.

A convenient way of representing $(k,X,Y)$-jets
is via joint parametrizations of $X,Y$
as provided by Lemma~\ref{joint-param} in the Appendix,
i.e.\ formal maps
\beq
\Label{A-joint}
A\colon (\R_{t'}^{s'}\times \R_{t''}^{s''},0)\to (\R^m,0)
\eeq
such that $A$ is a parametrization of $Y$
and $A|_{\R^{s'}}$ is a parametrization of $X$.
Then the pullback of $F\in \C[[x]]$ takes the form
$$
F\circ A = 
\sum_{\a'\a''} f_{\a'\a''} (t')^{\a'} (t'')^{\a''},
\quad
(t',t'')\in \R_{t'}^{s'}\times \R_{t''}^{s''},
$$
where $\a',\a''$ are multi-indices,
the pullback of the ideal $\C I(X)^{k+1} + \C I(Y)$ is
$$
A^*(\C I(X)^{k+1} + \C I(Y)) = \C I(\R_{t'}^{s'})^{k+1},
$$
and the pullback $(j^k_{X,Y}F)\circ A$ of the $(k,X,Y)$-jet of $F$
can be represented 
by polynomials of degree $\le k$ in $t''$
with formal power series coefficients in $t'$:
$$
(j^k_{X,Y}F)\circ A = j^k_{\R^{s'}, \R^{s'+s''}} (F\circ A) 
=
\sum_{\a',  |\a''|\le k} f_{\a'\a''} (t')^{\a'} (t'')^{\a''}.
$$

\subsection{Testing vanishing jets by their derivatives}
In proving Theorem~\ref{df}
we shall need a conclusion about 
a $(k,X,Y)$-jet being zero
when suitable lower jets are
zero for all derivatives 
along formal vector fields tangent to $Y$:

\bl
\Label{testing}
Let $X\subset Y\subset \R^m$ be
real formal submanifolds given respectively by 
their manifold ideals $ I(Y) \subset I(X)\subset\R[[x]]$.
Fix integers $0< k' \le k''$
and a formal power series 
$$
F\in I(X)^{k''} + I(Y).
$$
Assume that
\beq
\Label{der}
j^{k''-k'}_{X,Y} (L^{k'}\cdots L^1 F) =0
%L^k\ldots L^1 F \in I(X)^{k'-k+1} + I(Y)
\eeq
for all real formal vector fields $L^{k'},\ldots,L^1\in D_Y$
(i.e. tangent to $Y$).
Then $j^{k''}_{X,Y} F =0$, i.e.\ 
$F \in I(X)^{k''+1} + I(Y)$.
\el

\bpf
We shall use a joint formal parametrization $A$
as in \eqref{A-joint} such that 
$A|_{\R^{s'}}$ is a parametrization of $X$.
Then $F\circ A\in I(\R_{t'}^{s'})^{k''}$
and hence
$$
F\circ A = 
\sum_{\a',  |\a''|\ge k''} f_{\a'\a''} (t')^{\a'} (t'')^{\a''}.
$$
Assume by contradiction $F \notin I(X)^{k''+1} + I(Y)$,
which is equivalent to $F\circ A\notin I(\R_{t'}^{s'})^{k''+1}$.
Then there exists a multi-index $(\a',\a'')$ with
\beq
\Label{a}
f_{\a'\a''}\ne 0, \quad |\a''|=k''.
\eeq
Choosing a basis of $D_Y$ as in \eqref{do} in the Appendix,
we can rewrite \eqref{der} in the equivalent form
\beq
\Label{dtfa}
\d_{t_{j_k}}\cdots \d_{t_{j_1}} (F\circ A) \in I(\R_{t'}^{s'})^{k''-k'+1},
\quad
j_k,\ldots,j_1 \in \{1,\ldots,s'+s''\}.
\eeq
Now let $\a''$ be as in \eqref{a}
and choose indices $j_k,\ldots,j_1$ among the components 
of the multi-index $\a''$, so that
$$
g:=\d_{t_{j_{k'}}}\ldots \d_{t_{j_1}} (f_{\a'\a''} (t')^{\a'} (t'')^{\a''}) \ne 0.
%\quad
%g\in I(\R^{s''})^{k'-k}
$$
Then $g$ is a nontrivial monomial of order $k''-k'$ in $t''$
in the expansion of the left-hand side of \eqref{dtfa},
which hence cannot belong to $I(\R_{t'}^{s'})^{k''-k'+1}$,
a contradiction with \eqref{dtfa}
proving the desired conclusion.
\epf

\subsection{Taylor series of nonnegative functions}
In our proof of Theorem~\ref{df}, 
{\em pseudoconvexity} of the given hypersurface
yields nonegative functions
by means of evaluating the Levi form
along vector fields.
The lemma below captures
the way we use nonnegativity of those functions.
\bl\Label{pos}
Let 
$$
f(x,y)\ge 0,
\quad
(x,y)
\in
\R^m_x\times\R^n_y,
\quad m,n\ge 1,
$$
be a real smooth nonnegative function
in a neighborhood of $0$,
whose Taylor series $j_0^\infty f$ at $0$ satisfies
\beq
\Label{tf-cond}
j_0^\infty f\in I(\R^m_x)^k, \quad k\ge 0.
%(y)^k(z)^s + (y)^{k+1}(z)^{s-1}
%\sum_{k'\ge k,\; k'+s' = k+s}
%(y)^{k'}(z)^{s'},
%\quad k,k',s,s'\ge 0,
% I(\R^{m}_y)^{k} \cap I(\R^{n}_z)^{s},
\eeq
%where $k,k',s,s'$ are integers.
%where as before $I(O)$ denotes the formal ideal of a (formal) submanifold $O$, and
% we are using the standard identifications
%$$
%\R^{m}_y \cong \{0\}\times \R^{m}_y\times\{0\},
%\quad 
%\R^{n}_z \cong \{0\}\times\{0\}\times \R^{n}_z.
%$$
Assume that $k$ is odd.
Then $j_0^\infty f \in I(\R^m_x)^{k+1}$.
% of degree $k$ in $y$ vanish.
%, i.e.\
%$$
%j_0^\infty f\in 
%(y)^{k+1}(z)^s + (y)^k(z)^{s+1}
%.
%I(\R^{m})^{k+1} + I(\R^{n})^{s+1}.
%$$
\el

\bpf
We shall write $(x)$ and $(y)$ for the ideals
generated by all component coordinates of $x$ and $y$ respectively,
in particular, $(y)=I(\R^m_x)$.
Denote by $P$ the (formal) sum
of all terms of $j_0^\infty f$ of degree $k$ in $y$.
%By the assumptions, the Taylor series $j_0^\infty f$ can be written as
%\beq
%\label{Tf-exp}
%j_0^\infty f(x,y,z) = P(x,y,z) \mod 
%(y)^{k+1}(z)^s + (y)^k(z)^{s+1}
%%I(\R^{m})^{k+1} + I(\R^{n})^{s+1}
%,
%\eeq
%where 
Then
\beq\Label{fp}
j^\infty_0 f(x,y)=P(x,y) \mod (y)^{k+1},
\quad
P(x,y) \in \R[[x]][y],
\quad
P\in (y)^k,
\eeq
where $P$
is a homogeneous polynomial of degree $k$ in $y$
with formal power series coefficients in $x$.
It suffices to prove that $P=0$.

Assume by contradiction that $P\ne 0$.
Then we have the integer
$$
b:=\max\{a: P(x,y) \in (x)^a(y)^k  \}  \in \N_{\ge0},
$$
which also equals the minimum order in $x$ of a nonzero coefficient of $P$,
and $P$ can be written as
\beq\Label{P}
P(x,y) = P_0(x,y) \mod (x)^{b+1} (y)^k,
\eeq
where 
\beq
\Label{P0}
P_0\ne0\in \R[x,y],
\eeq 
is a bi-homogeneous polynomial in $(x,y)$
of the bi-degree $(b,k)$.

We restrict $f$ and $P$ to the real curves given by
$$
(x,y) = (At,Bt^{b+1}), 
\quad (A,B)\in \R^{m}_x\times\R^{n}_y,
\quad t\in \R.
$$
Then the bi-homogeneity of $P_0$
along with \eqref{P} imply
$$
P(At,Bt^{b+1}) = P_0(A,B)t^{b+k(b+1)}
\mod (t)^{b+k(b+1)+1},
$$
which, together with \eqref{tf-cond}, implies 
$$
(j_0^\infty f)(At,Bt^{b+1}) = P_0(A,B)t^{b+k(b+1)}
\mod (t)^{(k+1)(b+1)}.
$$
Since 
$$
(j_0^\infty f)(At,Bt^{b+1})$$ 
is the Taylor series at $0$ of
$f(At,Bt^{b+1})$ and $f\ge0$
in a neighborhood of $0$, we must have
\beq\Label{p0-pos}
P_0(A,B)\ge0,
\quad  (A,B)\in \R^{m}_x\times\R^{n}_y.
\eeq

Now we use our assumption that $k$ is odd.
Recall that $P_0(x,y)$ is homogeneous of degree $k$ in $y$.
Then \eqref{p0-pos} is only possible when $P_0=0$,
contradicting \eqref{P0}.
Hence our assumption $P\ne0$ is false and the proof is complete. 
\epf

%\br
%It is clear from the proof that 
%Lemma~\ref{pos} can be generalized
%with $(y,z)\in \R^m\times \R^n$ replaced by any number of
%groups of variables.
%We shall only need it for two groups $y,z$ as stated
%or just for one group $y$, corresponding to a special case of the lemma with
%$f$ independent of $z$ and $s=0$.
%\er

\section{Part I of Proof of Theorem~\ref{df}: Relative contact orders}
In the notation of Theorem~\ref{df},
we write $r$ for a smooth local defining function
in a neighborhood of $0$
of the smooth pseudoconvex real hypersurface $S\subset\cn$,
in particular, $dr\ne 0$.
Further, write $R=j^\infty_0 r$ for the Taylor series of $r$ at $0$.
Then $R$ is a generator of the ideal $I(j^\infty_0 S)$
of the formal hypersurface $j^\infty_0 S$
and the assumption $O\subset j^\infty_0 S$
implies $R\in I(O)$.
%Also $O\subset V$ implies $I(V)\subset I(O)$.

We also use a joint formal parametrization 
\beq\Label{A}
A\colon (\R^{s'}\times \R^{s''},0) \to (\R^{2n},0)
\eeq
as in \eqref{A-joint} such that 
$A$ is a parametrization of $V$ and
$A|_{\R^{s'}}$ is a parametrization of $O$.

\subsection{Relative contact orders along formal submanifolds}
%Let $r$ be a local defining function of $S$
%in a neighborhood of $0$ with $dr\ne 0$.
%Since our goal is to prove that $V$
%is tangent to $S$ of infinite order,
To prove Theorem~\ref{df},
we
assume by contradiction that $V$ is not tangent to $S$
of infinite order at $0$.
Then there exists the integer
\beq\Label{k}
k := 
%\min \{s: j^s_{O,V} Tr \ne 0 \} = 
\min \{s : R \in I(O)^{s} + I(V) \}  \in \N_{\ge0},
\eeq
that we call the {\em relative contact order}
between $S$
and the pair $(O,V)$ as mentioned in \S\ref{new-tech}.
Since $R\in I(O)$ in view of the assumption (2) in Theorem~\ref{df}, 
it follows that $k\ge 1$.
In fact, in the setup of Theorem~\ref{df}, it follows that $k\ge2$:

\bl\Label{kge2}
The assumption (1)--(3) in Theorem~\ref{df}
imply $k\ge 2$.
\el

\bpf
Since $HS$ and hence $j^\infty_0 HS$
is $J$-invariant,
where $J$ is the standard complex structure of $\cn$,
 (2) and (3) in Theorem~\ref{df}
imply
$$
D_V\subset j^\infty_0 HS \mod I(O),
$$ 
hence $LR = 0 \mod I(O)$
or, equivalently, $j^0_{O,V} (LR) =0$
 for all $L\in D_V$.
Then Lemma~\ref{testing} 
with $k'=k''=1$
implies
$R\in I(O)^2+I(V)$ 
and hence $k\ge 2$
as claimed.
\epf

\subsection{Vanishing odd jets of complex hessians}
We next study the complex hessian $\d\bar\d R$
evaluated along vector fields $L$ tangent to $V$.
Since $R\in I(O)^k + I(V)$, it follows that 
$\d\bar\d R(L,\1L) \in I(O)^{k-2} + I(V)$ for $L\in D^{10}_V$.
If $k$ is odd, we show that the lowest order jet of
$\d\bar\d R(L,\1L)$ must in fact vanish:

\bl\Label{odd-zero}
Assume $k$ is odd.
Then the assumption (1)--(3) in Theorem~\ref{df}
imply that 
$$j^{k-2}_{O,V} \d\bar\d R(L,\1L)=0,$$
or, equivalently,
\beq
\Label{rll}
\d\bar\d R(L,\1L) \in I(O)^{k-1} + I(V)
\eeq
for any $L\in D^{10}_V$.
\el

\bpf
Let $A$ be a joint parametrization given by \eqref{A},
so that
we have $R\circ A \in I(\R^{s'})^k$.
Then by Lemma~\ref{param-approx},
there exists an approximation by a smooth map germ
$a\colon (\R^{s'}\times \R^{s''},0) \to (S,0)$
with 
$$
j^\infty_0 a = A \mod I(\R^{s''})^k.
$$
Similarly, applying Lemma~\ref{L-approx} to $L\in D^{10}_V$,
we obtain an approximation by a smooth $(1,0)$ vector field germ
$l$ at $0$ in $\cn$ with 
$$
l|_S \in H^{10}S, \quad j^\infty_0 l = L \mod \C I(O)^{k-1} + \C I(V). 
$$
Using both approximations, we obtain
\beq
\Label{a-A}
j^\infty_0 (\d\bar \d r(l,\1l) \circ a) - \d\bar\d R(L,\1L) \circ A 
\in I(\R^{s'})^{k-1}.
\eeq
%as the left-hand side is real.
Since $S$ is {\em pseudoconvex}, we may assume
 $r$ to be chosen such that
$\d\bar\d r(l,\1l)|_S\ge 0$.
Since 
$$
j^\infty_0 (\d\bar \d r(l,\1l) \circ a) \in I(\R^{s'})^{k-2}
$$
when $k$ is odd, Lemma~\ref{pos}
implies
$j^\infty_0 (\d\bar \d r(l,\1l) \circ a) \in I(\R^{s'})^{k-1}$.
Together with \eqref{a-A}, this yields
$\d\bar\d R(L,\1L) \circ A \in I(\R^{s'})^{k-1}$
that is equivalent to the desired conclusion.
\epf

From Lemma~\ref{odd-zero}
we derive a vanishing property for $(k-2,O,V)$-jets of derivatives
along pairs of $(1,0)$ and $(0,1)$ vector fields:

\bc\Label{mixed}
Under the assumptions of Lemma~\ref{odd-zero},
$$
L^2 \1{L^1} R \in \C I(O)^{k-1} + \C I(V),
\quad L^2, L^1\in D^{10}_V.
$$
\ec

\bpf
Applying Lemma~\ref{odd-zero} along with polarization 
argument for \eqref{rll} yields
$$
\d\bar\d R(L^2,\1{L^1}) \in \C I(O)^{k-1} +  \C I(V),
\quad L^2,L^1\in D^{10}_V.
$$
By the Cartan formula for the exterior derivative applied to $\1\d R$,
$$
L^2 \1\d R(\1{L^1}) - \1{L^1} \; \1\d R(L^2) - \1\d R([L^2, \1{L^1}])
\in \C I(O)^{k-1} +  \C I(V).
$$
From here the conclusion follows,
because
the first term equals $L^2\1{L^1}R$, the second term
vanishes since $\1\d R$ is a $(0,1)$ form,
and the last term belongs to $\C I(O)^{k-1} +  \C I(V)$,
since $R\in I(O)^k + I(V)$ and 
$[L^2, \1{L^1}], J[L^2, \1{L^1}] \in \C D_V$.
\epf

\subsection{Pseudoconvexity implies that the relative contact order $k$ is even}

We used pseudoconvexity in Lemma~\ref{odd-zero}
together with the assumption that $k$ is odd
to obtain the vanishing of jets of the complex hessian of $r$.
We now derive a contradiction
with that assumption, proving that $k$ must be even.

\bl\Label{one-v}
Without any additional assumptions, 
$R\in I(O)^k+I(V)$ implies
$$
L^2 L^1 R \in \C I(O)^{k-1} + \C I(V),
\quad L^2, L^1\in \C D_V,
$$
whenever $L^j\in \C D_O$ for either $j=1$ or $j=2$.
\el

\bpf
Assume the first case $L^1\in \C D_O$.
Then $R\in I(O)^k+I(V)$ implies $L^1R\in I(O)^k+I(V)$,
from which the conclusion follows, 
since 
$$L^2  (\C I(O)^k) \subset  \C I(O)^{k-1},
\quad 
L^2 \C I(V) \subset \C I(V).
$$

Now assume the second case $L^2\in \C D_O$.
Then we write
$$
L^2L^1 R = L^1L^2 R + [L^2,L^1] R
$$
and use
our above argument 
in the first case and the property $[L^2,L^1]\in \C D_V$ 
to conclude that both terms belong to 
$\C I(O)^{k-1} + \C I(V)$,
implying the desired conclusion.
\epf

\bc\Label{must-even}
Pseudoconvexity of $S$ along with assumption (1)--(3) in Theorem~\ref{df}
imply that the relative contact order $k$ must be even.
\ec

\bpf
The assumption (1) in Theorem~\ref{df} that $O$ is generic in $V$
together with Corollary~\ref{dv-split} imply
$$
\C D_O\cap \C D_V + D^{10}_V = \C D_O\cap \C D_V + \1{D^{10}_V} = \C D_V \mod \C I(V),
$$
which in combination with Corollary~\ref{mixed} and 
Lemma~\ref{one-v} yields
$$
L^2 L^1 R \in \C I(O)^{k-1} + \C I(V),
\quad L^2, L^1\in \C D_V,
$$
whenever $k$ is odd, that is assumed in Corollary~\ref{mixed}.
Recall that $k\ge2$ by Lemma~\ref{kge2}.
Finally, applying Lemma~\ref{testing}
for $k'=2$ and $k''=k$, we obtain
$$
R\in I(O)^{k+1} + I(V),
$$
implying that the relative contact order is $\ge k+1$, 
a contradiction with our choice of $k$, see \eqref{k}.
Hence $k$ must be even as claimed.
\epf

\section{Part II of Proof of Theorem~\ref{df}: Supertangent vector fields}
\Label{sup}
\subsection{Definition and examples of complex-supertangent vector fields}

Our next new tool in the proof of Theorem~\ref{df}
is the notion that we call {\em supertangent}
and {\em complex-supertangent} vector fields to a (formal) hypersurface,
 which applied to a defining function of the hypersurface
 do not reduce the relative contact order.
As before, we write $R=j^\infty_o r$,
where $r$ is a smooth local 
defining function of the given real hypersurface $S$
near the origin.

\bd
\Label{super}
In the setup of Theorem~\ref{df},
let $k$ be the relative contact order
defined by \eqref{k}, in particular, $R\in I(O)^k + I(V)$.
A formal vector field $L\in \C D_V$ is 
said to be 
\ben
\item
{\em supertangent} to $S$ (or $j^\infty_0 S$) if 
$L R \in \C I(O)^k + \C I(V)$;
\item {\em complex-supertangent} to $S$ (or $j^\infty_0 S$)
if both $L$ and $JL$ are supertangent.
\een
\ed

As before, $J$ stands for the standard complex structure
of $\cn$.
Note that applying a general $L\in \C D_V$
to $R$ may potentially reduce the relative contact order,
i.e.\ we may only be able to conclude that 
$L R \in \C I(O)^{k-1} + \C I(V)$.
The meaning of the condition (1) in the above definition
is that supertangent vector fields
provide a better partial order of contact in comparison.

%the relative contact order of $LR$ is only known to be at least $(k-1)$.
%Hence, super-tangent vector fields provide 
%directions of derivatives of $R$ of a higher
%than expected relative contact order.

An immediate source of supertangent vector fields
is provided by the following simple lemma:

\bl
\Label{d10-super}
Any vector field $L\in (D^{10}_O\oplus \1{D^{10}_O}) \cap \C D_V$
is complex-supertangent.
\el

\bpf
Given $L$ as in the lemma, it follows that
$$
L, JL \in \C D_O \cap \C D_V,
$$
i.e., $L$ and $JL$ leave both ideals 
$\C I(O)$, $\C I(V)$ invariant.
Since $R\in I(O)^k + I(V)$, it follows that
$
LR, (JL)R \in I(O)^k + I(V)
$
proving the desired conclusion.
\epf

\subsection{Vanishing jets of complex hessians along supertangent vector fields}
In Lemma~\ref{odd-zero}
we showed the vanishing of $(k-2,O,V)$-jets
of complex hessians of $R$ when $k$ is odd.
Now that we established in Corollary~\ref{must-even}
that $k$ must be in fact even, 
we show a similar vanishing property of 
$(k-1,O,V)$-jets along
 {\em supertangent vector fields}:

\bl
\Label{vanish-even}
Assume $k$ is even. Then 
$$
j^{k-1}_{O,V} \d\1\d R(L,\1L)  = 0,
$$
or, equivalently,
$$
\d\1\d R(L,\1L) \in I(O)^k + I(V)
$$
hold
for any {\em supertangent} vector field $L\in D^{10}_V$.
\el

\bpf
The proof is analogous to that of Lemma~\ref{odd-zero}.
As in that proof,
let $A$ be a joint parametrization given by \eqref{A},
we have $R\circ A \in I(\R^{s'})^k$.
Then by Lemma~\ref{param-approx},
there exists an approximation by a smooth map germ
$a\colon (\R^{s'}\times \R^{s''},0) \to (S,0)$
with 
$$
j^\infty_0 a - A = 0 \mod I(\R^{s''})^k.
$$
Since $L$ is assumed to be {\em supertangent},
applying Lemma~\ref{L-approx} to $L\in D^{10}_V$,
we obtain an approximation by a smooth $(1,0)$ vector field germ
$l$ at $0$ in $\cn$ with
$$
l|_S \in HS, \quad j^\infty_0 l - L \in \C I(O)^{k} + \C I(V). 
$$
Using both approximations, we conclude
\beq
\Label{a-A'}
j^\infty_0 (\d\bar \d r(l,\1l) \circ a) - \d\bar\d R(L,\1L) \circ A 
\in I(\R^{s'})^{k}.
\eeq
%as the left-hand side is real.

By the Cartan formula for the exterior derivative,
%$$
\beq
\Label{ddr}
\d\1\d R(L,\1L) = L \1\d R(\1L) - \1L\; \1\d R(L) - \1\d R([L,\1L]).
\eeq
%$$
Since $L$ 
%and therefore $\1L$ are 
is supertangent and $R$ is real, 
$$
\1\d R(\1L) =\1{LR}  \in \C I(O)^k+ \C I(V),
$$
hence the first term in \eqref{ddr}
 belongs to $\C I(O)^{k-1} + \C I(V)$.
The second term vanishes since $L\in D^{10}_V$,
while the last term again belongs to $\C I(O)^{k-1} + \C I(V)$,
since $[L,\1L], J[L,\1L]\in D_V$. Hence
\beq
\Label{ddRk-1}
\d\1\d R(L,\1L) \in \C I(O)^{k-1} + \C I(V).
\eeq
Since $S$ is {\em pseudoconvex}, we may assume $r$ chosen such that
$\d\bar\d r(l,\1l)|_S\ge 0$,
where $l$ is the approximation of $L$ chosen above.
Then \eqref{a-A'} and \eqref{ddRk-1} yield
$$j^\infty_0 (\d\bar \d r(l,\1l) \circ a) \in I(\R^{s'})^{k-1}.$$
Now, since $k$ is even, $k-1$ is odd and Lemma~\ref{pos}
implies
$j^\infty_0 (\d\bar \d r(l,\1l) \circ a) \in I(\R^{s'})^{k}$.
Together with \eqref{a-A'}, this yields
$\d\bar\d R(L,\1L) \circ A \in I(\R^{s'})^{k}$
that is equivalent to the desired conclusion.
\epf

\subsection{Lie algebra property of complex-supertangent vector fields}

We now derive the key
{\em Lie algebra property
of complex-supertangent vector fields}
among those tangent to both $O$ and $V$:

\bc
\Label{lie-closed}
Under the assumptions of Lemma~\ref{vanish-even},
let $L^2,L^1 \in \C D_O\cap \C D_V$ 
be complex-supertangent vector fields, i.e.\
\beq
\Label{cs-tang}
L^jR, (JL^j)R \in \C I(O)^k + \C I(V), \quad j=1,2.
\eeq
Then their commutator $[L^2,L^1]$ is also complex-supertangent, i.e.\
$$
[L^2,L^1] R, (J[L^2,L^1])R \in  \C I(O)^k + \C I(V).
$$
Thus the complex-supertangent vector fields in
 $\C D_O\cap \C D_V$ form a Lie algebra.
\ec

\bpf
By the Cartan formula for the exterior derivative applied to the $1$-forms $\1\d R$, $\d R$,
\beq
\Label{cart}
\begin{split}
\d\1\d R(L^2,L^1) = L^2 \1\d R(L^1) - L^1\1\d R(L^2) - \1\d R([L^2,L^1]),
\\
\1\d\d R(L^2,L^1) = L^2 \d R(L^1) - L^1 \d R(L^2) - \d R([L^2,L^1]).
\end{split}
\eeq
Every $L^j$, $j=1,2$, can be represented in the standard way
as the sum of its $(1,0)$ and $(0,1)$ components:
\beq
\Label{splitting}
L^j = L^j_{10} + L^j_{01},
\quad 
L^j_{10}=
\frac12(L^j + iJ L^j) \in D^{10}_V, 
\quad 
L^j_{01}=
\frac12(L^j - iJ L^j) \in D^{01}_V = \1{D^{10}_V}
\eeq
and since $R$ is real,
our assumption \eqref{cs-tang} implies that 
both components $L^j_{10}, L^j_{01}$
are supertangent and 
\beq
\Label{dr}
\d R(L^j) = \d R(L^j_{10}) \in \C I(O)^k + \C I(V),
\quad
\1\d R(L^j) = \d R(L^j_{01}) \in \C I(O)^k + \C I(V).
\eeq
Applying Lemma~\ref{vanish-even} along with polarization,
we compute
$$
\d\1\d R(L^2,L^1) = - \1\d\d R(L^2,L^1)
=\d\1\d R(L^2_{10}, L^1_{01}) - \d\1\d R(L^2_{01}, L^1_{10})
\in \C I(O)^k + \C I(V),
$$
i.e.\ left-hand side terms in both equations \eqref{cart}  belong to 
$\C I(O)^k + \C I(V)$.
Applying $L^j$ to the identities in \eqref{dr}
and using our assumption $L^j \in \C D_O \cap \C D_V$ for $j=1,2$,
we compute
$$
 L^2 \1\d R(L^1), L^1\1\d R(L^2), L^2 \d R(L^1), L^1 \d R(L^2)
 \in \C I(O)^k + \C I(V).
$$
Hence \eqref{cart} imply
$$
\1\d R([L^2,L^1]), \d R([L^2,L^1]) \in \C I(O)^k + \C I(V).
$$
Splitting $[L^2,L^1]$ into $(1,0)$ and $(0,1)$ components 
similar to \eqref{splitting}, this proves the desired conclusion.
\epf

\subsection{Lie algebras generated by complex-supertangent vector fields}

The assumption (4) in Theorem~\ref{df}
is formulated for vector fields in $D^{10}_O\oplus \1{D^{10}_O}$.
To use supertangent vector fields,
we show a stronger conclusion involving
vector fields that are both in $D^{10}_O\oplus \1{D^{10}_O}$
and tangent to $V$:

\bl
\Label{lie-span}
The assumtion (4) in Theorem~\ref{df}, i.e.\ 
\beq
\Label{t0}
{\sf Lie}(D^{10}_O\oplus \1{D^{10}_O})(0) = \C T_0 O
\eeq
implies
$$
\C D_O =
{\sf span}_{\C[[z,\z]]}
{\sf Lie}((D^{10}_O\oplus \1{D^{10}_O}) \cap \C D_V) 
+ \C I(O) \cdot \C[[z,\z]]\otimes \C^{2n} ,
$$
where as before ${\sf Lie}$ denotes the generated Lie algebra
with respect to the Lie brackets and ${\sf span}_{\C[[z,\z]]}$
is the span over the ring $\C[[z,\z]]$.
\el

\bpf
Since $O$ is contained in the formal submanifold $V$,
it follows from Lemma~\ref{ov-ap} in the Appendix that
every vector field $L\in D^{10}_O\oplus \1{D^{10}_O}$
can be approximated modulo $\C I(O)$ by a vector field 
$\2L\in (D^{10}_O\oplus \1{D^{10}_O}) \cap \C D_V$,
i.e.\ 
$$
\2L = L \mod \C I(O).
$$
By \eqref{t0},
$\C T_0 O$ is spanned
by a finite collection of vector field values $L^1_0,\ldots,L^s_0$ at $0$,
where each $L^j$ is an iterated Lie bracket of vector fields 
$L\in D^{10}_O\oplus \1{D^{10}_O}$.
%Since each $L^j$ is an iterated Lie bracket of vector fields in 
%$D^{10}_O\oplus \1{D^{10}_O}$.
Since all such $L$ leave $\C I(O)$ invariant,
replacing them by their approximations $\2L$ modulo $\C I(O)$,
we obtain an approximation $\2L^j$ of 
each iterated bracket $L^j$ 
such that
\beq
\Label{lj-ap}
\2L^j = L^j \mod \C I(O),
\eeq
where each $\2L^j$ is an iterated Lie bracket of vector fields in 
$(D^{10}_O\oplus \1{D^{10}_O}) \cap \C D_V$,
i.e.\ 
$$
\2L^j \in 
{\sf Lie}
((D^{10}_O\oplus \1{D^{10}_O}) \cap \C D_V).
$$
By \eqref{lj-ap}, $\2L^j_0 = L^j_0$ for the values at $0$,
hence the $\2L^j_0$ span $\C T_0 O$ for $j=1,\ldots,s$,
and the desired conclusion follow from Lemma~\ref{span}
in the appendix.
\epf

We now use the assumtion (4) in Theorem~\ref{df}
to show that 
{\em all vector fields} in $\C D_O\cap \C D_V$ are in fact 
{\em complex-supertangent}.

\bc\Label{every-super}
In the setup of Theorem~\ref{df},
let $k$ be the relative contact order
defined by \eqref{k}
Assume $k$ is even.
Then the assumptions (1) and (4) in Theorem~\ref{df}
imply that all vector fields in $\C D_V$
are supertangent.
\ec

\bpf
By Lemma~\ref{d10-super}, 
vector fields in $(D^{10}_O\oplus \1{D^{10}_O}) \cap \C D_V$
are complex-supertangent.
In view of Lemma~\ref{lie-span},
iterated commutators of vector fields in 
$(D^{10}_O\oplus \1{D^{10}_O}) \cap \C D_V$
span the tangent space $\C T_0 O$.
Evidently, any such iterated Lie bracket is in $\C D_V$.
By Corollary~\ref{lie-closed},
all these iterated commutators are complex-supertangent.
Hence the span of all values at $0$ of 
complex-supertangent vector fields
contains the tangent space $\C T_0 O$.
But since $JL$ is supertangent whenever $L$ is complex-supertangent,
the genericity assumption (1) in Theorem~\ref{df}, i.e.\ 
$T_0 O + J T_0 O = T_0 V$ implies
that the span of all values at $0$ of 
supertangent vector fields equals $\C T_0 V$.

Now apply Lemma~\ref{span} for $V$
to conclude that $\C D_V$ is spanned 
by supertangent vector fields over the ring $\C[[z,\z]]$
modulo $\C I(V)$.
Then the desired conclusion follows directly from the definition
of the supertangent vector fields - Definition~\ref{super}.
\epf

\bpf[Proof of Theorem~\ref{df}]
The proof is completed
by obtaining the desired contradiction
from assuming the finiteness of $k$.
%\epf
%
%\bc
%\Label{k-inf}
%Under the assumptions of Theorem~\ref{df},
% the relative contact order $k$ must be infinite.
%\ec

As before, we assume by contradiction that $k$ is finite.
By Corollary~\ref{must-even}, $k$ must be even.
Then Corollary~\ref{every-super} implies
 that
$LR \in \C I(O)^k +\C I(V)$ for all $L\in \C D_V$.
Now we can apply Lemma~\ref{testing}
for $k'=1$, $k''=k$,
to conclude that $R\in I(O)^{k+1} + I(V)$,
which contradicts the choice of $k$.
\epf

%
%\subsection{Contact order along formal submanifolds}
%We use $(k,X,Y)$-jets for measuring
%contact orders between formal submanifolds
%along a common formal submanifold.
%
%\bd
%\Label{contact-def}
%Let $X,Y,Y'\subset \R^m$
%be formal submanifolds with $\dim Y\le \dim Y'$
%and
% such that $X$
%is contained in both $Y$ and $Y'$.
%The {\em contact order along $X$ between $Y$ and $Y'$}
%is 
%%the maximum $k\ge 1$ such that
%$$
%\sup 
%\{ k\ge 1:
%I(Y') \subset I(X)^k + I(Y)
%\}
%\in \N_{\ge1}\cup \{\infty\}
%.
%$$
%\ed
%
%
%\bc
%In the notation of Theorem~\ref{df}, assume only (1),
%let $r$ be a defining function of $S$ and let $k\ge 3$
%be an odd integer with
%$$
%j^\infty_0 r\in I(O)^k + I(V).
%$$
%Then for every $L\in D^{10}_V$,
%the lowest order jet of the complex hessian vanishes:
%$$
%j^{k-2}_{O,V} \d\bar\d r(L,\1L) =0.
%$$
%\ec
%
%
%\bpf
%
%\epf

\section{Completion of the proof of Theorem~\ref{reduce}}
\Label{pf-reg}
Here we complete the proof Theorem~\ref{reduce}
and thereby
all statements that still remain to be proved.

Recall that we already proved 
 Theorem~\ref{reduce} under condition (2$'$) and (3$'$)
 at the end of \S\ref{inf} and \S\ref{form-orb} respectively.
 Since (2) and (3) imply respectively (2$'$) and (3$'$),
 Thus it remains to prove the theorem under condition (1$'$).

Assume by contradiction that the conclusion of Theorem~\ref{reduce} fails,
i.e.\
the tower multitype at $p$ is $q$-infinite.
Then Proposition~\ref{thm-ass}
implies that the assumptions
of Theorem~\ref{df}
are satisfied
for the pair $(O,V)=(\6O^\R_E(p),\6O^\C_E(p)$
consisting of
respectively
 the real and complex formal orbits of $E$ at $p$.
Now by Theorem~\ref{df},
the complex orbit $V=\6O^\C_E(p)$
is tangent to $S$ of infinite order at $p$,
hence violates both conditions (1$'$) and (1), leading to the desired contradiction.

%
%Let $\Om$ be a pseudoconvex bounded domain 
%of finite regular type at every boudary pont
% as in Theorem~\ref{reg}.
%Our goal 
%is to prove that $\d\Om$ satisfies property (P).
%In view of Corollary~\ref{fin-m-type},
%it suffices to show that $\d\Om$
%has a {\em finite tower multitype} at every point.
%
%Assume by contradiction, 
%there exists a point $p\in \d\Om$
%of infinite tower multitype.
%Then by
%Proposition~\ref{thm-ass},
%there exist formal submanifolds 
%$O\subset V\subset \cn$ of positive dimension
%satisfying the assumptions of Theorem~\ref{df}.
%Finally, by Theorem~\ref{df},
%the complex formal submanifold $V$
%is tangent to $S$ of infinite order.
%Since $V$ has positive dimension,
%this contradicts the finite regular type at $p$,
%which completes the proof.

\section{Appendix: formal submanifolds}
\Label{formal-theory}
Here we summarize basic results on formal submanifolds
used throughout the paper.
%\cite{T92}

\subsection{Preliminaries}
We consider a field $\K$ that is either $\R$ or $\C$,
and
the ring 
$$
\K[[x]]=\K[[x_1,\ldots,x_m]],
\quad
x=(x_1,\ldots,x_m)\in \K^m,$$
of 
{\em formal power series} over $\K$
in $m$ indeterminates.
To every smooth complex function germ $f(x)$ at $p\in \K^m$ is associated
its formal Taylor series at $p$, given by
$$
Tf(x) = T_pf(x) = \sum_{\a\in \N^m}  
\frac1{\a!} \frac{\d f}{\d x^\a}(p) (x-p)^\a
\in \K[[x-p]],
$$
where $\a=(\a_1,\ldots,\a_m)$ are multiindices and $\N=\{0,1,\ldots\}$.
It follows directly from the Leibnitz rule that 
$T_p$ is a ring homomorphism
from the ring of germs of smooth functions
at $p$ into the ring of formal power series in $x-p$. 
We shall also use the jet notation and identify $T_pf$ with
the infinite order jet
$$
j^\infty_p f = T_p f.
$$

Similarly, 
by taking Taylor series of the coefficients,
to every germ at $p$ of a smooth $\K$-valued $q$-form 
$$
\om = \sum_J dx_J,
\quad
J=(j_1,\ldots,j_q), 
\quad
dx_J = dx_1\wedge\cdots\wedge dx_q
$$
is associated
its Taylor series
$$
T\om = j^\infty_p \om = \sum_J (T\om_J) dx_J \in 
\K[[x-p]]\otimes \L^q (\K^m)^*.
$$
Analogously, by taking Taylor series of the coefficients,
to every germ at $p$ of a smooth $\K$-valued vector field
$L=\sum_j a_j \d_{x_j}$ is associated
its Taylor series
$$
TL =  j^\infty_p L = \sum_j (Ta_j) \d_{x_j} \in 
\K[[x-p]]\otimes \K^m.
$$

%In the following we write $O$ for a formal submanifold
%to highlight our main use-case where $O$ 
%is a formal orbit of a Lie algebra of formal vector fields (see below).
For every smooth real submanifold $O\subset\R^m$ of codimension $k$ and $p\in O$,
the set $I_p(O)$ of germs at $p$ of smooth real functions in $\R^m$ vanishing on $O$
forms an ideal in the ring of all germs at $p$ of smooth real functions in $\R^m$.
The ideal $I_p(O)$ is
generated by a $k$-tuple $(f_1,\ldots,f_k)$ whose differentials at $p$
\beq\Label{dfs}
df_1(p),\ldots,df_k(p),
\eeq
are linearly independent.
Then the set of all Taylor series
$$
TI_p(O) = j^\infty_p I(O) = \{ T f(x) : f\in I_p(O) \}
\subset \R[[x-p]]
$$
forms an ideal generated by the $k$-tuple
$(Tf_1,\ldots,Tf_k)$ of the Taylor series of the generators $(f_1,\ldots,f_k)$
and the formal differentials 
$$
Tdf_1, \ldots, Tdf_k,
$$
when evaluated at $p$, coincide with \eqref{dfs},
and hence are also linearly independent.
Note that a formal power series in $(x-p)$ can only be evaluated at $x=p$.

Vice versa, for any ideal $I\subset \R[[x-p]]$
generated by $k$ elements $F_1\ldots,F_k\in I$
with linearly independent differentials
$$
dF_1(p),\ldots,dF_k(p),
$$
there is a smooth real submanifold $O$ through $p$ in $\R^m$
of codimension $k$, for which $I=TI_p(O)$.
Indeed, it suffices to choose a smooth function germ $f_j$ at $p$ 
with $Tf_j = F_j$ for each $j=1,\ldots,k$,
and take $O$ defined near $0$ by $f_1=\ldots=f_k=0$.
Here we use the classical 
{\em Borel theorem}
 on surjectivity of the Taylor series map $T$
for smooth functions.
%See \cite{BCP13, CDL20a, CDL20b} for more general
% situations when the
%surjectivity of $T$ (also called the ``Borel map") still holds.

\subsection{Manifold ideals and formal submanifolds}
\Label{mfd-ideals}
The above correspondence between submanifolds $O$
and their formal power series ideals $TI_p(O)=j^\infty_p I_p(O)$
motivates the following abstract definition:

\bd
\Label{m-ideal}
An ideal $I\subset \K[[x]]$, where $x=(x_1,\ldots,x_m)$, is a {\em manifold ideal} if
it is generated by $k$ elements $F_1,\ldots,F_k\in I$
with $F_1(0)=\ldots=F_k(0)$ and
 linearly independent differentials
$$
dF_1(0),\ldots,dF_k(0).
$$
A manifold ideal $I\subset \K[[x]]$ is said to define
a {\em formal submanifold} $O$ of $\K^n$
with its {\em formal ring} given by
$$
\K[[O]]:= \K[[x]]/I.
$$
\ed
One major source of formal submanifolds over $\K=\R$
is provided by $\infty$-jets of real submanifolds $S\subset\R^m$.
If $I_p(S)$ is the ideal of germs at $p\in S$ of smooth functions
in a neighborhood of $p$ in $\R^m$
vanishing on $S$, then
$$
j^\infty_p I(O) = \{ j^\infty_p f : f\in I(O) \} \subset \R[[x-p]]
$$
defines a manifold ideal as can be seen
by taking $\infty$-jets of a local 
system of smooth defining equations of $S$
with independent differentials.
We denote by $j^\infty_p S\subset\R^m$
the corresponding formal submanifold given by 
the ideal $j^\infty_p I(O)$.

By the implicit function theorem for formal power series,
after possible reordering coordinates, a manifold ideal $I\subset \K[[x]]$
with $k$ generators also
admits generators of the form
\beq\Label{spec-gen}
x_{j+m-k} - \phi_j(x'), \quad
\phi_j\in \K[[x']], \quad
x'=(x_1,\ldots,x_{m-k}),
\quad
 j=1,\ldots,k,
\eeq
and hence the formal power series map
$$
A(x') := (x',\phi_1(x'), \ldots, \phi_k(x'))
\colon (\K^{m-k},0) \to (\K^m,0)
$$
satisfies
%\beq
%\Label{par}
$$
\rk dA(0) = k,
\quad
F\circ A = 0 \text{ for all }F\in I.
$$
%\eeq
In particular, for any $F\in I$, its differential $dF(0)$ vanishes on the image $dA(\K^{m-k})\subset \K^{m}$, and hence the number of any generators of $I$ with linearly independent differentials at $0$ is always $k$.
We say that the formal submanifold $O$ defined by $I$
has {\em codimension $k$} and {\em dimension $m-k$}.

\bd
\Label{par}
Given a formal submanifold $O\subset\K^m$ of codimension $k$,
 any formal power series map
 $$
 H
\colon (\K^{m-k},0) \to (\K^m,0),
\quad H \in \K[[x']] \otimes \K^m,
 $$
satisfying \eqref{par}
is called a {\em (formal) parametrization of $O$}.
\ed

Again, as a consequence of the  implicit function theorem for formal power series,
any two formal parametrizations $A_1,A_2$ of $O$ satisfy
$$
A_2=A_1\circ \psi,
\quad
\psi\in (\K^k,0)\to (\K^k,0),
\quad
\rk d\psi(0) =k,
$$
for a suitable formal power series map $\psi$.

\br
\Label{par-smooth-formal}
If $S\subset\R^m$ is a smooth real submanifold
with $p\in S$, and $a\colon (\R^s,0)\to (\R^m,p)$
a germ of a local smooth parametrization of $S$,
then its Taylor series $A=j^\infty_p a$ is a formal 
parametrization of $j^\infty_p S$.
Vice versa,
given any formal parametrization
$A$ of $j^\infty_p S$, there exists
a germ $a\colon (\R^s,0)\to (\R^m,p)$
 of a local smooth parametrization of $S$
 with $A=j^\infty_0 a$.
 Indeed,
 writing $S$ locally as graph $x''=\phi(x')$
 in coordinates $(x',x'')$
 and splitting $A$ as $(A',A'')$ in these coordinates,
 we can construct $a$ as $(a', \phi\circ a')$
 where $a'$ is any smooth germ at $0$ with $A'=j^\infty_0 a'$.
\er

Since we cannot talk about ``points'' in a formal submanifold $O$,  
we only consider $O$ at the level of formal ideals,
i.e.\ $O$ is 
{\em identified with its manifold ideal}.
We write 
$
I=I(O)
$
for the ideal $I\subset \K[[x]]$ defining $O$.

Using the language of ideals, we can transfer
ideas and concepts from Algebraic Geometry
to formal submanifolds.
For example, given two formal submanifolds $X,Y\subset \K^n$,
their inclusion $X\subset Y$ is defined as $I(Y)\supset I(X)$.
In that situation, it will be useful to have 
{\em joint formal parametrizations}:
\bl[joint formal parametrization]
\Label{joint-param}
Let $X\subset Y\subset\K^m$ be formal submanifolds
and $s=\dim Y$, $s'=\dim X$ be their dimensions with 
$s'':=s-s'\ge 0$.
Then there exists a formal parametrization
$$
A(t',t'')\colon (\K^{s'}\times \K^{s''},0) \to (\K^m,0)
$$
of $Y$ such that the restriction
$A(t',0) = A|_{\K^{s'}}\colon (\K^{s'},0) \to (\K^m,0)$
is a formal parametrization of $X$.
\el

\bpf
Let $B\colon (\K^{s},0)\to (\K^m,0)$
be a formal parametrization of $Y$.
Then it is easy to see that the pullback
$B^*(I(X))\subset \K[[t]]$ is a manifold ideal,
hence admitting generators of the form \eqref{spec-gen}.
After a formal change of coordinates $t=(t',t'')\in \K^{s'}\times \K^{s''}$,
these generators become coordinate functions in
the variable $t''$,
leading to the desired property.
\epf

Note that, unlike manifold ideals in $\R[[x]]$ arising as Taylor series map images of ideals defined by smooth real submanifolds, a manifold ideal  $I\subset\C[[x]]$ corresponds to a complex submanifold
through $0$ in $\C^m$ if and only if $I$ has convergent power series as generators. 
%However, formal complex submanifolds defined by manifold ideals in $\C[[x]]$ do arise as {\em complexifications of real formal submanifolds}.

When $\K=\R$, for an ideal $I\subset\R[[x]]$, we also consider its complexification
$$
\C I := \left\{
f + ig: f,g\in I
 \right\} \subset \C[[x]],
$$
which is the ideal in $\C[[x]]$ generated by $I$,
that is invariant under conjugation.
If $I$ is a manifold ideal in $\R[[x]]$ with generators $F_1,\ldots,F_k$,
then $\C I$ is a manifold ideal in $\C[[x]]$ with the same generators.
Vice versa, if $\2I\subset\C[[x]]$ is a manifold ideal invariant under the conjugation, it equals the complexification $\C I$ 
of a real manifold ideal $I\subset \R[[x]]$.

For background material and further results on formal submanifolds, the reader is referred to \cite{BER00, BRZ01, BMR, BER03, MMZ03, DJL12}.

\subsection{Formal tangent spaces and vector fields}
\Label{form-cr}
Since we can evaluate formal power series at the origin,
we can define the tangent space $T_0O$ for every formal submanifold
as follows:

\bd[tangent space of a formal submanifold]
\Label{t-sp}
For a formal submanifold $O\subset \K^n$, define its 
{\em tangent space at $0$}
as
$$
T_0O := \{ \xi\in T_0\K^n : df(\xi)=0 \text{ for all } f\in I(O)\}
\subset T_0\K^n \cong \K^n
.
$$
\ed
Note that $0$ is the only point where we can evaluate
formal power series, hence only the tangent space at $0$
is defined for a formal submanifold.
Another standard notion is that of formal 
{\em tangent vector fields}
to a formal submanifold:

\bd\Label{formal-tangent}
If $O\subset \K^n$ is a formal submanifold,
a formal vector field $L\in \K[[x]]\otimes \K^n$
is said to be {\em tangent to $O$} whenever
$$
L(I(O))\subset I(O).
$$
We denote by $D_O$ the $\K[[x]]$-module of 
{\em all formal vector fields
tangent to $O$}, which is obviously a Lie algebra
with respect to the Lie bracket.
\ed
With generators of $I$ chosen as in \eqref{spec-gen},
the space $D_O$ is spanned, modulo $I(O)$, by
\beq
\Label{do}
%\textstyle
D_O = \span_{\K[[x]]}
\big\{
\d_{x_s} + \sum_j \frac{\d \phi_j}{\d x_s} \d_{x_{j+m-k}} : 
s=1,\ldots, m-k
\big\}
+ I(O)\cdot \K[[x]]\otimes \K^m,
\eeq
where as before $\K[[x]]\otimes \K^m$ is identified
with the space of all formal vector fields in $\K^m$.
In particular, evaluating at $0$, it is easy to see that
\beq
\Label{d-t}
D_O(0) :=  \{ L_0 : L\in D_O \} = T_0 O.
\eeq
%Note that unless $O=\K^m$, we cannot 
%the identity \eqref{do} does not hold 
%without 

More generally, to span $D_O$ over the ring $\K[[x]]$
by a collection of formal vector fields,
it suffices to span the tangent space at $0$:

\bl
\Label{span}
Let $O\subset\K^m$ be a formal submanifold and
$L^1,\ldots,L^l\in D_O$ be such that
$$
T_0 O = {\sf span}_\K\{L^1_0,\ldots,L^l_0\}.
$$
Then
$$
D_O =  {\sf span}_{\K[[x]]}\{L^1,\ldots,L^l\} 
+ I(O)\cdot \K[[x]]\otimes \K^m.
$$
\el

\bpf
Reordering coordinates if necessary, we may assume that
 $I(O)$ has generators 
as in \eqref{spec-gen}, so that \eqref{do} holds, i.e.
\beq
\Label{x-span}
D_O = {\sf span}_ {\K[[x]]}\{X^s: s=1,\ldots, m-k \}
+ I(O)\cdot \K[[x]]\otimes \K^m,
\quad
X^s = \d_{x_s} + \sum_j \frac{\d \phi_j}{\d x_s} \d_{x_{j+m-k}}.
\eeq
In particular, we can write 
$$
L^j = \sum_s A^j_s X^s
\mod I(O),
$$
where $A^j_s\in \K[[x]]$.
Since $L^j_0$ span $T_0O$,
the matrix $(A^j_s(0))$ evaluated at $0$ has rank $m-k$.
Reordering $L^j$ if necessary, we may assume the first
$(m-k)\times (m-k)$ minor of $(A^j_s(0))$ is nonzero.
Then the same minor of $(A^j_s)$ is invertible in $\K[[x]]$.
Taking the inverse, we conclude that each $X^j$ is contained in the span of $L^j$ modulo $I(O)$, which proves the desired conclusion
in view of \eqref{x-span}.
\epf

%
%\bc
%Let $A\colon (\K^s,0)\to (\R^m,0)$ be 
%any formal parametrization of a formal submanifold $O\subset\K^m$.
%Then 
%$$
%D_O = {\sf span}_{\K[[x]]} \{A_* \d_{t_1}, \ldots, A_* \d_{t_s} \} 
%+ I(O)\cdot \K[[x]]\otimes \K^m
%$$
%\ec

For a pair of formal submanifolds $O\subset V$,
we can approximate vector fields in $D_O$
by those in $D_V$ in the following sense:

\bl
\Label{ov-ap}
Let $O\subset V\subset\K[[x]]$ be formal submanifolds.
Then for every $L\in D_O$ there exists $\2L\in D_V$
with
\beq
\Label{2l-ap}
\2L = L \mod I(O).
\eeq
In particular, it follows that $\2L\in D_O\cap D_V$.
\el

\bpf
%By Lemma~\ref{joint-param},
%there exists a joing parametrization
%$A(t',t'')\colon (\K^{s'}\times \K^{s''},0) \to (\K^m,0)$
%of $V$ such that the restriction
%$A(t',0) = A|_{\K^{s'}}\colon (\K^{s'},0) \to (\K^m,0)$
%is a formal parametrization of $O$.
Reordering coordinates $x=(x',x'')\in \K^{m-k}\times \K^{k}$ 
if necessary, we may assume
that $I(V)$ has generators 
as in \eqref{spec-gen}, so that \eqref{do} holds for $D_V$, i.e.
%\beq
%\Label{x-span1}
$$
D_V = {\sf span}_ {\K[[x]]}\{X^s: s=1,\ldots, m-k \}
+ I(V)\cdot \K[[x]]\otimes \K^m,
\quad
X^s = \d_{x_s} + \sum_j \frac{\d \phi_j}{\d x_s} \d_{x_{j+m-k}}.
$$
%\eeq
Writing 
$
L = \sum_{j=1}^m a_j\d_{x_j}
$, $a_j \in \K[[x]]$,
we set
$
\2L := \sum_{s=1}^{m-k} a_s X^s \in D_V
$.
Since $O\subset V$, hence $I(V)\subset I(O)$,
we have $L(I(V))\subset I(O)$.
Writing this condition for the generators
$ x_{j+m-k} - \phi_j(x')$ of $(V)$, we obtain
$$
L( x_{j+m-k} - \phi_j(x'))  
= a_{j+m-k} - \sum_s a_s \frac{\d \phi_j}{\d x_s}
\in I(O)
$$
and hence
$$
L - \2L = 
\sum_{s=1}^{m-k} a_s\d_{x_s} + \sum_{j=1}^k a_{j+m-k}\d_{x_{j+m-k}}  - 
\sum_{s=1}^{m-k} a_s (\d_{x_s} + \sum_{j=1}^k \frac{\d \phi_j}{\d x_s} \d_{x_{j+m-k}})
\in I(O) \cdot \K[[x]]\otimes \K^m
$$
proving \eqref{2l-ap}.
The inclusion $\2L\in D_O$
follows immediately from \eqref{2l-ap}
and the definition of $D_O$ (Definition~\ref{formal-tangent}).
\epf

\subsection{Formal CR submanifolds}

It is customary to write a real formal power series $f$
in $\cn\cong \R^{2n}$ as series in $z$ and $\bar z$:
\beq\Label{fzz}
f(z,\bar z)=\sum_{\a\b} c_{\a\b} z^\a \bar z^\b
\in \C[[z,\bar z]],
\quad c_{\a\b} = \1{c_{\b\a}},
\quad z=(z_1,\ldots,z_n),
\eeq
which we identify with the complexified series
\beq\Label{real-formal}
\sum_{\a\b} c_{\a\b} z^\a \zeta^\b
\in \C[[z,\zeta]],
\quad c_{\a\b} = \1{c_{\b\a}},
\quad z=(z_1,\ldots,z_n),
\; \zeta=(\zeta_1,\ldots,\zeta_n).
\eeq
With that identification,
the complex conjugation $\s$ for series \eqref{fzz}
is
given by
\beq\Label{conjugation}
\s
\left(
{\sum_{\a\b} c_{\a\b} z^\a \zeta^\b}
\right)
=
\sum_{\a\b} \1{c_{\a\b}} \zeta^\a z^\b.
\eeq
%Similarly the conjugation of formal complex vector fields
%is given by
%\beq\Label{conjugation-vf}
%\s
%\left(
%\sum c_{\a\b j} z^\a \zeta^\b \d_{z_j}
%+\sum
%d_{\a\b j} z^\a \zeta^\b \d_{\zeta_j}
%\right)
%=
%\sum \1{c_{\a\b j}} \zeta^\a z^\b \d_{\zeta_j}
%+\sum
%\1{d_{\a\b j}} \zeta^\a z^\b \d_{z_j}
%.
%\eeq
%It is straightforward to see that manifold ideals in $\R[[x]]$, 
%$x\in \R^{2n}\cong\cn$, are thereby identified
%with manifold ideals in $\C[[z,\zeta]]$ that are invariant under conjugation, where $z,\zeta \in \C^n$.
Then a real
formal submanifold $O\subset \C^n$
has its ideal 
$
I=I(O)
$
in the ring  
$$
\{f\in  \C[[z,\zeta]] : \sigma f = f\} \cong 
 \R[[x,y]], 
 \quad x,y\in \R^n,
$$
 whose complexification 
$$
\C I = \{ f + ig : f, g\in I\}
$$
is a $\sigma$-invariant ideal in $\C[[z,\zeta]]$.
Conversely, a $\sigma$-invariant ideal in $\C[[z,\zeta]]$
is a complexification of a real ideal in
the ring $\{f\in  \C[[z,\zeta]] : \sigma f = f\}$,
and a $\sigma$-invariant manifold ideal in $\C[[z,\zeta]]$
is a complexification of a real manifold ideal in $\{f\in  \C[[z,\zeta]] : \sigma f = f\}$.
We call 
$$\C I =\C I(O)\subset \C[[z,\zeta]]$$
the {\em complexified ideal} of $O$.

For a complex formal submanifold $O\subset \cn\times \cn$
given by a manifold ideal $I\subset \C[[z,\zeta]]$,
%and its complexified ideal $I=\C I(O)\subset\C[[z,\zeta]]$,
define the spaces of $(1,0)$ and $(0,1)$ formal tangent vector fields by
\beq
\Label{d10}
%D_I = \{ L\in \6D : L(I)\subset I \}
%\quad
\6D^{10}_O = 
\left\{ L = \sum a_j(z,\zeta) \d_{z_j}
%\in \C[[z,\zeta]]\otimes \C^n
 :  L(I)\subset I 
\right\},
\quad
\6D^{01}_O = 
\left\{ L = \sum b_j(z,\zeta) \d_{\zeta_j} :  L(I)\subset I 
\right\},
\eeq
where $a_j(z,\zeta), b_j(z,\zeta)\in \C[[z,\zeta]]$,
and the space of complex formal tangent vector fields by
$$
\C D_O = 
\left\{
L = \sum (a_j(z,\zeta) \d_{z_j} + b_j(z,\zeta) \d_{\zeta_j}) : L(I)\subset I
\right\},
$$
%Note that $\C D_O$ can be identified with the complexification
while the space $D_O$ of real formal tangent vector fields is
$$
D_O = \{ L \in \C D_O : \s L = L\}.
$$

If $O$ is a {\em real formal submanifold}, i.e.\ defined 
by a real submanifold ideal $I(O)\subset\{f\in  \C[[z,\zeta]] : \sigma f = f\}$, 
the corresponding complexification $\C I(O)$  in $\C[[z,\zeta]]$
is generated by the $\sigma$-invariant formal power series \eqref{real-formal} and one has
\beq\Label{d-conj}
\6D^{01}_O = \s{\6D^{10}_O}=
\left\{
\sum \bar a_j(\zeta, z)\d_{\zeta_j}
:
\sum a_j(z,\zeta) \d_{z_j}
\in \6D^{10}_O
\right\},
\eeq
where we keep the notation $\s$
for complex formal vector fields
given by
$$
\s
\left(
\sum (a_j \d_{z_j} + b_j \d_{\zeta_j}) 
\right)
=
\sum (\s a_j \d_{\zeta_j} + \s b_j \d_{z_j}),
$$
or, equivalently, by
$$
\s(L) f= \s (L (\s f)),
\quad f\in \C[[z,\zeta]].
$$
%
% by complexification of
%the standard conjugation
%on formal complex vector fields in $\cn$.

It is easy to see that the evaluations 
$$
\6D^{10}_O(0) = 
\left\{  
\sum a_j(0) \d_{z_j} : \sum a_j(z,\zeta) \d_{z_j} \in \6D^{10}_O
\right \},
$$
$$
\6D^{01}_O(0) = 
\left\{  
\sum a_j(0) \d_{z_j} : \sum a_j(z,\zeta) \d_{\zeta_j} \in \6D^{01}_O
\right \},
$$
%
%$\6D^{10}_O(0)$
%and $\6D^{01}_O(0)$ 
are contained in the corresponding 
$(1,0)$ and $(0,1)$ tangent spaces at $0$:
$$
H^{10}_0 O = H^{10}_0 \cn \cap \C T_0O, 
\quad
H^{01}_0 O = H^{01}_0 \cn \cap \C T_0O,
$$
i.e.\
$$
\6D^{10}_O(0) \subset H^{10}_0 O,
\quad
\6D^{01}_O(0) \subset H^{01}_0 O.
$$

The opposite inclusions do not hold in general
and, in the case $O$ is real-analytic,
are equivalent to $O$ being CR in a neighborhood of $0$.
This motivates the following definition (cf. \cite[p. 105]{DJL12}:

\bd\Label{CR-formal}
A real formal submanifold $O\subset \cn$ 
%defined by an ideal in $\C[[z,\zeta]]$
is said to be CR if
$$
 H^{10}_0 O \subset \6D^{10}_O(0) ,
$$
or, equivalently, if
every vector in $H^{10}_0 O$
extends to a formal vector field in $\6D^{10}_O$.
\ed
In view of \eqref{d-conj},
Definition~\ref{CR-formal} can be restated for $(0,1)$ vector fields,
i.e.\
a real formal submanifold $O\subset \C^{n}$ is CR
if and only if
$\6D^{10}_O(0) = H^{10}_0 O$ or, equivalently,
 $\6D^{01}_O(0) = H^{01}_0 O$.
We quote the following structure result for formal CR submanifolds
 from Della Salla-Juhlin-Lamel~\cite[Lemma~1]{DJL12}:
\bl\Label{djl}
Let $O\subset \cn$ be a formal CR submanifold of codimension $d$
with its complexified ideal $I=\C I(O)\subset \C[[z,\zeta]]$.
Then there exists integers $d_1,d_2,d_3\ge0$ with
$
n=d_1+d_2+d_3$,
$
d=2d_1+d_3,
$
such that after renumbering coordinates $z$ and $\zeta$
if necessary, we can write
$$
z=(z^1,z^2,z^3),\; 
\zeta=(\zeta^1,\zeta^2,\zeta^3) \in
\C^{d_1}\times \C^{d_2}\times \C^{d_3},
$$
and choose a set of generators of $I$ given by
the components of
$$
z^1-\phi(z^2,z^3),
\quad
\zeta^1 - \s\phi(\zeta^2,\zeta^3),
\quad
\zeta^2 - R(z^2,z^3,\zeta^3)
$$
for some vector-valued formal power series
$\phi\in\C[[z^2,z^3]]^{d_1}$,
$R\in\C[[z^2,z^3,\zeta^3]]$.
\el

It is well-known that real-analytic CR submanifolds $O$
are generic in their local 
{\em intrinsic complexifications}
defined as minimal complex submanifolds (locally) containing $O$,
see e.g.\ \cite[Chapter I]{BERbook}.
We shall need the following formal analogue:

\bc
\Label{cr-cx}
Let
$O\subset \cn$ be a formal CR submanifold.
Then the maximal formal holomorphic ideal $I'\subset \C[[z]]$
 contained in the complexification 
 $\C I(O) = I(O) + i I(O)\subset \C[[z,\zeta]]$, i.e.\ 
 $$
 I':=\C I(O)\cap \C[[z]],
 $$
is a complex manifold ideal,
(where $\C[[z]]$ is identified with the subring of all power series
in $\C[[z,\zeta]]$
independent of $\zeta$).
\ec

\bpf
Using the notation and generators of $I(O)$ 
in Lemma~\ref{djl},
it follows that $I'$ is generated over $\C[[z]]$ by the components 
of $z^1-\phi(z^2,z^3)$, hence is a manifold ideal as claimed.
\epf

\bd\Label{i-c}
The {\em intrinsic complexification} of
a formal CR submanifold $O\subset \cn$
is the complex
 formal submanifold $V\subset\cn$
defined by the ideal
$$
I'=
\C
I(O)\cap \C[[z]].
$$
\ed
To any formal complex submanifold $V\subset\cn$
given by an ideal $I(V)\subset \C[[z]]$,
there is also the associated real formal submanifold 
$V^\R\subset \cn$ 
given by the real formal ideal 
$I(V^\R) \subset \{ f \in \C[[z,\zeta]]: f=\s f \}$
whose complexification $\C I(V^\R)$ is
generated by
$$I(V) \oplus \s I(V) = \{ f + \s g : f, g\in I(V) \} \subset \C[[z,\zeta]],$$
and using these generators, it follows that
$$
\C D_{V^\R} = D_{V^\R}^{10} \oplus D_{V^\R}^{01}
\subset \C[[z,\zeta]]\otimes \cn.
$$
%and 
%$$
%D_{V^\R}^{10}(0) = H \oplus D_{V^\R}^{01}
%\subset \C[[z,\zeta]]\otimes \cn
%$$
%in view of \eqref{d-t}.

In the sequel, we slightly abuse the notation
writing $V$ for both complex formal submanifold
and its associated real formal submanifold $V^\R$.
We have the following consequence of Lemma~\ref{djl}:

\bc
\Label{i-c-prop}
If $V$ is the intrinsic complexification of
a formal CR submanifold $O\subset\cn$,
it follows that $O$ is {\em generic in} $V$ in the sense that
\beq
\Label{generic}
T_0 O + J T_0 O = T_0 V,
%\quad
%D_O\cap D_V + J({D_O}\cap D_V) = D_{V} \mod I(V),
\eeq
where $J$ is the standard complex structure acting on 
the standard basis vectors by 
$$
J\d_{z_j}=i\d_{z_j},
\quad
J\d_{\zeta_j}=-i\d_{\zeta_j},
\quad
j=1,\ldots,n
.
$$
\ec

%
%
%$$
%L = (L +iJL) + (L-iJL)
%$$

As an application,
we obtain a convenient
expression of tangent vector fields to 
the complexification $V$ in terms
of vector fields tangent to both $O$ and $V$:

\bc
\Label{dv-split}
If $O\subset V\subset\cn$ are formal submanifolds
satisfying \eqref{generic},
it follows that
\beq
\Label{dv-j}
D_V = (D_O\cap D_V) + J(D_O\cap D_V)
\mod I(V).
\eeq
In particular, \eqref{dv-j}
holds when $O$ is a formal CR submanifold and
 $V$ is the intrinsic complexification of $O$.
\ec

\bpf
In view of \eqref{d-t},
there exists a collection of vector fields
$L^1,\ldots, L^s \in D_O$ such that 
$$
T_0 O = \span_{\C} \{L^1_0,\ldots,L^s_0\}.
$$
By Lemma~\ref{2l-ap}, there exists approximations
$\2L^1,\ldots, \2L^s \in D_O\cap D_V$
with $\2L^j=L^j \mod I(O)$.
In particular, the values $\2L^1_0,\ldots, \2L^s_0$
also span $T_0 O$.
Then by Corollary~\ref{i-c-prop},
$$
T_0 V = \span_{\C} \{L^1_0,\ldots,L^s_0, JL^1_0, \ldots, JL^s_0\}.
$$ 
Finally, since $V$ is a complex formal submanifold,
each $JL^j$ belongs to $D_V$,
and the desired conclusion follow from 
Lemma~\ref{span} applied
to $V$.
\epf

\subsection{Formal orbits}

We make use of
a result of Baouendi-Ebenfelt-Rothschild \cite[Proposition~5.1]{BER03}
(an adaptation to the formal case
 of an unpublished proof of Nagano's theorem by the author)
to obtain a formal variant of Nagano's theorem:
%While their original result is stated over $\C$,
%we need the following analogous statement over $\R$,
%which can be proved along the same lines:

\bp[existence of a formal orbit]
\Label{orb-exist}
Let $\5g$ be Lie algrebra of formal vector fields over $\K$ in
$n$ indeterminates $x=(x_1,\ldots,x_n)$.
%which is also a $\R[[x]]$-submodule.
Then there exists an unique
formal submanifold $O\subset \K^n$
with $\dim O =\dim\5g(0)$
such that any vector field $L\in\5g$ is tangent to $O$,
i.e.\ $L (I(O))\subset I(O)$.
Moreover if $O'\subset\K^n$ is 
any other formal manifold such that
all vector fields in $\5g$ are tangent to $O'$,
then $O\subset O'$.
\ep

Here $\5g(0)\subset T_0\R^n$ denotes the subspace
consisting of all values $L_0$ at $0$ of $L\in \5g$.
We call the formal submanifold $O$ in Proposition~\ref{orb-exist}
the {\em orbit} of $\5g$.

\br\Label{tangent-eq}
Since every value $L_0$ is contained in $T_0O$,
the equality $\dim O =\dim\5g(0)$ holds
if and only if $T_0O = \5g(0)$.
\er

\bpf[Proof of Proposition~\ref{orb-exist}]
In case $\K=\C$, the $\C[[x]]$-submodule $\2{\5g}\subset \C[[x]]\otimes\C^n$
generated by $\5g$ 
is a Lie algebra satisfying the assumptions of \cite[Proposition~5.1]{BER03}, 
which yields the existence and uniqueness of a $\2{\5g}$-orbit $O$.
Then the same $O$ is also a $\5g$-orbit, proving the statement for $\K=\C$.
The uniqueness follows since any $\5g$-orbit is also a $\2{\5g}$-orbit.
Finally if $O'$ is as in the proposition, then all vector fields
in $\2{\5g}$ are tangent to $O'$ and $O\subset O'$
by the second statement of \cite[Proposition~5.1]{BER03}.

In case $\K=\R$, consider
the complexification 
$$
\{X+iY: X,Y\in \5g\} \subset\C[[x]]\otimes \C^n,
$$
which is a Lie subalgebra over $\C$
that is invariant under conjugations.
By the statement for $\K=\C$,
there exists a $\2{\5g}$-orbit $\2O\subset \C^n$,
which is a formal submanifold over $\C$.
Since $\2{\5g}$ is invariant under conjugations
and $L\1f = \1{ \1L f}$ for $L\in \2{\5g}$, $f\in I(\2O)$,
it follows that $\1{I(\2O)}$ is also a submanifold ideal
defining a formal submanifold of the same dimension as $\2O$, 
to which $\2{\5g}$ is also tangent. 
Hence by the uniqueness of the $\2{\5g}$-orbit, 
$\1{I(\2O)} = I(\2O)$, which implies that $I(\2O)$
is the complexification of a real manifold ideal $I=I(O)\subset\R[[x]]$
defining a formal manifold $O\subset\R^n$.
Then any $L\in \5g\subset \2{\5g}$
satisfies $L(I(\2O))\subset I(\2O)$
and hence $L(I(O))\subset I(O)$
since $L$ is real and $I(O)=I(\2O)\cap \R[[x]]$.
If $O'$ is any other $\5g$-orbit,
its complexification $\2O'$ is a $\2{\5g}$-orbit,
hence $\2O'=\2O$ by the uniqueness for $\K=\C$,
and then intersecting with $\R[[x]]$ yields $O'=O$.
Finally, if $O'$ is as in the proposition, 
all vector fields in $\2g$ are tangent to
its complexification $\2O'$,
hence by the case $\K=\C$ of the proposition,
$\2O\subset \2O'$, i.e.\ $I(\2O)\supset I(\2O')$.
Intersecting with $\R[[x]]$, we obtain
$I(O)\supset I(O')$, i.e.\ $O\subset O'$ as claimed.
%Finally, the uniqueness 
\epf

Taking generators of $O$ of the form \eqref{spec-gen}
and corresponding basis of $D_O$ given by \eqref{do},
and using the assumption that $\5g$ is a $\K[[x]]$-submodule,
one can prove the following property of the orbit:
\bl\Label{orb-tangent}
Let $\5g$ be as in Proposition~\ref{orb-exist}
that is also a $\K[[x]]$-module.
Then every formal vector field
that is tangent to the orbit
$O$ of $\5g$ is necessarily contained in $\5g$ modulo $O$,
in fact
$$
D_O = \5g \mod I(O).
$$
\el

This ends the Appendix section.

\section*{Acknowledgements}
%\bigskip
%{\bf Acknowledgements.}
The author is deeply grateful 
to Professor Joe Kohn,
to whose memory this work is dedicated,
and   
without whose constant interest and encouragement
this work might not have been possible.
The author would also like to thank
G.~Bharali,
D.W.~Catlin, J.P.~D'Angelo, S.~Fu,
A.C.~Nicoara and E.J.~Straube
for inspiring discussions on the subject of this work,
 B.~Lamel and J.~Lebl
for sharing Example~\ref{l-l}
and X.~Huang for his interest and patience in 
verifying many details of the proofs.

\end{document}